*FLORENTIN SMARANDACHE*

# *SEQUENCES OF NUMBERS INVOLVED IN UNSOLVED PROBLEMS*

1
141
1

1
1 8 1
1 8 40 8 1
1 8 1
1

1
1 12 1
1 12 108 12 1
1 12 108 540 108 12 1
1 12 108 12 1
1 12 1
1

1
1 16 1
1 16 208 16 1
1 16 208 1872 208 16 1
1 16 208 1872 9360 1872 208 16 1
1 16 208 1872 208 16 1
1 16 208 16 1
1 16 1
1

# Introduction

Over 300 sequences and many unsolved problems and conjectures related to them are presented herein.

These notions, definitions, unsolved problems, questions, theorems corollaries, formulae, conjectures, examples, mathematical criteria, etc. ( on integer sequences, numbers, quotients, residues, exponents, sieves, pseudo-primes/squares/cubes/factorials, almost primes, mobile periodicals, functions, tables, prime/square/factorial bases, generalized factorials, generalized palindromes, etc. ) have been extracted from the Archives of American Mathematics (University of Texas at Austin) and Arizona State University (Tempe): "The Florentin Smarandache papers" special collections, University of Craiova Library, and Arhivele Statului (Filiala Vâlcea, Romania).
It is based on the old article "Properties of Numbers" (1975), updated many times.

Special thanks to C. Dumitrescu & V. Seleacu from the University of Craiova (see their edited book "Some Notions and Questions in Number Theory", Erhus Univ. Press, Glendale, 1994), M. Perez, J. Castillo, M. Bencze, L. Tutescu, E, Burton who helped in collecting and editing this material.

The Author

## Sequences of Numbers Involved in Unsolved Problems

Here it is a long list of sequences, functions, unsolved problems, conjectures, theorems, relationships, operations, etc. Some of them are inter-connected.

1) Consecutive Sequence:

1,12,123,1234,12345,123456,1234567,12345678,123456789,12345678910,
1234567891011,123456789101112,12345678910111213,...

How many primes are there among these numbers?

In a general form, the Consecutive Sequence is considered in an arbitrary numeration base B.

2) Circular Sequence:

1,12,21,123,231,312,1234,2341,3412,4123,
12345,23451,34512,45123,51234,

1 2 3 4 5

123456,234561,345612,456123,561234,612345,1234567,2345671,3456712,...

6 7

3) Symmetric Sequence:

1,11,121,1221,12321,123321,1234321,12344321,123454321,1234554321,

12345654321,123456654321,1234567654321,12345677654321,123456787654321,
1234567887654321,12345678987654321,123456789987654321,

12345678910987654321,1234567891010987654321,123456789101110987654321,
12345678910111110987654321, ...

How many primes are there among these numbers?

In a general form, the Symmetric Sequence is considered in an arbitrary numeration base B.

4) Deconstructive Sequence:

$\underbrace{1,23,456,789}\ \underbrace{1,23456,789}\ \underbrace{123,456789}\ \underbrace{1,23456789,}\ \underbrace{123456789,}\ \underbrace{123456789}\ \underbrace{1,...}$

5) Mirror Sequence:

1,212,32123,4321234,543212345,65432123456,7654321234567,876543212345678,
98765432123456789,109876543212345678910,1110987654321234567891011,...

Question: How many of them are primes?

6) Permutation Sequence:

12,1342,135642,13578642,13579108642,135791112108642,1357911131412108642,
135791113151614121086421,35791113151718161412108642,
13579111315171920181614121086642,...

Question: Is there any perfect power among these numbers?
(Their last digit should be:
either 2 for exponents of the form 4k+1,
either 8 for exponents of the form 4k+3, where $k \geq 0$ .)

We conjecture: no!

7) Generalized Permutation Sequence:

If g(n), as a function, gives the number of digits of a(n), and F is a permutation of g(n) elements, then:

$$a(n) = \overline{F(1)F(2)...F(g(n))}$$

8) Mobile Periodicals (I):

...000000000000000000000000000000000100000000000000000000000000000000000000...
...000000000000000000000000000000001110000000000000000000000000000000000000...
...000000000000000000000000000000011011000000000000000000000000000000000000...
...000000000000000000000000000000001110000000000000000000000000000000000000...
...000000000000000000000000000000000100000000000000000000000000000000000000...
...000000000000000000000000000000001110000000000000000000000000000000000000...
...000000000000000000000000000000011011000000000000000000000000000000000000...
...000000000000000000000000000000110001100000000000000000000000000000000000...
...000000000000000000000000000000011011000000000000000000000000000000000000...
...000000000000000000000000000000001110000000000000000000000000000000000000...
...000000000000000000000000000000000100000000000000000000000000000000000000...
...000000000000000000000000000000001110000000000000000000000000000000000000...
...000000000000000000000000000000011011000000000000000000000000000000000000...
...000000000000000000000000000000110001100000000000000000000000000000000000...
...000000000000000000000000000001100000110000000000000000000000000000000000...
...000000000000000000000000000000110001100000000000000000000000000000000000...
...000000000000000000000000000000011011000000000000000000000000000000000000...
...000000000000000000000000000000001110000000000000000000000000000000000000...
...000000000000000000000000000000000100000000000000000000000000000000000000...
...000000000000000000000000000000001110000000000000000000000000000000000000...
...000000000000000000000000000000011011000000000000000000000000000000000000...
...000000000000000000000000000000110001100000000000000000000000000000000000...
...000000000000000000000000000001100000110000000000000000000000000000000000...
...000000000000000000000000000011000000011000000000000000000000000000000000...
...000000000000000000000000000001100000110000000000000000000000000000000000...
...000000000000000000000000000000110001100000000000000000000000000000000000...
...000000000000000000000000000000011011000000000000000000000000000000000000...
...000000000000000000000000000000001110000000000000000000000000000000000000...
...000000000000000000000000000000000100000000000000000000000000000000000000...
...000000000000000000000000000000001110000000000000000000000000000000000000...
...000000000000000000000000000000011011000000000000000000000000000000000000...
...000000000000000000000000000000110001100000000000000000000000000000000000...
...000000000000000000000000000001100000110000000000000000000000000000000000...
...000000000000000000000000000011000000011000000000000000000000000000000000...
...000000000000000000000000000110000000001100000000000000000000000000000000...
...000000000000000000000000000011000000011000000000000000000000000000000000...
...000000000000000000000000000001100000110000000000000000000000000000000000...
...000000000000000000000000000000110001100000000000000000000000000000000000...
...000000000000000000000000000000011011000000000000000000000000000000000000...
...000000000000000000000000000000001110000000000000000000000000000000000000...
...000000000000000000000000000000000100000000000000000000000000000000000000...
...000000000000000000000000000000001110000000000000000000000000000000000000...
...000000000000000000000000000000011011000000000000000000000000000000000000...

...0000000000000000000000000000110001100000000000000000000000000000000...
...000000000000000000000000000110000011000000000000000000000000000000...
...00000000000000000000000000110000000110000000000000000000000000000...
...0000000000000000000000000110000000001100000000000000000000000000...
...000000000000000000000000110000000000011000000000000000000000000000...
..........................................................................

This sequence has the form

$$1,\underbrace{111,11011,111,}_{5}\underbrace{1,111,11011,1100011,11011,}_{7}\underbrace{111,1,111,11011,1100011,110000011}_{9},...$$

9) Mobile Periodicals (II):

...000000000000000000000000000000010000000000000000000000000000000000...
...000000000000000000000000000000111000000000000000000000000000000000...
...000000000000000000000000000001121100000000000000000000000000000000...
...000000000000000000000000000000111000000000000000000000000000000000...
...000000000000000000000000000000010000000000000000000000000000000000...
...000000000000000000000000000000111000000000000000000000000000000000...
...000000000000000000000000000001121100000000000000000000000000000000...
...000000000000000000000000000011232110000000000000000000000000000000...
...000000000000000000000000000001121100000000000000000000000000000000...
...000000000000000000000000000000111000000000000000000000000000000000...
...000000000000000000000000000000010000000000000000000000000000000000...
...000000000000000000000000000000111000000000000000000000000000000000...
...000000000000000000000000000001121100000000000000000000000000000000...
...000000000000000000000000000011232110000000000000000000000000000000...
...000000000000000000000000000112343211000000000000000000000000000000...
...000000000000000000000000000011232110000000000000000000000000000000...
...000000000000000000000000000001121100000000000000000000000000000000...
...000000000000000000000000000000111000000000000000000000000000000000...
...000000000000000000000000000000010000000000000000000000000000000000...
...000000000000000000000000000000111000000000000000000000000000000000...
...000000000000000000000000000001121100000000000000000000000000000000...
...000000000000000000000000000011232110000000000000000000000000000000...
...000000000000000000000000000112343211000000000000000000000000000000...
...000000000000000000000000001123454321100000000000000000000000000000...
...000000000000000000000000000112343211000000000000000000000000000000...
...000000000000000000000000000011232110000000000000000000000000000000...
...000000000000000000000000000001121100000000000000000000000000000000...
...000000000000000000000000000000111000000000000000000000000000000000...
...000000000000000000000000000000010000000000000000000000000000000000...
...000000000000000000000000000000111000000000000000000000000000000000...
...000000000000000000000000000001121100000000000000000000000000000000...
...000000000000000000000000000011232110000000000000000000000000000000...
...000000000000000000000000000112343211000000000000000000000000000000...
...000000000000000000000000001123454321100000000000000000000000000000...
...000000000000000000000000011234565432110000000000000000000000000000...
...000000000000000000000000001123454321100000000000000000000000000000...
...000000000000000000000000000112343211000000000000000000000000000000...
...000000000000000000000000000011232110000000000000000000000000000000...
...000000000000000000000000000001121100000000000000000000000000000000...
...000000000000000000000000000000111000000000000000000000000000000000...
...000000000000000000000000000000010000000000000000000000000000000000...
...000000000000000000000000000000111000000000000000000000000000000000...
...000000000000000000000000000001121100000000000000000000000000000000...

...0000000000000000000000000000112321100000000000000000000000000000000...
...000000000000000000000000000112343211000000000000000000000000000000...
...0000000000000000000000000011234543211000000000000000000000000000000...
...00000000000000000000000001123456543211000000000000000000000000000...
...000000000000000000000000112345676543211000000000000000000000000000...
................................................................

This sequence has the form

1,111,11211,111,1,111,11211,1123211,11211,111,1,111,11211,1123211,112343211,...

5 7 9

10) Infinite Numbers (I):

...11111111111111111111111111111110111111111111111111111111111111111...

...111111111111111111111111111111000111111111111111111111111111111111...

...11111111111111111111111111111001001111111111111111111111111111111111...

...111111111111111111111111111111000111111111111111111111111111111111...

...11111111111111111111111111111110111111111111111111111111111111111...

...111111111111111111111111111111000111111111111111111111111111111111...

...11111111111111111111111111111001001111111111111111111111111111111111...

...111111111111111111111111111100111001111111111111111111111111111111...

...11111111111111111111111111111001001111111111111111111111111111111111...

...111111111111111111111111111111000111111111111111111111111111111111...

...11111111111111111111111111111110111111111111111111111111111111111...

...111111111111111111111111111111000111111111111111111111111111111111...

...11111111111111111111111111111001001111111111111111111111111111111111...

...111111111111111111111111111100111001111111111111111111111111111111...

...111111111111111111111111111001111100111111111111111111111111111111...

...111111111111111111111111111100111001111111111111111111111111111111...

...11111111111111111111111111111001001111111111111111111111111111111111...

...111111111111111111111111111111000111111111111111111111111111111111...

...11111111111111111111111111111110111111111111111111111111111111111...

...111111111111111111111111111111000111111111111111111111111111111111...

...11111111111111111111111111111001001111111111111111111111111111111111...

...111111111111111111111111111100111001111111111111111111111111111111...

...111111111111111111111111111001111100111111111111111111111111111111...

...111111111111111111111111110011111110011111111111111111111111111111...

...111111111111111111111111111001111100111111111111111111111111111111...

...111111111111111111111111111100111001111111111111111111111111111111...

...11111111111111111111111111111001001111111111111111111111111111111111...

...111111111111111111111111111111000111111111111111111111111111111111...

...11111111111111111111111111111110111111111111111111111111111111111...

...111111111111111111111111111111000111111111111111111111111111111111...

...11111111111111111111111111111001001111111111111111111111111111111111...

...111111111111111111111111111100111001111111111111111111111111111111...

...111111111111111111111111111001111100111111111111111111111111111111...

...111111111111111111111111110011111110011111111111111111111111111111...

...111111111111111111111111100111111111001111111111111111111111111111...

...111111111111111111111111110011111110011111111111111111111111111111...

...111111111111111111111111111001111100111111111111111111111111111111...

...111111111111111111111111111100111001111111111111111111111111111111...

...11111111111111111111111111111001001111111111111111111111111111111111...

...111111111111111111111111111111000111111111111111111111111111111111...

...11111111111111111111111111111110111111111111111111111111111111111...

...111111111111111111111111111111000111111111111111111111111111111111...

...11111111111111111111111111111001001111111111111111111111111111111111...

...1111111111111111111111111111001100111111111111111111111111111111...
...111111111111111111111111111001111001111111111111111111111111111...
...11111111111111111111111111001111110011111111111111111111111111111...
...1111111111111111111111111001111111100111111111111111111111111111...
...111111111111111111111111001111111111001111111111111111111111111...
........................................................................

11) Infinite Numbers (II):

...111111111111111111111111111111121111111111111111111111111111111111...
...111111111111111111111111111111222111111111111111111111111111111111...
...111111111111111111111111111112232211111111111111111111111111111111...
...111111111111111111111111111111222111111111111111111111111111111111...
...111111111111111111111111111111121111111111111111111111111111111111...
...111111111111111111111111111111222111111111111111111111111111111111...
...111111111111111111111111111112232211111111111111111111111111111111...
...111111111111111111111111111122343221111111111111111111111111111111...
...111111111111111111111111111112232211111111111111111111111111111111...
...111111111111111111111111111111222111111111111111111111111111111111...
...111111111111111111111111111111121111111111111111111111111111111111...
...111111111111111111111111111111222111111111111111111111111111111111...
...111111111111111111111111111112232211111111111111111111111111111111...
...111111111111111111111111111122343221111111111111111111111111111111...
...111111111111111111111111111223454322111111111111111111111111111111...
...111111111111111111111111111122343221111111111111111111111111111111...
...111111111111111111111111111112232211111111111111111111111111111111...
...111111111111111111111111111111222111111111111111111111111111111111...
...111111111111111111111111111111121111111111111111111111111111111111...
...111111111111111111111111111111222111111111111111111111111111111111...
...111111111111111111111111111112232211111111111111111111111111111111...
...111111111111111111111111111122343221111111111111111111111111111111...
...111111111111111111111111111223454322111111111111111111111111111111...
...111111111111111111111111112234565432211111111111111111111111111111...
...111111111111111111111111111223454322111111111111111111111111111111...
...111111111111111111111111111122343221111111111111111111111111111111...
...111111111111111111111111111112232211111111111111111111111111111111...
...111111111111111111111111111111222111111111111111111111111111111111...
...111111111111111111111111111111121111111111111111111111111111111111...
...111111111111111111111111111111222111111111111111111111111111111111...
...111111111111111111111111111112232211111111111111111111111111111111...
...111111111111111111111111111122343221111111111111111111111111111111...
...111111111111111111111111111223454322111111111111111111111111111111...
...111111111111111111111111112234565432211111111111111111111111111111...
...111111111111111111111111122345676543221111111111111111111111111111...
...111111111111111111111111112234565432211111111111111111111111111111...
...111111111111111111111111111223454322111111111111111111111111111111...
...111111111111111111111111111122343221111111111111111111111111111111...
...111111111111111111111111111112232211111111111111111111111111111111...
...111111111111111111111111111111222111111111111111111111111111111111...
...111111111111111111111111111111121111111111111111111111111111111111...
...111111111111111111111111111111222111111111111111111111111111111111...
...111111111111111111111111111112232211111111111111111111111111111111...

...111111111111111111111111111223432211111111111111111111111111111111...
...11111111111111111111111111223454322111111111111111111111111111111...
...1111111111111111111111111223456543221111111111111111111111111111...
...111111111111111111111111223456765432211111111111111111111111111...
...11111111111111111111111223456787654322111111111111111111111111111...
.........................................................................

12) Numerical Car:

...000000000000000000000000000000000000000000000000000000000000000...
...0000000000000000011111111111111111111111100000000000000000000000...
...00000000000000000111111111111111111111111110000000000000000000000...
...00000000000000011000000000000000000000000011000000000000000000000...
...000000000000001100000000000000000000000000001100000000000000000000...
...0000000111111111000000000000000000000000000000011111111111110000000...
...000000111111111000000000000000000000000000000001111111111111000000...
...0000001100000000000000000000000000000000000000000000000000011200000...
...0000001100000000000000000000000000000000000000000000000000011000000...
...0000001100000444000000000000000000000000000000000004440011000000...
...00000011111144444111111111111111111111111111111111444441111200000...
...0000000111144444441111111111111111111111111111111444444410000000...
...000000000000444440000000000000000000000000000000000444440000000000000...
...0000000000000444000000000000000000000000000000000000044400000000000000...
...000000000000000000000000000000000000000000000000000000000000000...
..........................................................................

13) Finite Lattice:

...000000000000000000000000000000000000000000000000000000000000000...
...0777000000000000700000007777777700777777770077007777777700777777770...
...0777000000000007770000007777777700777777770077007777777700777777770...
...0777000000000077077000000000770000000007700000770077700000000770000000...
...0777000000000770007700000000770000000007700000770077700000000777770000...
...0777000000007777777700000000770000000007700000770077700000000770000000...
...0777770007700000000077000007700000000077000007700777777777700777777770...
...0777770077000000000077000077000000000077000007700777777777700777777770...
...000000000000000000000000000000000000000000000000000000000000000...
..........................................................................

14) Infinite Lattice:

...1111111111111111111111111111111111111111111111111111111111111111...
...17771111111111117111111177777777117777777711771177777777117777777771...
...17771111111111177711111177777777117777777711771177777777117777777771...
...177711111111117717711111111177111111111771111177117771111111177111111...
...177711111111177111771111111177111111111771111177117771111111177777111...
...1777111111117777777711111117711111111177111111771177711111111177111111...
...17777771117711111111177111117711111111177111117711777777777117777777771...
...17777771177111111111177111177111111111177111117711777777777117777777771...
...1111111111111111111111111111111111111111111111111111111111111111...
..........................................................................

Remark: Of course, it's interesting to "design" a large variety of
numerical <object sequences> in the same way.
Their numbers may be finite if the picture's background is zeroed, or infinite if the picture's background is not zeroed -- as for the previous examples.

15) Simple Numbers:

2,3,4,5,6,7,8,9,10,11,13,14,15,17,19,21,22,23,25,26,27,29,31,33,34,35,37,38, 39,41,43,45,46,47,49,51,53,55,57,58,61,62,65,67,69,71,73,74,77,78,79,82,83, 85,86,87,89,91,93,94,95,97,101,103,...

(A number n is called <simple number> if the product of its proper divisors is less than or equal to n.)

Generally speaking, n has the form:

$n = p$, or $p^2$, or $p^3$, or $p^q$, where p and q are distinct primes.

16) Digital Sum:

$\underbrace{0,1,2,3,4,5,6,7,8,9,}\ \underbrace{1,2,3,4,5,6,7,8,9,10,}\ \underbrace{2,3,4,5,6,7,8,9,10,11,}$

$\underbrace{3,4,5,6,7,8,9,10,11,12,}\ \underbrace{4,5,6,7,8,9,10,11,12,13,}\ \underbrace{5,6,7,8,9,10,11,12,13,14,}\ \cdots$

( $d_s(n)$ is the sum of digits.)

17) Digital Products:

$\underbrace{0,1,2,3,4,5,6,7,8,9,}\ \underbrace{0,1,2,3,4,5,6,7,8,9,}\ \underbrace{0,2,4,6,8,19,12,14,16,18,}$

$\underbrace{0,3,6,9,12,15,18,21,24,27,}\ \underbrace{0,4,8,12,16,20,24,28,32,36,}\ \underbrace{0,5,10,15,20,25,...}$

( $d_p(n)$ is the product of digits.)

18) Code Puzzle:

151405,202315,2008180505,06152118,06092205,190924,1905220514,0509070820,
14091405,200514,051205220514,...

Using the following letter-to-number code:

| | | | | | | | | | | | | |
|---|---|---|---|---|---|---|---|---|---|---|---|---|
| A | B | C | D | E | F | G | H | I | J | K | L | M |
| 01 | 05 | 06 | 07 | 08 | 09 | 10 | 11 | 12 | 13 | 14 | 15 | 16 |
| N | O | P | Q | R | S | T | U | V | W | X | Y | Z |
| 14 | 15 | 16 | 17 | 18 | 19 | 20 | 21 | 22 | 23 | 24 | 25 | 26 |

then $c_p(n)$= the numerical code for the spelling of n in English language; for example: 1 = ONE = 151405, etc.

19)) Pierced Chain:

101,1010101,10101010101,101010101010101,1010101010101010101,
10101010101010101010101,101010101010101010101010101,...

$$(c(n) = 101 * 1\underbrace{0001}_{1}\,\underbrace{0001}_{2}\ldots\underbrace{0001}_{n-1}, \text{ for } n \geq 1.)$$

How many c(n)/101 are primes ?

20) Divisor Products:

1,2,3,8,5,36,7,64,27,100,11,1728,13,196,225,1024,17,5832,19,8000,441,484,
23,331776,125,676,729,21952,29,810000,31,32768,1089,1156,1225,100776
96,37,1444,1521,2560000,41,...

( $P_d(n)$ is the product of all positive divisors of n.)

21) Proper Divisor Products:

1,1,1,2,1,6,1,8,3,10,1,144,1,14,15,64,1,324,1,400,21,22,1,13824,5,26,27, 784,1,27000,1,1024,33,34,35,279936,1,38,39,64000,1,...

( $P_d(n)$ is the product of all positive divisors of n but n.)

22) Cube Free Sieve:

2,3,4,5,6,7,9,10,11,12,13,14,15,17,18,19,20,21,22,23,25,26,28,29,30,31,33, 34,35,36,37,38,39,41,42,43,44,45,46,47,49,50,51,52,53,55,57,58,59,60,61,62, 63,65,66,67,68,69,70,71,73,...

Definition: from the set of natural numbers (except 0 and 1):

- take off all multiples of $2^3$ (i.e. 8, 16, 24, 32, 40, ...)
- take off all multiples of $3^3$
- take off all multiples of $5^3$

... and so on (take off all multiples of all cubic primes).

(One obtains all cube free numbers.)

23) m-Power Free Sieve:

Definition: from the set of natural numbers (except 0 and 1)

take off all multiples of $2^m$, afterwards all multiples of $3^m$, ...

and so on (take off all multiples of all m-power primes, $m \geq 2$).

(One obtains all m-power free numbers.)

24) Irrational Root Sieve:

2,3,5,6,7,10,11,12,13,14,15,17,18,19,20,21,22,23,24,26,28,29,30,31,33,34, 35,37,38,39,40,41,42,43,44,45,46,47,48,50,51,52,53,54,55,56,57,58,59,60,61, 62,63,65,66,67,68,69,70,71,72,73,...

Definition: from the set of natural numbers (except 0 and 1):

- take off all powers of $2^k$, $k \geq 2$, (i.e. 4, 8, 16, 32, 64, ...)
- take off all powers of $3^k$, $k \geq 2$;
- take off all powers of $5^k$, $k \geq 2$;
- take off all powers of $6^k$, $k \geq 2$;
- take off all powers of $7^k$, $k \geq 2$;
- take off all powers of $10^k$, $k \geq 2$;

... and so on (take off all k-powers, $k \geq 2$, of all square free numbers).

We got all square free numbers by the following method (sieve):

from the set of natural numbers (except 0 and 1):

- take off all multiples of $2^2$ (i.e. 4, 8, 12, 16, 20, ...)
- take off all multiples of $3^3$
- take off all multiples of $5^3$

... and so on (take off all multiples of all square primes);
one obtains, therefore:

2,3,5,6,7,10,11,13,14,15,17,19,21,22,23,26,29,30,31,33,34,35,37,38,39, 41,42,43,46,47,51,53,55,57,58,59,61,62,65,66,67,69,70,71,... ,

which are used for irrational root sieve.

(One obtains all natural numbers those m-th roots, for any $m \geq 2$, are irrational.)

25) Odd Sieve:

7,13,19,23,25,31,33,37,43,47,49,53,55,61,63,67,73,75,79,83,85,91,93, 97,...

(All odd numbers that are not equal to the difference of two primes.)

A sieve is used to get this sequence:

- subtract 2 from all prime numbers and obtain a temporary sequence;
- choose all odd numbers that do not belong to the temporary one.

26) Binary Sieve:

1,3,5,9,11,13,17,21,25,27,29,33,35,37,43,49,51,53,57,59,65,67,69,73,75,77, 81,85,89,91,97,101,107,109,113,115,117,121,123,129,131,133,137,139,145, 149,...

(Starting to count on the natural numbers set at any step from 1:

- delete every 2-nd numbers
- delete, from the remaining ones, every 4-th numbers

... and so on: delete, from the remaining ones, every ($2^k$) –th numbers, k = 1, 2, 3, ... .)

Conjectures:
- there are an infinity of primes that belong to this sequence;
- there are an infinity of numbers of this sequence which are not prime.

27) Trinary Sieve:

1,2,4,5,7,8,10,11,14,16,17,19,20,22,23,25,28,29,31,32,34,35,37,38,41,43,46, 47,49,50,52,55,56,58,59,61,62,64,65,68,70,71,73,74,76,77,79,82,83,85,86,88, 91,92,95,97,98,100,101,103,104,106,109,110,112,113,115,116,118,119,122, 124,125,127,128,130,131,133,137,139,142,143,145,146,149,...

(Starting to count on the natural numbers set at any step from 1:
- delete every 3-rd numbers
- delete, from the remaining ones, every 9-th numbers

... and so on: delete, from the remaining ones, every ($3^k$) -th numbers, k = 1, 2, 3, ... .)

Conjectures:
- there are an infinity of primes that belong to this sequence;
- there are an infinity of numbers of this sequence which are not prime.

28) n-ary Power Sieve (generalization, $n \geq 2$):

(Starting to count on the natural numbers set at any step from 1:
- delete every n- th numbers
- delete, from the remaining ones, every ($n^2$) –th numbers

... and so on: delete, from the remaining ones, every ($n^k$) –th numbers, k = 1, 2, 3, ... .)

Conjectures:
- there are an infinity of primes that belong to this sequence;
- there are an infinity of numbers of this sequence which are not prime.

29) k-ary Consecutive Sieve (second version):

1, 2, 4, 7, 9, 14, 20, 25, 31, 34, 44, ... .

Keep the first k numbers, skip the k+1 numbers, for k = 2, 3, 4, ... .

Question: How many terms are prime?

30) Consecutive Sieve:

1,3,5,9,11,17,21,29,33,41,47,57,59,77,81,101,107,117,131,149,153,173,191, 209,213,239,257,273,281,321,329,359,371,401,417,441,435,491,...

(From the natural numbers set:

- keep the first number,
  delete one number out of 2 from all remaining numbers;
- keep the first remaining number,
  delete one number out of 3 from the next remaining numbers;
- keep the first remaining number,
  delete one number out of 4 from the next remaining numbers;

... and so on, for step k ($k \geq 2$):

- keep the first remaining number,
  delete one number out of k from the next remaining numbers;

... .)

This sequence is much less dense than the prime number sequence, and their ratio tends to $p_n$: n as n tends to infinity.

For this sequence we chose to keep the first remaining number at all steps,
but in a more general case:
the kept number may be any among the remaining k-plet (even at random).

31) General-Sequence Sieve:

Let $u_i > 1$, for i = 1, 2, 3, ..., a strictly increasing positive integer sequence. Then:

From the natural numbers set:

- keep one number among 1, 2, 3, ..., $u_1$-1,
  and delete every ($u_1$)–th numbers;
- keep one number among the next $u_2$-1 remaining numbers,
  and delete every ($u_2$)–th numbers;

... and so on, for step k ($k \geq 1$):

- keep one number among the next $u_k$-1 remaining numbers,
  and delete every $u_k$ -th numbers;

...

Problem: study the relationship between sequence $u_i$, i = 1, 2, 3, ..., and the remaining sequence resulted from the general sieve.

$u_i$, previously defined, is called sieve generator.

32) More General-Sequence Sieve:

For i = 1, 2, 3, ..., let $u_i > 1$, be a strictly increasing positive integer sequence, and $v_i < u_i$ another positive integer sequence. Then:

From the natural numbers set:

- keep the ($v_1$ )–th number among 1, 2, 3, ..., u - 1, and delete every ($u_1$ )–th numbers;
- keep the($v_2$ )–th number among the next $u_2$ -1 remaining numbers, and delete every $u_2$-th numbers;

... and so on, for step k (k ≥ 1):

- keep the ($v_k$ )–th number among the next $u_k$-1 remaining numbers, and delete every ($u_k$ )–th numbers;

... .

Problem: study the relationship between sequences $u_i$, $v_i$, i = 1, 2, 3, ..., and the remaining sequence resulted from the more general sieve.

$u_i$ and $v_i$ previously defined, are called *sieve generators*.

33) Digital Sequences:

(This a particular case of sequences of sequences.)

General definition:

in any numeration base B, for any given infinite integer or rational sequence $S_1$ , $S_2$, $S_3$, ..., and any digit D from 0 to B-1, it's built up a new integer sequence witch associates to $S_1$ the number of digits D of $S_1$ in base B, to $S_2$ the number of digits D of $S_2$ in base B, and so on...

For example, considering the prime number sequence in base 10, then the number of digits 1 (for example) of each prime number following their order is: 0,0,0,0,2,1,1,1,0,0,1,0,...

(Digit-1 Prime Sequence).

Second example if we consider the factorial sequence n! in base 10, then the number of digits 0 of each factorial number following their order is: 0,0,0,0,0,1,1,2,2,1,3,...

Digit-0 Factorial Sequence).

Third example if we consider the sequence $n^n$ in base 10, n=1,2,..., then the number of digits 5 of each term $1^1$, $2^2$, $3^3$, ...,

following their order is: 0,0,0,1,1,1,1,0,0,0,...
(Digit-5 $n^n$-Sequence).

34) Construction Sequences:

(This a particular case of sequences of sequences.)
General definition:
in any numeration base B, for any given infinite integer or rational sequence $S_1$, $S_2$, $S_3$, ..., and any digits $D_1$, $D_2$, ..., $D_K$ ($k < B$), it's built up a new integer sequence such that
each of its terms $Q_1 < Q_2 < Q_3 < \ldots$ is formed by these digits $D_1$, $D_2$, ..., $D_K$ only (all these digits are used), and matches a term $S_i$ of the previous sequence.

For example, considering in base 10 the prime number sequence, and the digits 1 and 7 (for example),
we construct a written-only-with-these-digits (all these digits are used) prime number new sequence: 17,71,...
(Digit-1-7-only Prime Sequence).

Second example, considering in base 10 the multiple of 3 sequence, and the digits 0 and 1,
we construct a written-only-with-these-digits (all these digits are used) multiple of 3 new sequence: 1011,1101,1110,10011,10101,10110, 11001,11010,11100,...
(Digit-0-1-only Multiple-of-3 Sequence).

35) General Residual Sequence:

$(x + C_1)...(x + C_{F(m)})$, m = 2, 3, 4, ...,

where $C_i$, $1 \leq i \leq F(m)$, forms a reduced set of residues mod m, x is an integer, and F is Euler's totient.

The General Residual Sequence is induced from the The Residual Function (see <Libertas Mathematica>):

Let $L : Z \times Z \longrightarrow Z$ be a function defined by

$$L(x, m)=(x + C_1)...(x + C_{F(m)}),$$

where $C_i$, $1 \leq i \leq F(m)$, forms a reduced set of residues mod m, $m \geq 2$, x is an integer, and F is Euler's totient.

The Residual Function is important because it generalizes the classical theorems by Wilson, Fermat, Euler, Wilson, Gauss, Lagrange, Leibnitz, Moser, and Sierpinski all together.

For x=0 it's obtained the following sequence:

$L(m) = C_1 ... C_{F(m)}$ , where m = 2, 3, 4, ...

(the product of all residues of a reduced set mod m):

1,2,3,24,5,720,105,2240,189,3628800,385,479001600,19305,896896,2027025,

20922789888000,85085,6402373705728000,8729721,47297536000,1249937325,...

which is found in "The Encyclopedia of Integer Sequences", by N. J. A. Sloane, Academic Press, San Diego, New York, Boston, London, Sydney, Tokyo, Toronto, 1995.

The Residual Function extends it.

36) Inferior f-Part of x.

Let f: Z ↦ Z be a strictly increasing function and $x \in R$. Then:
ISf(x) is the smallest k such that $f(k) \leq x < f(k+1)$.

37) Superior f-Part of x.

Let f: Z ↦ Z be a strictly increasing function and $x \in R$. Then:
SSf(x) is the smallest k such that $f(k) < x \leq f(k+1)$.

38) (Inferior) Prime Part:

2,3,3,5,5,7,7,7,7,11,11,13,13,13,13,17,17,19,19,19,19,23,23,23,23,23,23,
29,29,31,31,31,31,31,31,37,37,37,37,41,41,43,43,43,43,47,47,47,47,47,47,
53,53,53,53,53,53,59, ...

(For any positive real number n one defines $p_p(n)$ as the largest prime number less than or equal to n.)

39) (Superior) Prime Part:

2,2,2,3,5,5,7,7,11,11,11,11,13,13,17,17,17,17,19,19,23,23,23,23,29,29,29,
29,29,29,31,31,37,37,37,37,37,37,41,41,41,41,43,43,47,47,47,47,53,53,53,
53,53,53,59,59,59,59,59,59,61,...

(For any positive real number n one defines $p_p(n)$ as the smallest prime number greater than or equal to n.)

40) (Inferior) Square Part:

0,1,1,1,4,4,4,4,4,9,9,9,9,9,9,9,16,16,16,16,16,16,16,16,16,25,25,25,25,25,
25,25,25,25,25,25,36,36,36,36,36,36,36,36,36,36,36,36,36,49,49,49,49,49,
49,49,49,49,49,49,49,49,49,49,64,64,...

(The largest square less than or equal to n.)

41) (Superior) Square Part:

0,1,4,4,4,9,9,9,9,9,16,16,16,16,16,16,16,25,25,25,25,25,25,25,25,25,36,36,

36,36,36,36,36,36,36,36,36,49,49,49,49,49,49,49,49,49,49,49,49,49,64,64,64,
64,64,64,64,64,64,64,64,64,64,64,64,64,81,81,...
(The smallest square greater than or equal to n.)

42) (Inferior) Cube Part:
0,1,1,1,1,1,1,1,8,8,8,8,8,8,8,8,8,8,8,8,8,8,8,8,8,8,8,27,27,27,27,27,27,
27,27,27,27,27,27,27,27,27,27,27,27,27,27,27,27,27,27,27,27,27,27,27,27,27,
27,27,27,27,27,27,27,64,64,64,...
(The largest cube less than or equal to n.)

43) (Superior) Cube Part:
0,1,8,8,8,8,8,8,8,27,27,27,27,27,27,27,27,27,27,27,27,27,27,27,27,27,27,27,
27,64,64,64,64,64,64,64,64,64,64,64,64,64,64,64,64,64,64,64,64,64,64,64,64,
64,64,64,64,64,64,64,64,64,64,64,64,64,64,125,125,125,..
(The smallest cube greater than or equal to n.)

44) (Inferior) Factorial Part:
1,2,2,2,2,6,6,6,6,6,6,6,6,6,6,6,6,6,6,6,6,6,6,24,24,24,24,24,24,24,24,24,
24,24,24,24,24,24,24,24,24,...
( $F_p(n)$ is the largest factorial less than or equal to n.)

45) (Superior) Factorial Part:
1,2,6,6,6,6,24,24,24,24,24,24,24,24,24,24,24,24,24,24,24,24,24,24,120,120,
120,120,120,120,120,120,120,120,120,...
( $f_p(n)$ is the smallest factorial greater than or equal to n.)

46) Inferior Fractional f-Part of x.
Let f: $Z \mapsto Z$ be a strictly increasing function and $x \in R$. Then:
IFSf(x) = x - ISf(x),
where ISf(x) is the Inferior f-Part of x defined above.
For example: Inferior Fractional Cubic Part of 12.501
= 12.501 - 8 = 4.501.

47) Superior Fractional f-Part of x.
Let f: $Z \mapsto Z$ be a strictly increasing function and $x \in R$. Then:
SFSf(x) = SSf(x)-x,
where SSf(x) is the Superior f-Part of x defined above.
For example: Superior Fractional Cubic Part of 12.501

= 27 - 12.501 = 14.499.

The Fractional f-Parts are generalizations of the inferior and respectively superior fractional part of a number.

Particular cases:

48) Fractional Prime Part:

FSp(x) = x - ISp(x),

where ISp(x) is the Inferior Prime Part defined above.

Example: FSp(12.501) = 12.501 - 11 = 1.501.

49) Fractional Square Part:

FSs(x) = x - ISs(x),

where ISs(x) is the Inferior Square Part defined above.

Example: FSs(12.501) = 12.501 - 9 = 3.501.

50) Fractional Cubic Part:

FSc(x) = x - ISc(x),

where ISc(x) is the Inferior Cubic Part defined above.

Example: FSc(12.501) = 12.501 - 8 = 4.501.

51) Fractional Factorial Part:

FSf(x) = x - ISf(x),

where ISf(x) is the Inferior Factorial Part defined above.

Example: FSf(12.501) = 12.501 - 6 = 6.501.

52) Smarandacheian Complements.

Let g: A ↦ be a strictly increasing function, and let "~" be a given internal law on A.

Then f: A ↦ A is a smarandacheian complement with respect to the function g and the internal law "~" if:

f(x) is the smallest k such that there exists a z in A so that x~k = g(z).

53) Square Complements:

1,2,3,1,5,6,7,2,1,10,11,3,14,15,1,17,2,19,5,21,22,23,6,1,26,3,7,29,30,31, 2,33,34,35,1,37,38,39,10,41,42,43,11,5,46,47,3,1,2,51,13,53,6,55,14,57,58, 59,15,61,62,7,1,65,66,67,17,69,70,71,2,...

Definition:

for each integer n to find the smallest integer k such that

n·k is a perfect square..
(All these numbers are square free.)

54) Cubic Complements:
1,4,9,2,25,36,49,1,3,100,121,18,169,196,225,4,289,12,361,50,441,484,529, 9,5,676,1,841,900,961,2,1089,1156,1225,6,1369,1444,1521,25,1681,1764,1849, 242,75,2116,2209,36,7,20,...

Definition:
for each integer n to find the smallest integer k such that n·k is a perfect cub.
(All these numbers are cube free.)

55) m-Power Complements:

Definition:
for each integer n to find the smallest integer k such that n·k is a perfect m-power ($m \geq 2$).
(All these numbers are m-power free.)

56) Double Factorial Complements:
1,1,1,2,3,8,15,1,105,192,945,4,10395,46080,1,3,2027025,2560,34459425,192,

5,3715891200,13749310575,2,81081,1961990553600,35,23040,213458046676875,

128,6190283353629375,12,...

(For each n to find the smallest k such that n·k is a double factorial,
i.e. n·k = either 1·3·5·7·9·...·n if n is odd,
either 2·4·6·8·...·n if n is even.)

57) Prime Additive Complements:
1,0,0,1,0,1,0,3,2,1,0,1,0,3,2,1,0,1,0,3,2,1,0,5,4,3,2,1,0,1,0,5,4,3,2,1,0,
3,2,1,0,1,0,3,2,1,0,5,4,3,2,1,0,...
(For each n to find the smallest k such that n+k is prime.)
Remark: Is it possible to get as large as we want but finite decreasing k, k-1, k-2, ..., 2, 1, 0 (odd k) sequence included in the previous sequence, i.e. for any even integer are there two primes those difference is equal to it? We conjecture the answer is negative.

58) Prime Base:

0,1,10,100,101,1000,1001,10000,10001,10010,10100,100000,100001,1000000,

1000001,1000010,1000100,10000000,10000001,100000000,100000001,100000010,

100000100,1000000000,1000000001,1000000010,1000000100,1000000101,...
(Each number n written in the prime base.)

(We defined over the set of natural numbers the following infinite base: $p_0 = 1$, and for $k \geq 1$ $p_k$ is the k-th prime number.)
We proved that every positive integer A may be uniquely written in the prime base as:

$$A = (\overline{a_n \cdots a_1 a_0})_{(SP)} \overset{def}{=} \sum_{i=1}^{n} a_i p_i \text{, with all } a_i = 0 \text{ or } 1, \text{ (of course } a_n = 1),$$

in the following way:

- if $p_n \leq A < p_{n+1}$ then $A = p_n + r_1$ ;
- if $p_m \leq r_1 < p_{m+1}$ then $r_1 = p_m + r_2$, $m < n$;

and so on until one obtains a rest $r_j = 0$.

Therefore, any number may be written as a sum of prime numbers + e, where e = 0 or 1.

If we note by p(A) the superior part of A (i.e. the largest prime less than or equal to A), then A is written in the prime base as:

$$A = p(A) + p(A-p(A)) + p(A-p(A)-p(A-p(A))) + \ldots .$$

This base is important for partitions with primes.

59) Square Base:

0,1,2,3,10,11,12,13,20,100,101,102,103,110,111,112,1000,1001,1002,1003,

1010,1011,1012,1013,1020,10000,10001,10002,10003,10010,10011,10012,10013,

10020,10100,10101,100000,100001,100002,100003,100010,100011,100012,100013,

100020,100100,100101,100102,100103,100110,100111,100112,101000,101001,

101002,101003,101010,101011,101012,101013,101020,101100,101101,101102,

1000000,...

(Each number n written in the square base.)

(We defined over the set of natural numbers the following infinite base: for $k \geq 0$ $s_k = k^2$.)

We proved that every positive integer A may be uniquely written in the square base as:

$$A = (\overline{a_n \cdots a_1 a_0})_{(S2)} \overset{def}{=} \sum_{i=0}^{n} a_i s_i \text{, with } a_i = 0 \text{ or } 1 \text{ for } i \geq 2,$$

$0 \leq a_0 \leq 3$, $0 \leq a_1 \leq 2$, and of course $a_n = 1$,

in the following way:

- if $s_n \leq A < s_{n+1}$ then $A = s_n + r_1$;
- if $s_m \leq r_1 < p_{m+1}$ then $r_1 = s_m + r_2$, $m < n$;

and so on until one obtains a rest $r_j = 0$.

Therefore, any number may be written as a sum of squares (1 not counted as a square -- being obvious) + e, where e = 0, 1, or 3.

If we note by s(A) the superior square part of A (i.e. the largest square less than or equal to A), then A is written in the square base as:

A = s(A) + s(A-s(A)) + s(A-s(A)-s(A-s(A))) + ... .

This base is important for partitions with squares.

60) m-Power Base (generalization):
(Each number n written in the m-power base, where m is an integer ≥ 2.)

(We defined over the set of natural numbers the following infinite m-power base: for $k \geq 0$ $t_k = k^m$ )

We proved that every positive integer A may be uniquely written in the m-power base as:

$$A = (\overline{a_n \cdots a_1 a_0})_{(SM)} \stackrel{def}{=} \sum_{i=0}^{n} a_i t_i$$, with $a_i = 0$ or 1 for $i \geq m$,

$$0 \leq a_i \leq \left[ \frac{(i+2)^m - 1}{(i+1)^m} \right]$$ (integer part)

for i = 0, 1, ..., m-1, $a_i = 0$ or 1 for i >= m, and of course $a_n = 1$, in the following way:

- if $t_n \leq A < t_{n+1}$ then $A = t_n + r_1$;
- if $t_m \leq r_1 < t_{m+1}$ then $r_1 = t_m + r_2$, m< n;

and so on until one obtains a rest $r_j = 0$.

Therefore, any number may be written as a sum of m-powers (1 not counted as an m-power -- being obvious) + e, where e = 0, 1, 2, ..., or $2^m$-1.

If we note by t(A) the superior m-power part of A (i.e. the largest m-power less than or equal to A), then A is written in the m-power base as:

A = t(A) + t(A-t(A)) + t(A-t(A)-t(A-t(A))) + ...

This base is important for partitions with m-powers.

61) Factorial Base:
0,1,10,11,20,21,100,101,110,111,120,121,200,201,210,211,220,221,300,301,310,311,320,321,1000,1001,1010,1011,1020,1021,1100,1101,1110,1111,1120,1121,1200,...
(Each number n written in the factorial base.)

(We defined over the set of natural numbers the following infinite base: for $k \geq 1$ $f_k = k!$)

We proved that every positive integer A may be uniquely written in the square base as:

$$A = \overline{(a_n \cdots a_2 a_1)}_{(F)} \overset{def}{=} \sum_{i=1}^{n} a_i f_i \text{, with all } a_i = 0, 1, \ldots, i \text{ for } i \geq 1,$$

in the following way:
- if $f_n \leq A < f_{n+1}$ then $A = f_n + r_1$;
- if $f_m \leq r_1 < f_{m+1}$ then $r_1 = f_m + r_2$, $m < n$

and so on until one obtains a rest $r_j = 0$.

What's very interesting: $a_1 = 0$ or 1; $a_2 = 0$, 1, or 2; $a_3 = 0$, 1, 2, or 3, and so on...

If we note by f(A) the superior factorial part of A (i.e. the largest factorial less than or equal to A), then A is written in the factorial base as:

$$A = f(A) + f(A-f(A)) + f(A-f(A)-f(A-f(A))) + \ldots .$$

Rules of addition and subtraction in factorial base:
For each digit $a_i$ we add and subtract in base i+1, for $i \geq 1$.
For example, an addition:

```
base   5 4 3 2
         2 1 0 +
         2 2 1
       1 1 0 1
```

because: 0+1= 1 (in base 2);
1+2=10 (in base 3), therefore we write 0 and keep 1;
2+2+1=11 (in base 4).

Now a subtraction:

```
          base   5 4 3 2
         ----------------
                 1 0 0 1 –
                   3 2 0
         ----------------
                 = = 1 1
```

because: 1-0=1 (in base 2);

0-2=? it's not possible (in base 3),

go to the next left unit, which is 0 again (in base 4),

go again to the next left unit, which is 1 (in base 5),

therefore 1001 ⟶ 0401 ⟶ 0331

and then 0331-320=11.

Find some rules for multiplication and division.

In a general case:

if we want to design a base such that any number

$$A = (\overline{a_n \cdots a_2 a_1})_{(B)} \overset{def}{=} \sum_{i=1}^{n} a_i b_i$$

, with all $a_i = 0, 1, \ldots, t_i$ for $i \geq 1$, where all $t_i \geq 1$,

then:

this base should be

$b_1 = 1$, $b_{i+1} = (t_i + 1) * b_i$ for $i \geq 1$.

62) Double Factorial Base:

1,10,100,101,110,200,201,1000,1001,1010,1100,1101,1110,1200, 10000,10001,10010,10100,10101,10110,10200,10201,11000,11001, 11010,11100,11101,11110,11200,11201,12000,...

( Numbers written in the double factorial base, defined as follows: df(n) = n!! )

63) Triangular Base:

1,2,10,11,12,100,101,102,110,1000,1001,1002,1010,1011,10000, 10001,10002,10010,10011,10012,100000,100001,100002,100010, 100011,100012,100100,1000000,1000001,1000002,1000010,1000011, 1000012,1000100,...

( Numbers written in the triangular base, defined as follows: $t(n) = \frac{n(n+1)}{2}$, for $n \geq 1$. )

64) Generalized Base:

(Each number n written in the generalized base.)

(We defined over the set of natural numbers the following infinite generalized base: $1 = g_0 < g_1 < ... < g_k < ... .$)

We proved that every positive integer A may be uniquely written in the generalized base as:

$$A = (\overline{a_n \cdots a_1 a_0})_{(SG)} \overset{def}{=} \sum_{i=0}^{n} a_i g_i \text{, with } 0 \le a_i \le \left[ \frac{(g_{i+1} - 1)}{g_i} \right]$$

(integer part) for i = 0, 1, ..., n, and of course $a_n \ge 1$,

in the following way:

- if $g_n \le A < g_{n+1}$ then $A = g_n + r_1$;
- if $g_m \le r_1 < g_{m+1}$ then $r_1 = g_m + r_2$, $m < n$

and so on until one obtains a rest $r_j = 0$.

If we note by g(A) the superior generalized part of A (i.e. the largest $g_i$ less than or equal to A), then A is written in the m-power base as:

$$A = g(A) + g(A-g(A)) + g(A-g(A)-g(A-g(A))) + ...$$

This base is important for partitions: the generalized base may be any infinite integer set (primes, squares, cubes, any m-powers, Fibonacci/Lucas numbers, Bernoully numbers, Smarandache numbers, etc.) those partitions are studied.

A particular case is when the base verifies: $2g_i \ge g_{i+1}$ any i, and $g_0 = 1$, because all coefficients of a written number in this base will be 0 or 1.

Remark: another particular case: if one takes $g_i = p^{i-1}$, i = 1, 2, 3, ..., p an integer ≥ 2, one gets the representation of a number in the numerical base p {p may be 10 (decimal), 2 (binary), 16 (hexadecimal), etc.}.

65) Smarandache numbers:

1,2,3,4,5,3,7,4,6,5,11,4,13,7,5,6,17,6,19,5,7,11,23,4,10,13,9,7,29, 5,31,8,11,17,7,6,37,19,13,5,41,7,43,11,5,23,47,6,14,10,17,13,53,9,11, 7,19,29,59,5,61,31,7,8,13,...

(S(n) is the smallest integer such that S(n)! is divisible by n.)
Remark: S(n) are the values of Smarandache function.

66) Smarandache quotients:
1,1,2,6,24,1,720,3,80,12,3628800,2,479001600,360,8,45,20922789888000,
40,6402373705728000,6,240,1814400,1124000727777607680000,1,145152,
239500800,13440,180,304888344611713860501504000000,...
(For each n to find the smallest k such that n·k is a factorial number.)

67) Double Factorial Numbers:
1,2,3,4,5,6,7,4,9,10,11,6,13,14,5,6,17,12,19,10,7,22,23,6,15,26,9,14,29,
10,31,8,11,34,7,12,37,38,13,10,41,14,43,22,9,46,47,6,21,10,...
($d_f(n)$ is the smallest integer such that $d_f(n)!!$ is a multiple of n.)
where m!! is the double factorial of m, i.e.
m!! = 2·4·6·...·m if m is even,
and m!! = 1·3·5·...·m if m is odd.

68) Primitive Numbers (of power 2):

2,4,4,6,8,8,8,10,12,12,14,16,16,16,16,18,20,20,22,24,24,24,26,28,28,30,32, 32,32,32,32,34,36,36,38,40,40,40,42,44,44,46,48,48,48,48,50,52,52,54,56,56, 56,58,60,60,62,64,64,64,64,64,64,66,...

( $S_2(n)$ is the smallest integer such that $S_2(n)!$ is divisible by $2^n$.)

Curious property: This is the sequence of even numbers, each number being repeated as many times as its exponent (of power 2) is.

This is one of irreducible functions, noted $S_2(k)$, which helps to calculate the Smarandache function (called also Smarandache numbers in "The Encyclopedia of Integer Sequences", by N. J. A. Sloane and S. Plouffe, Academic Press, San Diego, New York, Boston, London, Sydney, Tokyo, Toronto, 1995).

69) Primitive Numbers (of power 3):

3,6,9,9,12,15,18,18,21,24,27,27,27,30,33,36,36,39,42,45,45,48,51,54,54,54, 57,60,63,63,66,69,72,72,75,78,81,81,81,81,84,87,90,90,93,96,99,99,102,105, 108,108,108,111,...

( $S_3(n)$ is the smallest integer such that $S_3(n)!$ is divisible by $3^n$.)

Curious property: this is the sequence of multiples of 3, each number being repeated as many times as its exponent (of power 3) is.

This is one of irreducible functions, noted $S_3(k)$, which helps to calculate the Smarandache function (called also Smarandache numbers in "The Encyclopedia of Integer Sequences", by N. J. A. Sloane and S. Plouffe, Academic Press, San Diego, New York, Boston, London, Sydney, Tokyo, Toronto, 1995).

70) Primitive Numbers (of power p, p prime) {generalization}:

( $S_p(n)$ is the smallest integer such that $S_p(n)!$ is divisible by $p^n$.)

Curious property: this is the sequence of multiples of p, each number being repeated as many times as its exponent (of power p) is.

These are the irreducible functions, noted $S_p(k)$, for any prime number p, which helps to calculate the Smarandache function (called also Smarandache numbers in "The Encyclopedia of Integer Sequences", by N. J. A. Sloane and S. Plouffe, Academic Press, San Diego, New York, Boston, London, Sydney, Tokyo, Toronto, 1995).

71) Smarandache Functions of the First Kind:

$S_n : N^* \mapsto N^*$

i) If $n = u^r$ (with u = 1, or u = p prime number), then
$S_n(a) = k$, where k is the smallest positive integer such that
k! is a multiple of $u^{r \cdot a}$ ;

ii) If $n = p1^{r1} \cdot p2^{r2} \cdot \ldots \cdot pt^{rt}$ , then
$S_n(a) = \max_{1 \le j \le t} \{ S_{pj^\wedge rj}(a) \}$.

72) Smarandache Functions of the Second Kind:

$S^k : N^* \mapsto N^*$, $S^k(n) = S_n(k)$ for $k \in N^*$,
where $S_n$ are the Smarandache functions of the first kind.

73) Smarandache Function of the Third Kind:

$S_a^b(n) = S_{an}(b_n)$, where $S_{an}$ is the Smarandache function of the first kind, and the sequences $(a_n)$ and $(b_n)$ are different from the following situations:

i) $a_n = 1$ and $b_n = n$, for $n \in N^*$;

ii) $a_n = n$ and $b_n = 1$, for $n \in N^*$.

74) Pseudo-Smarandache Numbers Z(n):

Z(n) is the smallest integer such that 1 + 2 + ... + Z(n) is divisible by n.

For example:

| n | 1 | 2 | 3 | 4 | 5 | 6 | 7 |
|---|---|---|---|---|---|---|---|
| Z(n) | 1 | 3 | 2 | 3 | 4 | 3 | 6 |

75) Square Residues:

1,2,3,2,5,6,7,2,3,10,11,6,13,14,15,2,17,6,19,10,21,22,23,6,5,26,3,14,29,30,
31,2,33,34,35,6,37,38,39,10,41,42,43,22,15,46,47,6,7,10,51,26,53,6,14,57,58,
59,30,61,62,21,...

( $S_r(n)$ is the largest square free number which divides n.)
Or, $S_r(n)$ is the number n released of its squares:

if n = $(P_1^{a_1}) * ... * (P_r^{a_r})$, with all $p_i$ primes and all $a_i \geq 1$,

then $S_r(n) = p_1 * ... * p_r$ .

Remark: at least the $(2^2)$*k-th numbers (k = 1, 2, 3, ...) are released of their squares;
and more general: all $(p^2)$*k-th numbers (for all p prime, and k = 1, 2, 3, ...) are released of their squares.

76) Cubical Residues:
1,2,3,4,5,6,7,4,9,10,11,12,13,14,15,4,17,18,19,20,21,22,23,12,25,26,9,28, 29,30,31,4,33,34,35,36,37,38,39,20,41,42,43,44,45,46,47,12,49,50,51,52,53, 18,55,28,...
( $c_r(n)$ is the largest cube free number which divides n.)

Or, $c_r(n)$ is the number n released of its cubicals:

if n = $(P_r^{a_1}) * ... * (P_r^{a_r})$, with all $p_i$ primes and all $a_i \geq 1$,

then $c_r(n) = (P_r^{b_1}) * ... * (P_r^{b_r})$, with all $b_i = \min \{2, a_i\}$.

Remark: at least the $(2^3)$*k-th numbers (k = 1, 2, 3, ...) are released of their cubicals;
and more general: all $(p^3)$*k-th numbers (for all p prime, and k = 1, 2, 3, ...) are released of their cubicals.

77) m-Power Residues (generalization):
( $m_r(n)$ is the largest m-power free number which divides n.)

Or, $m_r(n)$ is the number n released of its m-powers:

if n = $(P_r^{a_1}) * ... * (P_r^{a_r})$, with all $p_i$ primes and all $a_i \geq 1$,

then $m_r(n) = (P_r^{b_1}) * ... * (P_r^{b_r})$, with all $b_i = \min \{ m-1, a_i \}$.

Remark: at least the $(2^m)$*k-th numbers (k = 1, 2, 3, ...) are released of their m-powers;
and more general: all $(p^m)$*k-th numbers (for all p prime, and k = 1, 2, 3, ...) are released of their m-powers.

78) Exponents (of power 2):

0,1,0,2,0,1,0,3,0,1,0,2,0,1,0,4,0,1,0,2,0,1,0,2,0,1,0,2,0,1,0,5,0,1,0,2,0,
1,0,3,0,1,0,2,0,1,0,3,0,1,0,2,0,1,0,3,0,1,0,2,0,1,0,6,0,1,...

( $e_2(n)$ is the largest exponent (of power 2) which divides n.)

Or, $e_2(n) = k$ if $2^k$ divides n but $2^{k+1}$ does not.

79) Exponents (of power 3):

0,0,1,0,0,1,0,0,2,0,0,1,0,0,1,0,0,2,0,0,1,0,0,1,0,0,3,0,0,1,0,0,1,0,0,2,0,
0,1,0,0,1,0,0,2,0,0,1,0,0,1,0,0,2,0,0,1,0,0,1,0,0,2,0,0,1,0,...

( $e_3(n)$ is the largest exponent (of power 3) which divides n.)

Or, $e_3(n)= k$ if $3^k$ divides n but $3^{k+1}$ does not.

80) Exponents (of power p) {generalization}:

( $e_p(n)$ is the largest exponent (of power p) which divides n, where p is an integer $\geq 2$.)

Or, $e_p(n)= k$ if $p^k$ divides n but $p^{k+1}$ does not.

81) Pseudo-Primes of First Kind:

2,3,5,7,11,13,14,16,17,19,20,23,29,30,31,32,34,35,37,38,41,43,47,50,53,59,
61,67,70,71,73,74,76,79,83,89,91,92,95,97,98,101,103,104,106,107,109,110,
112,113,115,118,119,121,124,125,127,128,130,131,133,134,136,137,139,140,
142,143,145,146, ...

(A number is a pseudo-prime of first kind if some permutation of the digits is a prime number, including the identity permutation.)

(Of course, all primes are pseudo-primes of first kind,
but not the reverse!)

82) Pseudo-Primes of Second Kind:

14,16,20,30,32,34,35,38,50,70,74,76,91,92,95,98,104,106,110,112,115,118, 119,121,124,125,128,130,133,134,136,140,142,143,145,146, ...

(A composite number is a pseudo-prime of second kind if some permutation of the digits is a prime number.)

83) Pseudo-Primes of Third Kind:

11,13,14,16,17,20,30,31,32,34,35,37,38,50,70,71,73,74,76,79,91,92,95,97,98, 101,103,104,106,107,109,110,112,113,115,118,119,121,124,125,127,128,130, 131,133,134,136,137,139,140,142,143,145,146, ...

(A number is a pseudo-prime of third kind if some nontrivial permutation of the digits is a prime number.)

Question: How many pseudo-primes of third kind are prime numbers? (We conjecture: an infinity).

(There are primes which are not pseudo-primes of third kind, and the reverse: there are pseudo-primes of third kind which are not primes.)

84) Almost Primes of First Kind:

Let $a_1 \geq 2$, and for $n \geq 1$ $a_{n+1}$ is the smallest number that is not divisible by any of the previous terms (of the sequence) $a_1, a_2, ..., a_n$.

Example for $a_1 = 10$:

10,11,12,13,14,15,16,17,18,19,21,23,25,27,29,31,35,37,41,43,47,49,53,57, 61,67,71,73,...

If one starts by $a_1 = 2$, it obtains the complete prime sequence and only it.

If one starts by $a_1 > 2$, it obtains after a rank r, where $a_r = p(a_1)^2$ with p(x) the strictly superior prime part of x, i.e. the largest prime strictly less than x, the prime sequence:

- between $a_1$ and $a_r$, the sequence contains all prime numbers of this interval and some composite numbers;
- from $a_{r+1}$ and up, the sequence contains all prime numbers greater than $a_r$ and no composite numbers.

85) Almost Primes of Second Kind:

$a_1 \geq 2$, and for $n \geq 1$ $a_{n+1}$ is the smallest number that is coprime with all of the previous terms (of the sequence) , $a_1, a_2, ..., a_n$ .

This second kind sequence merges faster to the prime numbers than the first kind sequence.

Example for $a_1 = 10$:

10,11,13,17,19,21,23,29,31,37,41,43,47,53,57,61,67,71,73,...

If one starts by $a_1 = 2$, it obtains the complete prime sequence and only it.

If one starts by $a_1 > 2$, it obtains after a rank r, where $a_1 = p_i p_j$ with $p_i$ and $p_j$ prime numbers strictly less than and not dividing $a_1$ , the prime sequence:

- between $a_1$ and $a_r$ , the sequence contains all prime numbers of this interval and some composite numbers;
- from $a_{r+1}$ and up, the sequence contains all prime numbers greater than $a_r$ and no composite numbers.

86) Pseudo-Squares of First Kind:

1,4,9,10,16,18,25,36,40,46,49,52,61,63,64,81,90,94,100,106,108,112,121,136, 144,148,160,163,169,180,184,196,205,211,225,234,243,250,252,256,259,265, 279,289,295,297,298,306,316,324,342,360,361,400,406,409,414,418,423,432, 441,448,460,478,481,484,487,490,502,520,522,526,529,562,567,576,592,601, 603,604,610,613,619,625,630,631,640,652,657,667,675,676,691,729,748,756, 765,766,784,792,801,810,814,829,841,844,847,874,892,900,904,916,925,927, 928,940,952,961,972,982,1000, ...

(A number is a pseudo-square of first kind if some permutation of the digits is a perfect square, including the identity permutation.)

(Of course, all perfect squares are pseudo-squares of first kind, but not the reverse!)

One listed all pseudo-squares of first kind up to 1000.

87) Pseudo-Squares of Second Kind:

10,18,40,46,52,61,63,90,94,106,108,112,136,148,160,163,180,184,205,211,234, 243,250,252,259,265,279,295,297,298,306,316,342,360,406,409,414,418,423,

432,448,460,478,481,487,490,502,520,522,526,562,567,592,601,603,604,610,
613,619,630,631,640,652,657,667,675,691,748,756,765,766,792,801,810,814,
829,844,847,874,892,904,916,925,927,928,940,952,972,982,1000, ...

(A non-square number is a pseudo-square of second kind if some permutation of the digits is a square.)

One listed all pseudo-squares of second kind up to 1000.

88) Pseudo-Squares of Third Kind:

10,18,40,46,52,61,63,90,94,100,106,108,112,121,136,144,148,160,163,169,180,
184,196,205,211,225,234,243,250,252,256,259,265,279,295,297,298,306,316,
342,360,400,406,409,414,418,423,432,441,448,460,478,481,484,487,490,502,
520,522,526,562,567,592,601,603,604,610,613,619,625,630,631,640,652,657,
667,675,676,691,748,756,765,766,792,801,810,814,829,844,847,874,892,900,
904,916,925,927,928,940,952,961,972,982,1000,...

(A number is a pseudo-square of third kind if some nontrivial permutation of the digits is a square.)

Question: How many pseudo-squares of third kind are square numbers? (We conjecture: an infinity).

(There are squares which are not pseudo-squares of third kind, and the reverse:
there are pseudo-squares of third kind which are not squares.)

One listed all pseudo-squares of third kind up to 1000.

89) Pseudo-Cubes of First Kind:

1,8,10,27,46,64,72,80,100,125,126,152,162,207,215,216,251,261,270,279,297,
334,343,406,433,460,512,521,604,612,621,640,702,720,729,792,800,927,972,
1000,...

(A number is a pseudo-cube of first kind if some permutation of the digits is a cube, including the identity permutation.)

(Of course, all perfect cubes are pseudo-cubes of first kind, but not the reverse!)

One listed all pseudo-cubes of first kind up to 1000.

90) Pseudo-Cubes of Second Kind:

10,46,72,80,100,126,152,162,207,215,251,261,270,279,297,334,406,433,460,
521,604,612,621,640,702,720,792,800,927,972,...

(A non-cube number is a pseudo-cube of second kind if some permutation of the digits is a cube.)

One listed all pseudo-cubes of second kind up to 1000.

91) Pseudo-Cubes of Third Kind:

10,46,72,80,100,125,126,152,162,207,215,251,261,270,279,297,334,343, 406,433,460,512,521,604,612,621,640,702,720,792,800,927,972,1000,...

(A number is a pseudo-cube of third kind if some nontrivial permutation of the digits is a cube.)

Question: How many pseudo-cubes of third kind are cubes? (We conjecture: an infinity).

(There are cubes which are not pseudo-cubes of third kind, and the reverse: there are pseudo-cubes of third kind which are not cubes.)

One listed all pseudo-cubes of third kind up to 1000.

92) Pseudo-m-Powers of First Kind:

(A number is a pseudo-m-power of first kind if some permutation of the digits is an m-power, including the identity permutation; $m \geq 2$.)

93) Pseudo-m-powers of second kind:

(A non m-power number is a pseudo-m-power of second kind if some permutation of the digits is an m-power; $m \geq 2$.)

94) Pseudo-m-Powers of Third Kind:

(A number is a pseudo-m-power of third kind if some nontrivial permutation of the digits is an m-power; $m \geq 2$.)

Question: How many pseudo-m-powers of third kind are m-power numbers? (We conjecture: an infinity).

(There are m-powers which are not pseudo-m-powers of third kind, and the reverse: there are pseudo-m-powers of third kind which are not m-powers.)

95) Pseudo-Factorials of First Kind:

1,2,6,10,20,24,42,60,100,102,120,200,201,204,207,210,240,270,402,420,600,

702,720,1000,1002,1020,1200,2000,2001,2004,2007,2010,2040,2070,2100,2400,

2700,4002,4005,4020,4050,4200,4500,5004,5040,5400,6000,7002,7020,7200,...

(A number is a pseudo-factorial of first kind if some permutation of the digits is a factorial number, including the identity permutation.)

(Of course, all factorials are pseudo-factorials of first kind, but not the reverse!)

One listed all pseudo-factorials of first kind up to 10000.

Procedure to obtain this sequence:

- calculate all factorials with one digit only (1!=1, 2!=2, and 3!=6),
  this is line_1 (of one digit pseudo-factorials):
  1,2,6;
- add 0 (zero) at the end of each element of line_1,
  calculate all factorials with two digits (4!=24 only)
  and all permutations of their digits:
  this is line_2 (of two digits pseudo-factorials):
  10,20,60; 24, 42;
- add 0 (zero) at the end of each element of line_2 as well as anywhere
  in between their digits,
  calculate all factorials with three digits (5!=120, and 6!=720)
  and all permutations of their digits:
  this is line_3 (of three digits pseudo-factorials):

100,200,600,240,420,204,402; 120,720, 102,210,201,702,270,720;
and so on ...
to get from line_k to line_(k+1) do:
- add 0 (zero) at the end of each element of line_k as well as anywhere
  in between their digits,
  calculate all factorials with (k+1) digits
  and all permutations of their digits;
The set will be formed by all line_1 to the last line elements
in an increasing order.

The pseudo-factorials of second kind and third kind can
be deduced from the first kind ones..

96) Pseudo-Factorials of Second Kind:
10,20,42,60,100,102,200,201,204,207,210,240,270,402,420,600,
702,1000,1002,1020,1200,2000,2001,2004,2007,2010,2040,2070,2100,2400,
2700,4002,4005,4020,4050,4200,4500,5004,5400,6000,7002,7020,7200,...
(A non-factorial number is a pseudo-factorial of second kind if
some permutation of the digits is a factorial number.)

97) Pseudo-Factorials of Third Kind:
10,20,42,60,100,102,200,201,204,207,210,240,270,402,420,600,
702,1000,1002,1020,1200,2000,2001,2004,2007,2010,2040,2070,2100,2400,
2700,4002,4005,4020,4050,4200,4500,5004,5400,6000,7002,7020,7200,...
(A number is a pseudo-factorial of third kind if some nontrivial
permutation of the digits is a factorial number.)

Question: How many pseudo-factorials of third kind are
factorial numbers? (We conjectured: none! ... that means the
pseudo-factorials of second kind set and pseudo-factorials of
third kind set coincide!).

(Unfortunately, the second and third kinds of pseudo-factorials
coincide.)

98) Pseudo-Divisors of First Kind:
1,10,100,1,2,10,20,100,200,1,3,10,30,100,300,1,2,4,10,20,40,100,200,400,
1,5,10,50,100,500,1,2,3,6,10,20,30,60,100,200,300,600,1,7,10,70,100,700,
1,2,4,8,10,20,40,80,100,200,400,800,1,3,9,10,30,90,100,300,900,1,2,5,10,
20,50,100,200,500,1000,...
(The pseudo-divisors of first kind of n.)

(A number is a pseudo-divisor of first kind of n if
some permutation of the digits is a divisor of n, including the
identity permutation.)

(Of course, all divisors are pseudo-divisors of first kind,
but not the reverse!)

A strange property: any integer has an infinity of pseudo-divisors of first kind !!
because 10...0 becomes 0...01 = 1, by a circular permutation of its digits,
and 1 divides any integer !

One listed all pseudo-divisors of first kind up to 1000
for the numbers 1, 2, 3, ..., 10.

Procedure to obtain this sequence:
- calculate all divisors with one digit only,
  this is line_1 (of one digit pseudo-divisors);
- add 0 (zero) at the end of each element of line_1,
  calculate all divisors with two digits
  and all permutations of their digits:
  this is line_2 (of two digits pseudo-divisors);
- add 0 (zero) at the end of each element of line_2 as well as anywhere
  in between their digits,
  calculate all divisors with three digits
  and all permutations of their digits:
  this is line_3 (of three digits pseudo-divisors);

and so on ...
to get from line_k to line_(k+1) do:
- add 0 (zero) at the end of each element of line_k as well as anywhere
  in between their digits,
  calculate all divisors with (k+1) digits
  and all permutations of their digits;

The set will be formed by all line_1 to the last line elements
in an increasing order.

The pseudo-divisors of second kind and third kind can
be deduced from the first kind ones.

99) Pseudo-Divisors of Second Kind:
10,100,10,20,100,200,10,30,100,300,10,20,40,100,200,400,10,50,100,500,10,
20,30,60,100,200,300,600,10,70,100,700,10,20,40,80,100,200,400,800,10,30,

90,100,300,900,20,50,100,200,500,1000,...
(The pseudo-divisors of second kind of n.)

(A non-divisor of n is a pseudo-divisor of second kind of n
if some permutation of the digits is a divisor of n.)

100) Pseudo-Divisors of Third Kind:
10,100,10,20,100,200,10,30,100,300,10,20,40,100,200,400,10,50,100,500,10,
20,30,60,100,200,300,600,10,70,100,700,10,20,40,80,100,200,400,800,10,30,
90,100,300,900,10,20,50,100,200,500,1000,...
(The pseudo-divisors of third kind of n.)

(A number is a pseudo-divisor of third kind of n if some
nontrivial permutation of the digits is a divisor of n.)

A strange property: any integer has an infinity of pseudo-divisors of third kind !!
because 10...0 becomes 0...01 = 1, by a circular permutation of its digits,
and 1 divides any integer !

There are divisors of n which are not pseudo-divisors of
third kind of n,
and the reverse:
there are pseudo-divisors of third kind of n which are not
divisors of n.

101) Pseudo-Odd Numbers of First Kind:
1,3,5,7,9,10,11,12,13,14,15,16,17,18,19,21,23,25,27,29,30,31,32,33,34,35,
36,37,38,39,41,43,45,47,49,50,51,52,53,54,55,56,57,58,59,61,63,65,67,69,70,
71,72,73,74,75,76,...
(Some permutation of digits is an odd number.)

102) Pseudo-odd Numbers of Second Kind:
10,12,14,16,18,30,32,34,36,38,50,52,54,56,58,70,72,74,76,78,90,92,94,96,98,
100,102,104,106,108,110,112,114,116,118,...
(Even numbers such that some permutation of digits is an odd number.)

103) Pseudo-Odd Numbers of Third Kind:
10,11,12,13,14,15,16,17,18,19,30,31,32,33,34,35,36,37,38,39,50,51,52,53,54,
55,56,57,58,59,70,71,72,73,74,75,76,...

(Nontrivial permutation of digits is an odd number.)

104) Pseudo-Triangular Numbers:
1,3,6,10,12,15,19,21,28,30,36,45,54,55,60,61,63,66,78,82,87,91,...
(Some permutation of digits is a triangular number.)

A triangular number has the general form: $\frac{n(n+1)}{2}$.

105) Pseudo-Even Numbers of First Kind:
0,2,4,6,8,10,12,14,16,18,20,21,22,23,24,25,26,27,28,29,30,32,34,36,38,40,
41,42,43,44,45,46,47,48,49,50,52,54,56,58,60,61,62,63,64,65,66,67,68,69,70,
72,74,76,78,80,81,82,83,84,85,86,87,88,89,90,92,94,96,98,100,...
(The pseudo-even numbers of first kind.)

(A number is a pseudo-even number of first kind if some permutation of the digits is a even number, including the identity permutation.)

(Of course, all even numbers are pseudo-even numbers of first kind, but not the reverse!)

A strange property: an odd number can be a pseudo-even number!

One listed all pseudo-even numbers of first kind up to 100.

106) Pseudo-Even Numbers of Second Kind:
21,23,25,27,29,41,43,45,47,49,61,63,65,67,69,81,83,85,87,89,101,103,105,
107,109,121,123,125,127,129,141,143,145,147,149,161,163,165,167,169,181,
183,185,187,189,201,...
(The pseudo-even numbers of second kind.)

(A non-even number is a pseudo-even number of second kind if some permutation of the digits is a even number.)

107) Pseudo-Even Numbers of Third Kind:
20,21,22,23,24,25,26,27,28,29,40,41,42,43,44,45,46,47,48,49,60,61,62,63,64,
65,66,67,68,69,80,81,82,83,84,85,86,87,88,89,100,101,102,103.104,105,106,
107,108,109,110,120,121,122,123,124,125,126,127,128,129,130,...

(The pseudo-even numbers of third kind.)

(A number is a pseudo-even number of third kind if some nontrivial permutation of the digits is a even number.)

108) Pseudo-Multiples of First Kind (of 5):
0,5,10,15,20,25,30,35,40,45,50,51,52,53,54,55,56,57,58,59,60,65,70,75,80, 85,90,95,100,101,102,103,104,105,106,107,108,109,110,115,120,125,130,135, 140,145,150,151,152,153,154,155,156,157,158,159,160,165,...
(The pseudo-multiples of first kind of 5.)

(A number is a pseudo-multiple of first kind of 5 if some permutation of the digits is a multiple of 5, including the identity permutation.)

(Of course, all multiples of 5 are pseudo-multiples of first kind, but not the reverse!)

109) Pseudo-Multiples of Second Kind (of 5):
51,52,53,54,56,57,58,59,101,102,103,104,106,107,108,109,151,152,153,154, 156,157,158,159,201,202,203,204,206,207,208,209,251,252,253,254,256,257, 258,259,301,302,303,304,306,307,308,309,351,352...
(The pseudo-multiples of second kind of 5.)

(A non-multiple of 5 is a pseudo-multiple of second kind of 5 if some permutation of the digits is a multiple of 5.)

110) Pseudo-Multiples of Third kind (of 5):
50,51,52,53,54,55,56,57,58,59,100,101,102,103,104,105,106,107,108,109,110, 115,120,125,130,135,140,145,150,151,152,153,154,155,156,157,158,159,160, 165,170,175,180,185,190,195,200,...
(The pseudo-multiples of third kind of 5.)

(A number is a pseudo-multiple of third kind of 5 if some nontrivial permutation of the digits is a multiple of 5.)

111) Pseudo-Multiples of first kind of p (p is an integer $\geq 2$) {Generalizations}:
(The pseudo-multiples of first kind of p.)

(A number is a pseudo-multiple of first kind of p if
some permutation of the digits is a multiple of p, including the identity permutation.)

(Of course, all multiples of p are pseudo-multiples of first kind, but not the reverse!)

Procedure to obtain this sequence:

- calculate all multiples of p with one digit only (if any), this is line_1 (of one digit pseudo-multiples of p);
- add 0 (zero) at the end of each element of line_1, calculate all multiples of p with two digits (if any) and all permutations of their digits: this is line_2 (of two digits pseudo-multiples of p);
- add 0 (zero) at the end of each element of line_2 as well as anywhere in between their digits, calculate all multiples with three digits (if any) and all permutations of their digits: this is line_3 (of three digits pseudo-multiples of p);

and so on ...
to get from line_k to line_(k+1) do:

- add 0 (zero) at the end of each element of line_k as well as anywhere in between their digits, calculate all multiples with (k+1) digits (if any) and all permutations of their digits;

The set will be formed by all line_1 to the last line elements in an increasing order.

The pseudo-multiples of second kind and third kind of p can be deduced from the first kind ones.

112) Pseudo-Multiples of Second Kind of p (p is an integer $\geq 2$):
(The pseudo-multiples of second kind of p.)

(A non-multiple of p is a pseudo-multiple of second kind of p if some permutation of the digits is a multiple of p.)

113) Pseudo-multiples of third kind of p (p is an integer $\geq 2$):
(The pseudo-multiples of third kind of p.)

(A number is a pseudo-multiple of third kind of p if some nontrivial permutation of the digits is a multiple of p.)

114) Constructive Set (of digits 1,2):
1,2,11,12,21,22,111,112,121,122,211,212,221,222,1111,1112,1121,1122,1211,
1212,1221,1222,21112112,2121,2122,2211,2212,2221,2222,...
(Numbers formed by digits 1 and 2 only.)

Definition:
a1) 1, 2 belong to S;

a2) if a, b belong to S, then $\overline{ab}$ belongs to S too;

a3) only elements obtained by rules a1) and a2) applied a finite number of times belong to S.

Remark:
- there are $2^k$ numbers of k digits in the sequence, for k = 1, 2, 3, ... ;
- to obtain from the k-digits number group the (k+1)-digits number group, just put first the digit 1 and second the digit 2 in the front of all k-digits numbers.

115) Constructive Set (of digits 1,2,3):
1,2,3,11,12,13,21,22,23,31,32,33,111,112,113,121,122,123,131,132,133,211,
212,213,221,222,223,231,232,233,311,312,313,321,322,323,331,332,333,...
(Numbers formed by digits 1, 2, and 3 only.)

Definition:
a1) 1, 2, 3 belong to S;

a2) if a, b belong to S, then $\overline{ab}$ belongs to S too;

a3) only elements obtained by rules a1) and a2) applied a finite number of times belong to S.

Remark:

- there are $3^k$ numbers of k digits in the sequence, for k = 1, 2, 3, ... ;
- to obtain from the k-digits number group the (k+1)-digits number group, just put first the digit 1, second the digit 2, and third the digit 3 in the front of all k-digits numbers.

116) Generalized constructive set:

(Numbers formed by digits $d_1$ , $d_2$, ..., $d_m$ only, all $d_i$ being different each other, $1 \leq m \leq 9$.)

Definition:

a1) $d_1$ , $d_2$, ..., $d_m$ belong to S;

a2) if a, b belong to S, then $\overline{ab}$ belongs to S too;

a3) only elements obtained by rules a1) and a2) applied a finite number of times belong to S.

Remark:

- there are $m^k$ numbers of k digits in the sequence, for k = 1, 2, 3, ... ;
- to obtain from the k-digits number group the (k+1)-digits number group, just put first the digit $d_1$ , second the digit $d_2$, ..., and the m-th time digit $d_m$ in the front of all k-digits numbers.

More general: all digits $d_i$ can be replaced by numbers as large as we want (therefore of many digits each), and also m can be as large as we want.

117) Square Roots:

0,1,1,1,2,2,2,2,2,3,3,3,3,3,3,3,4,4,4,4,4,4,4,4,4,5,5,5,5,5,5,5,5,5,5,5,
6,6,6,6,6,6,6,6,6,6,6,6,6,7,7,7,7,7,7,7,7,7,7,7,7,7,7,7,8,8,8,8,8,8,8,8,8,
8,8,8,8,8,8,8,8,9,9,9,9,9,9,9,9,9,9,9,9,9,9,9,9,9,9,9,10,10,10,10,10,10,10,
10,10,10,10,10,10,10,10,10,10,10,10,10,10,...

( $s_q(n)$ is the superior integer part of square root of n.)

Remark: this sequence is the natural sequence, where each number is repeated 2n+1 times,
because between $n^2$ (included) and $(n+1)^2$ (excluded) there are $(n+1)^2 - n^2$ different numbers.

118) Cubical Roots:

0,1,1,1,1,1,1,1,2,2,2,2,2,2,2,2,2,2,2,2,2,2,2,2,2,2,2,3,3,3,3,3,3,3,3,3,3,
3,3,3,3,3,3,3,3,3,3,3,3,3,3,3,3,3,3,3,3,3,3,3,3,3,3,3,4,4,4,4,4,4,4,4,4,4,
4,4,4,4,4,4,4,4,4,4,4,4,4,4,4,4,4,4,4,4,4,4,4,4,4,4,4,4,4,4,4,4,4,4,4,4,4,
4,4,4,4,4,4,4,4,4,4,4,4,4,4,...

( $c_q(n)$ is the superior integer part of cubical root of n.)

Remark: this sequence is the natural sequence, where each number is repeated $3n^2 + 3n + 1$ times,
because between $n^3$ (included) and $(n+1)^3$ (excluded) there are $(n+1)^3 - n^3$ different numbers.

119) m-Power Roots:

( $m_q(n)$ is the superior integer part of m-power root of n.)

Remark: this sequence is the natural sequence, where each number is repeated $(n+1)^m - n^m$ times.

120) Numerical Carpet:

has the general form

.
.
.
1
1a1
1aba1
1abcba1
1abcdcba1
1abcdedcba1
1abcdefedcba1
...1abcdefgfedcba1...
1abcdefedcba1
1abcdedcba1
1abcdcba1
1abcba1
1aba1
1a1
1
.
.
.

On the border of level 0, the elements are equal to "1";
they form a rhomb.
Next, on the border of level 1, the elements are equal to "a",
where "a" is the sum of all elements of the previous border;
the "a"s form a rhombus too inside the previous one.
Next again, on the border of level 2, the elements are equal to "b",
where "b" is the sum of all elements of the previous border;
the "b"s form a rhombus too inside the previous one.
And so on...
The carpet is symmetric and esthetic, in its middle g is the
sum of all carpet numbers (the core).

Look at a few terms of the Numerical Carpet:

1

1
141
1

1
1 8 1
1 8 40 8 1
1 8 1
1

1
1 12 1
1 12 108 12 1
1 12 108 540 108 12 1
1 12 108 12 1
1 12 1
1

1
1 16 1
1 16 208 16 1
1 16 208 1872 208 16 1
1 16 208 1872 9360 1872 208 16 1
1 16 208 1872 208 16 1
1 16 208 16 1
1 16 1
1

1
1 20 1
1 20 340 20 1
1 20 340 4420 340 20 1
1 20 340 4420 39780 4420 340 20 1
1 20 340 4420 39780 198900 39780 4420 340 20 1
1 20 340 4420 39780 4420 340 20 1
1 20 340 4420 340 20 1
1 20 340 20 1
1 20 1
1
.
.

Or, under other form:

| | | | | | | | | |
|---|---|---|---|---|---|---|---|---|
| 1 | | | | | | | | |
| 1 | 4 | | | | | | | |
| 1 | 8 | 40 | | | | | | |
| 1 | 12 | 108 | 504 | | | | | |
| 1 | 16 | 208 | 1872 | 9360 | | | | |
| 1 | 20 | 340 | 4420 | 39780 | 198900 | | | |
| 1 | 24 | 504 | 8568 | 111384 | 1002456 | 5012280 | | |
| 1 | 28 | 700 | 14700 | 249900 | 3248700 | 29238300 | 146191500 | |
| 1 | 32 | 928 | 23200 | 487200 | 8282400 | 107671200 | 969040800 | 4845204000 |

..................................................................

.

.

.

General Formula:

$$C(n,k) = 4n\prod_{i=1}^{K}(4n - 4i + 1), \text{ for } 1 \le k \le n,$$

and C(n,0) = 1.

121) Goldbach-Smarandache table:
6,10,14,18,26,30,38,42,42,54,62,74,74,90,...
( t(n) is the largest even number such that any other even number not exceeding it is the sum of two of the first n odd primes.)

It helps to better understand Goldbach's conjecture:
- if t(n) is unlimited, then the conjecture is true;
- if t(n) is constant after a certain rank, then the conjecture is false.

Also, the table gives how many times an even number is written as a sum of two odd primes, and in what combinations -- which can be found in the "Encyclopedia of Integer Sequences" by N. J. A. Sloane and S. Plouffe,

Of course, $t(n) \leq 2p_n$ , where $p_n$ is the n-th odd prime, n = 1, 2, 3, ... .

Here is the table:

| + | 3 | 5 | 7 | 11 | 13 | 17 | 19 | 23 | 29 | 31 | 37 | 41 | 43 | 47 |
|---|---|---|---|---|---|---|---|---|---|---|---|---|---|---|
| 3 | 6 | 8 | 10 | 14 | 16 | 20 | 22 | 26 | 32 | 34 | 40 | 44 | 46 | 50 . |
| 5 | | 10 | 12 | 16 | 18 | 22 | 24 | 28 | 34 | 36 | 42 | 46 | 48 | 52 . |
| 7 | | | 14 | 18 | 20 | 24 | 26 | 30 | 36 | 38 | 44 | 48 | 50 | 54 . |
| 11 | | | | 22 | 24 | 28 | 30 | 34 | 40 | 42 | 48 | 52 | 54 | 58 . |
| 13 | | | | | 26 | 30 | 32 | 36 | 42 | 44 | 50 | 54 | 56 | 60 . |
| 17 | | | | | | 34 | 36 | 40 | 46 | 48 | 54 | 58 | 60 | 64 . |
| 19 | | | | | | | 38 | 42 | 48 | 50 | 56 | 60 | 62 | 66 . |
| 23 | | | | | | | | 46 | 52 | 54 | 60 | 64 | 66 | 70 . |
| 29 | | | | | | | | | 58 | 60 | 66 | 70 | 72 | 76 . |
| 31 | | | | | | | | | | 62 | 68 | 72 | 74 | 78 . |
| 37 | | | | | | | | | | | 74 | 78 | 80 | 84 . |
| 41 | | | | | | | | | | | | 82 | 84 | 88 . |
| 43 | | | | | | | | | | | | | 86 | 90 . |
| 47 | | | | | | | | | | | | | | 94 . |

. .
. .
. .

122) Smarandache-Vinogradov table:

9,15,21,29,39,47,57,65,71,93,99,115,129,137,...

(v(n) is the largest odd number such that any odd number $\geq 9$ not exceeding it is the sum of three of the first n odd primes.)

It helps to better understand Goldbach's conjecture for three primes:

- if v(n) is unlimited, then the conjecture is true;
- if v(n) is constant after a certain rank, then the conjecture is false.

(Vinogradov proved in 1937 that any odd number greater than $3^{3^{15}}$

satisfies this conjecture.

But what about values less than $3^{3^{15}}$ ?)

Also, the table gives you in how many different combinations an odd number is written as a sum of three odd primes, and in what combinations.

Of course, $v(n) \leq 3p_n$, where $p_n$ is the n-th odd prime, n = 1, 2, 3, ... .
It is also generalized for the sum of m primes,
and how many times a number is written as a sum of m primes (m > 2).

This is a 3-dimensional 14x14x14 table, that we can expose only as 14 planar 14x14 tables (using Goldbach-Smarandache table):

| 3 + | 3 | 5 | 7 | 11 | 13 | 17 | 19 | 23 | 29 | 31 | 37 | 41 | 43 | 47 . |
|---|---|---|---|---|---|---|---|---|---|---|---|---|---|---|
| 3 | 9 | 11 | 13 | 17 | 19 | 23 | 25 | 29 | 35 | 37 | 43 | 47 | 49 | 53 . |
| 5 | | 13 | 15 | 19 | 21 | 25 | 27 | 31 | 37 | 39 | 45 | 49 | 51 | 55 . |
| 7 | | | 17 | 21 | 23 | 27 | 29 | 33 | 39 | 41 | 47 | 51 | 53 | 57 . |
| 11 | | | | 25 | 27 | 31 | 33 | 37 | 43 | 45 | 51 | 55 | 57 | 61 . |
| 13 | | | | | 29 | 33 | 35 | 39 | 45 | 47 | 53 | 57 | 59 | 63 . |
| 17 | | | | | | 37 | 39 | 43 | 49 | 51 | 57 | 61 | 63 | 67 . |
| 19 | | | | | | | 41 | 45 | 51 | 53 | 59 | 63 | 65 | 69 . |
| 23 | | | | | | | | 49 | 55 | 57 | 63 | 67 | 69 | 73 . |
| 29 | | | | | | | | | 61 | 63 | 69 | 73 | 75 | 79 . |
| 31 | | | | | | | | | | 65 | 71 | 75 | 77 | 81 . |
| 37 | | | | | | | | | | | 77 | 81 | 83 | 87 . |
| 41 | | | | | | | | | | | | 85 | 87 | 91 . |
| 43 | | | | | | | | | | | | | 89 | 93 . |
| 47 | | | | | | | | | | | | | | 97 . |

.

.

.

| 5 + | 3 | 5 | 7 | 11 | 13 | 17 | 19 | 23 | 29 | 31 | 37 | 41 | 43 | 47 . |
|---|---|---|---|---|---|---|---|---|---|---|---|---|---|---|
| 3 | 11 | 13 | 15 | 19 | 21 | 25 | 27 | 31 | 37 | 39 | 45 | 49 | 51 | 55 . |
| 5 | | 15 | 17 | 21 | 23 | 27 | 29 | 33 | 39 | 41 | 47 | 51 | 53 | 57 . |
| 7 | | | 19 | 23 | 25 | 29 | 31 | 35 | 41 | 43 | 49 | 53 | 55 | 59 . |
| 11 | | | | 27 | 29 | 33 | 35 | 39 | 45 | 47 | 53 | 57 | 59 | 63 . |
| 13 | | | | | 31 | 35 | 37 | 41 | 47 | 49 | 55 | 59 | 61 | 65 . |
| 17 | | | | | | 39 | 41 | 45 | 51 | 53 | 59 | 63 | 65 | 69 . |
| 19 | | | | | | | 43 | 47 | 53 | 55 | 61 | 65 | 67 | 71 . |
| 23 | | | | | | | | 51 | 57 | 59 | 65 | 69 | 71 | 75 . |
| 29 | | | | | | | | | 63 | 65 | 71 | 75 | 77 | 81 . |
| 31 | | | | | | | | | | 67 | 73 | 77 | 79 | 83 . |
| 37 | | | | | | | | | | | 79 | 83 | 85 | 89 . |
| 41 | | | | | | | | | | | | 87 | 89 | 93 . |
| 43 | | | | | | | | | | | | | 91 | 95 . |
| 47 | | | | | | | | | | | | | | 99 . |

.

.

.

| 7 + | 3 | 5 | 7 | 11 | 13 | 17 | 19 | 23 | 29 | 31 | 37 | 41 | 43 | 47 . |
|---|---|---|---|---|---|---|---|---|---|---|---|---|---|---|
| 3 | 13 | 15 | 17 | 21 | 23 | 27 | 29 | 33 | 39 | 41 | 47 | 51 | 53 | 57 . |
| 5 | | 17 | 19 | 23 | 25 | 29 | 31 | 35 | 41 | 43 | 49 | 53 | 55 | 59 . |
| 7 | | | 21 | 25 | 27 | 31 | 33 | 37 | 43 | 45 | 51 | 55 | 57 | 61 . |
| 11 | | | | 29 | 31 | 35 | 37 | 41 | 47 | 49 | 55 | 59 | 61 | 65 . |
| 13 | | | | | 33 | 37 | 39 | 43 | 49 | 51 | 57 | 61 | 63 | 67 . |
| 17 | | | | | | 41 | 43 | 47 | 53 | 55 | 61 | 65 | 67 | 71 . |
| 19 | | | | | | | 45 | 49 | 55 | 57 | 63 | 67 | 69 | 73 . |
| 23 | | | | | | | | 53 | 59 | 61 | 67 | 71 | 73 | 77 . |
| 29 | | | | | | | | | 65 | 67 | 73 | 77 | 79 | 83 . |
| 31 | | | | | | | | | | 69 | 75 | 79 | 81 | 85 . |
| 37 | | | | | | | | | | | 81 | 85 | 87 | 91 . |
| 41 | | | | | | | | | | | | 89 | 91 | 95 . |
| 43 | | | | | | | | | | | | | 93 | 97 . |
| 47 | | | | | | | | | | | | | | 101 . |

.
.
.

| 11 + | 3 | 5 | 7 | 11 | 13 | 17 | 19 | 23 | 29 | 31 | 37 | 41 | 43 | 47 . |
|---|---|---|---|---|---|---|---|---|---|---|---|---|---|---|
| 3 | 17 | 19 | 21 | 25 | 27 | 31 | 33 | 37 | 43 | 45 | 51 | 55 | 57 | 61 . |
| 5 | | 21 | 23 | 27 | 29 | 33 | 35 | 39 | 45 | 47 | 53 | 57 | 59 | 63 . |
| 7 | | | 25 | 29 | 31 | 35 | 37 | 41 | 47 | 49 | 55 | 59 | 61 | 65 . |
| 11 | | | | 33 | 35 | 39 | 41 | 45 | 51 | 53 | 59 | 63 | 65 | 69 . |
| 13 | | | | | 37 | 41 | 43 | 47 | 53 | 55 | 61 | 65 | 67 | 71 . |
| 17 | | | | | | 45 | 47 | 51 | 57 | 59 | 65 | 69 | 71 | 75 . |
| 19 | | | | | | | 49 | 53 | 59 | 61 | 67 | 71 | 73 | 77 . |
| 23 | | | | | | | | 57 | 63 | 65 | 71 | 75 | 77 | 81 . |
| 29 | | | | | | | | | 69 | 71 | 77 | 81 | 83 | 87 . |
| 31 | | | | | | | | | | 73 | 79 | 83 | 85 | 89 . |
| 37 | | | | | | | | | | | 85 | 89 | 91 | 95 . |
| 41 | | | | | | | | | | | | 93 | 95 | 99 . |
| 43 | | | | | | | | | | | | | 97 | 101. |
| 47 | | | | | | | | | | | | | | 105. |

.
.
.

| 13 + | 3 | 5 | 7 | 11 | 13 | 17 | 19 | 23 | 29 | 31 | 37 | 41 | 43 | 47 . |
|---|---|---|---|---|---|---|---|---|---|---|---|---|---|---|
| 3 | 19 | 21 | 23 | 27 | 29 | 33 | 35 | 39 | 45 | 47 | 53 | 57 | 59 | 63. |
| 5 | | 23 | 25 | 29 | 31 | 35 | 37 | 41 | 47 | 49 | 55 | 59 | 61 | 65 . |
| 7 | | | 27 | 31 | 33 | 37 | 39 | 43 | 49 | 51 | 57 | 61 | 63 | 67 . |
| 11 | | | | 35 | 37 | 41 | 43 | 47 | 53 | 55 | 61 | 65 | 67 | 71 . |
| 13 | | | | | 39 | 43 | 45 | 49 | 55 | 57 | 63 | 67 | 69 | 73 . |
| 17 | | | | | | 47 | 49 | 53 | 59 | 61 | 67 | 71 | 73 | 77 . |
| 19 | | | | | | | 51 | 55 | 61 | 63 | 69 | 73 | 75 | 79 . |
| 23 | | | | | | | | 59 | 65 | 67 | 73 | 77 | 79 | 83 . |
| 29 | | | | | | | | | 71 | 73 | 79 | 83 | 85 | 89 . |
| 31 | | | | | | | | | | 75 | 81 | 85 | 87 | 91 . |
| 37 | | | | | | | | | | | 87 | 91 | 93 | 97 . |
| 41 | | | | | | | | | | | | 95 | 97 | 101 . |
| 43 | | | | | | | | | | | | | 99 | 103. |
| 47 | | | | | | | | | | | | | | 107 . |

.
.
.

| 17 + | 3 | 5 | 7 | 11 | 13 | 17 | 19 | 23 | 29 | 31 | 37 | 41 | 43 | 47 . |
|---|---|---|---|---|---|---|---|---|---|---|---|---|---|---|
| 3 | 23 | 25 | 27 | 31 | 33 | 37 | 39 | 43 | 49 | 51 | 57 | 61 | 63 | 67 . |
| 5 | | 27 | 29 | 33 | 35 | 39 | 41 | 45 | 51 | 53 | 59 | 63 | 65 | 69 . |
| 7 | | | 31 | 35 | 37 | 41 | 43 | 47 | 53 | 55 | 61 | 65 | 67 | 71 . |
| 11 | | | | 39 | 41 | 45 | 47 | 51 | 57 | 59 | 65 | 69 | 71 | 75 . |
| 13 | | | | | 43 | 47 | 49 | 53 | 59 | 61 | 67 | 71 | 73 | 77 . |
| 17 | | | | | | 51 | 53 | 57 | 63 | 65 | 71 | 75 | 77 | 81 . |
| 19 | | | | | | | 55 | 59 | 65 | 67 | 73 | 77 | 79 | 83 . |
| 23 | | | | | | | | 63 | 69 | 71 | 77 | 81 | 83 | 87 . |
| 29 | | | | | | | | | 75 | 77 | 83 | 87 | 89 | 93 . |
| 31 | | | | | | | | | | 79 | 85 | 89 | 91 | 95 . |
| 37 | | | | | | | | | | | 91 | 95 | 97 | 101. |
| 41 | | | | | | | | | | | | 99 | 101 | 105 . |
| 43 | | | | | | | | | | | | | 103 | 107 . |
| 47 | | | | | | | | | | | | | | 111. |

.
.
.

| 19 + | 3 | 5 | 7 | 11 | 13 | 17 | 19 | 23 | 29 | 31 | 37 | 41 | 43 | 47 . |
|---|---|---|---|---|---|---|---|---|---|---|---|---|---|---|
| 3 | 25 | 27 | 29 | 33 | 35 | 39 | 41 | 45 | 51 | 53 | 59 | 63 | 65 | 69 . |
| 5 | | 29 | 31 | 35 | 37 | 41 | 43 | 47 | 53 | 55 | 61 | 65 | 67 | 71 . |
| 7 | | | 33 | 37 | 39 | 43 | 45 | 49 | 55 | 57 | 63 | 67 | 69 | 73 . |
| 11 | | | | 41 | 43 | 47 | 49 | 53 | 59 | 61 | 67 | 71 | 73 | 77 . |
| 13 | | | | | 45 | 49 | 51 | 55 | 61 | 63 | 69 | 73 | 75 | 79 . |
| 17 | | | | | | 53 | 55 | 59 | 65 | 67 | 73 | 77 | 79 | 83 . |
| 19 | | | | | | | 57 | 61 | 67 | 69 | 75 | 79 | 81 | 85 . |
| 23 | | | | | | | | 65 | 71 | 73 | 79 | 83 | 85 | 89 . |
| 29 | | | | | | | | | 77 | 79 | 85 | 89 | 91 | 95 . |
| 31 | | | | | | | | | | 81 | 87 | 91 | 93 | 97 . |
| 37 | | | | | | | | | | | 93 | 97 | 99 | 103. |
| 41 | | | | | | | | | | | | 101 | 103 | 107. |
| 43 | | | | | | | | | | | | | 105 | 109. |
| 47 | | | | | | | | | | | | | | 113. |

.
.
.

| 23 + | 3 | 5 | 7 | 11 | 13 | 17 | 19 | 23 | 29 | 31 | 37 | 41 | 43 | 47 . |
|---|---|---|---|---|---|---|---|---|---|---|---|---|---|---|
| 3 | 29 | 31 | 33 | 37 | 39 | 43 | 45 | 49 | 55 | 57 | 63 | 67 | 69 | 73 . |
| 5 | | 33 | 35 | 39 | 41 | 45 | 47 | 51 | 57 | 59 | 65 | 69 | 71 | 75 . |
| 7 | | | 37 | 41 | 43 | 47 | 49 | 53 | 59 | 61 | 67 | 71 | 73 | 77 . |
| 11 | | | | 45 | 47 | 51 | 53 | 57 | 63 | 65 | 71 | 75 | 77 | 81 . |
| 13 | | | | | 49 | 53 | 55 | 59 | 65 | 67 | 73 | 77 | 79 | 83 . |
| 17 | | | | | | 57 | 59 | 63 | 69 | 71 | 77 | 81 | 83 | 87 . |
| 19 | | | | | | | 61 | 65 | 71 | 73 | 79 | 83 | 85 | 89 . |
| 23 | | | | | | | | 69 | 75 | 77 | 83 | 87 | 89 | 93 . |
| 29 | | | | | | | | | 81 | 83 | 89 | 93 | 95 | 99 . |
| 31 | | | | | | | | | | 85 | 91 | 95 | 97 | 101. |
| 37 | | | | | | | | | | | 97 | 101 | 103 | 107. |
| 41 | | | | | | | | | | | | 105 | 107 | 111. |
| 43 | | | | | | | | | | | | | 109 | 113. |
| 47 | | | | | | | | | | | | | | 117. |

.
.
.

| 29 + | 3 | 5 | 7 | 11 | 13 | 17 | 19 | 23 | 29 | 31 | 37 | 41 | 43 | 47 . |
|---|---|---|---|---|---|---|---|---|---|---|---|---|---|---|
| 3 | 35 | 37 | 39 | 43 | 45 | 49 | 51 | 55 | 61 | 63 | 69 | 73 | 75 | 79 . |
| 5 | | 39 | 41 | 45 | 47 | 51 | 53 | 57 | 63 | 65 | 71 | 75 | 77 | 81 . |
| 7 | | | 43 | 47 | 49 | 53 | 55 | 59 | 65 | 67 | 73 | 77 | 79 | 83 . |
| 11 | | | | 51 | 53 | 57 | 59 | 63 | 69 | 71 | 77 | 81 | 83 | 87 . |
| 13 | | | | | 55 | 59 | 61 | 65 | 71 | 73 | 79 | 83 | 85 | 89 . |
| 17 | | | | | | 63 | 65 | 69 | 75 | 77 | 83 | 87 | 89 | 93 . |
| 19 | | | | | | | 67 | 71 | 77 | 79 | 85 | 89 | 91 | 95 . |
| 23 | | | | | | | | 75 | 81 | 83 | 89 | 93 | 95 | 99 . |
| 29 | | | | | | | | | 87 | 89 | 95 | 99 | 101 | 105 . |
| 31 | | | | | | | | | | 91 | 97 | 101 | 103 | 107. |
| 37 | | | | | | | | | | | 103 | 107 | 109 | 113 . |
| 41 | | | | | | | | | | | | 111 | 113 | 117 . |
| 43 | | | | | | | | | | | | | 115 | 119 . |
| 47 | | | | | | | | | | | | | | 123 . |

.
.
.

| 31 + | 3 | 5 | 7 | 11 | 13 | 17 | 19 | 23 | 29 | 31 | 37 | 41 | 43 | 47 . |
|---|---|---|---|---|---|---|---|---|---|---|---|---|---|---|
| 3 | 37 | 39 | 41 | 45 | 47 | 51 | 53 | 57 | 63 | 65 | 71 | 75 | 77 | 81 . |
| 5 | | 41 | 43 | 47 | 49 | 53 | 55 | 59 | 65 | 67 | 73 | 77 | 79 | 83 . |
| 7 | | | 45 | 49 | 51 | 55 | 57 | 61 | 67 | 69 | 75 | 79 | 81 | 85 . |
| 11 | | | | 53 | 55 | 59 | 61 | 65 | 71 | 73 | 79 | 83 | 85 | 89 . |
| 13 | | | | | 57 | 61 | 63 | 67 | 73 | 75 | 81 | 85 | 87 | 91 . |
| 17 | | | | | | 65 | 67 | 71 | 77 | 79 | 85 | 89 | 91 | 95 . |
| 19 | | | | | | | 69 | 73 | 79 | 81 | 87 | 91 | 93 | 97 . |
| 23 | | | | | | | | 77 | 83 | 85 | 91 | 95 | 97 | 101. |
| 29 | | | | | | | | | 89 | 91 | 97 | 101 | 103 | 107. |
| 31 | | | | | | | | | | 93 | 99 | 103 | 105 | 109. |
| 37 | | | | | | | | | | | 105 | 109 | 111 | 115. |
| 41 | | | | | | | | | | | | 113 | 115 | 119. |
| 43 | | | | | | | | | | | | | 117 | 121. |
| 47 | | | | | | | | | | | | | | 125. |

.
.
.

| 37 + | 3 | 5 | 7 | 11 | 13 | 17 | 19 | 23 | 29 | 31 | 37 | 41 | 43 | 47 . |
|---|---|---|---|---|---|---|---|---|---|---|---|---|---|---|
| 3 | 43 | 45 | 47 | 51 | 53 | 57 | 59 | 63 | 69 | 71 | 77 | 81 | 83 | 87 . |
| 5 | | 47 | 49 | 53 | 55 | 59 | 61 | 65 | 71 | 73 | 79 | 83 | 85 | 89 . |
| 7 | | | 51 | 55 | 57 | 61 | 63 | 67 | 73 | 75 | 81 | 85 | 87 | 91 . |
| 11 | | | | 59 | 61 | 65 | 67 | 71 | 77 | 79 | 85 | 89 | 91 | 95 . |
| 13 | | | | | 63 | 67 | 69 | 73 | 79 | 81 | 87 | 91 | 93 | 97 . |
| 17 | | | | | | 71 | 73 | 77 | 83 | 85 | 91 | 95 | 97 | 101. |
| 19 | | | | | | | 75 | 79 | 85 | 87 | 93 | 97 | 99 | 103. |
| 23 | | | | | | | | 83 | 89 | 91 | 97 | 101 | 103 | 107. |
| 29 | | | | | | | | | 95 | 97 | 103 | 107 | 109 | 113. |
| 31 | | | | | | | | | | 99 | 105 | 109 | 111 | 115. |
| 37 | | | | | | | | | | | 111 | 115 | 117 | 121. |
| 41 | | | | | | | | | | | | 119 | 121 | 125. |
| 43 | | | | | | | | | | | | | 123 | 127. |
| 47 | | | | | | | | | | | | | | 131. |

.
.
.

| 41 + | 3 | 5 | 7 | 11 | 13 | 17 | 19 | 23 | 29 | 31 | 37 | 41 | 43 | 47 . |
|---|---|---|---|---|---|---|---|---|---|---|---|---|---|---|
| 3 | 47 | 49 | 51 | 55 | 57 | 61 | 63 | 67 | 73 | 75 | 81 | 85 | 87 | 91 . |
| 5 | | 51 | 53 | 57 | 59 | 63 | 65 | 69 | 75 | 77 | 83 | 87 | 89 | 93 . |
| 7 | | | 55 | 59 | 61 | 65 | 67 | 71 | 77 | 79 | 85 | 89 | 91 | 95 . |
| 11 | | | | 63 | 65 | 69 | 71 | 75 | 81 | 83 | 89 | 93 | 95 | 99 . |
| 13 | | | | | 67 | 71 | 73 | 77 | 83 | 85 | 91 | 95 | 97 | 101. |
| 17 | | | | | | 75 | 77 | 81 | 87 | 89 | 95 | 99 | 101 | 105. |
| 19 | | | | | | | 79 | 83 | 89 | 91 | 97 | 101 | 103 | 107. |
| 23 | | | | | | | | 87 | 93 | 95 | 101 | 105 | 107 | 111. |
| 29 | | | | | | | | | 99 | 101 | 107 | 111 | 113 | 117. |
| 31 | | | | | | | | | | 103 | 109 | 113 | 115 | 119. |
| 37 | | | | | | | | | | | 115 | 119 | 121 | 125. |
| 41 | | | | | | | | | | | | 123 | 125 | 129. |
| 43 | | | | | | | | | | | | | 127 | 131. |
| 47 | | | | | | | | | | | | | | 135. |

.
.
.

| 43 + | 3 | 5 | 7 | 11 | 13 | 17 | 19 | 23 | 29 | 31 | 37 | 41 | 43 | 47 . |
|---|---|---|---|---|---|---|---|---|---|---|---|---|---|---|
| 3 | 49 | 51 | 53 | 57 | 59 | 63 | 65 | 69 | 75 | 77 | 83 | 87 | 89 | 93 . |
| 5 | | 53 | 55 | 59 | 61 | 65 | 67 | 71 | 77 | 79 | 85 | 89 | 91 | 95 . |
| 7 | | | 57 | 61 | 63 | 67 | 69 | 73 | 79 | 81 | 87 | 91 | 93 | 97 . |
| 11 | | | | 65 | 67 | 71 | 73 | 77 | 83 | 85 | 91 | 95 | 97 | 101. |
| 13 | | | | | 69 | 73 | 75 | 79 | 85 | 87 | 93 | 97 | 99 | 103. |
| 17 | | | | | | 77 | 79 | 83 | 89 | 91 | 97 | 101 | 103 | 107. |
| 19 | | | | | | | 81 | 85 | 91 | 93 | 99 | 103 | 105 | 109. |
| 23 | | | | | | | | 89 | 95 | 97 | 103 | 107 | 109 | 113. |
| 29 | | | | | | | | | 101 | 103 | 109 | 113 | 115 | 119. |
| 31 | | | | | | | | | | 105 | 111 | 115 | 117 | 121. |
| 37 | | | | | | | | | | | 117 | 121 | 123 | 127. |
| 41 | | | | | | | | | | | | 125 | 127 | 131. |
| 43 | | | | | | | | | | | | | 129 | 133. |
| 47 | | | | | | | | | | | | | | 137. |

.
.
.

| 47 + | 3 | 5 | 7 | 11 | 13 | 17 | 19 | 23 | 29 | 31 | 37 | 41 | 43 | 47 . |
|---|---|---|---|---|---|---|---|---|---|---|---|---|---|---|
| 3 | 53 | 55 | 57 | 61 | 63 | 67 | 69 | 73 | 79 | 81 | 87 | 91 | 93 | 97 . |
| 5 | | 57 | 59 | 63 | 65 | 69 | 71 | 75 | 81 | 83 | 89 | 93 | 95 | 99 . |
| 7 | | | 61 | 65 | 67 | 71 | 73 | 77 | 83 | 85 | 91 | 95 | 97 | 101. |
| 11 | | | | 69 | 71 | 75 | 77 | 81 | 87 | 89 | 95 | 99 | 101 | 105. |
| 13 | | | | | 73 | 77 | 79 | 83 | 89 | 91 | 97 | 101 | 103 | 107. |
| 17 | | | | | | 81 | 83 | 87 | 93 | 95 | 101 | 105 | 107 | 111. |
| 19 | | | | | | | 85 | 89 | 95 | 97 | 103 | 107 | 109 | 113. |
| 23 | | | | | | | | 93 | 99 | 101 | 107 | 111 | 113 | 117. |
| 29 | | | | | | | | | 105 | 107 | 113 | 117 | 119 | 123. |
| 31 | | | | | | | | | | 109 | 115 | 119 | 121 | 125. |
| 37 | | | | | | | | | | | 121 | 125 | 127 | 131. |
| 41 | | | | | | | | | | | | 129 | 131 | 135. |
| 43 | | | | | | | | | | | | | 133 | 137. |
| 47 | | | | | | | | | | | | | | 141. |

.
.
.

123) Smarandache-Vinogradov sequence:

0,0,0,0,1,2,4,4,6,7,9,10,11,15,17,16,19,19,23,25,26,26,28,33,32,35,43,39,
40,43,43,...

(a(2k+1) represents the number of different combinations such that 2k+1 is written as a sum of three odd primes.)

This sequence is deduced from the Smarandache-Vinogradov table.

124) Smarandache Paradoxist Numbers:

There exist a few "Smarandache" number sequences.
A number n is called a "Smarandache paradoxist number" if and only if n doesn't belong to any of the Smarandache defined numbers.
Question: find the Smarandache paradoxist number sequence.

Solution:

If a number k is a Smarandache paradoxist number, then k doesn't belong to any of the Smarandache defined numbers,
therefore k doesn't belong to the Smarandache paradoxist numbers too!

If a number k doesn't belong to any of the Smarandache defined numbers, then k is a Smarandache paradoxist number,
therefore k belongs to a Smarandache defined numbers (because Smarandache paradoxist numbers is also in the same category) - contradiction.

Dilemma: Is the Smarandache paradoxist number sequence empty ?

125) Non-Smarandache numbers:

A number n is called a "non-Smarandache number" if and only if n is neither a Smarandache paradoxist number nor any of the Smarandache-defined numbers.

Question: find the non-Smarandache number sequence.

Dilemma 1: is the non-Smarandache number sequence empty, too?
Dilemma 2: is a non-Smarandache number equivalent to a Smarandache paradoxist number?? (this would be another paradox !! ... because a non-Smarandache number is not a Smarandache paradoxist number).

126) The paradox of Smarandache numbers:
Any number is a Smarandache number, the non-Smarandache number too.
(This is deduced from the following paradox (see the reference):
"All is possible, the impossible too!")

127) Romanian Multiplication:
Another algorithm to multiply two integer numbers, A and B:
- let k be an integer $\geq 2$;
- write A and B on two different vertical columns: c(A), respectively c(B);
- multiply A by k, and write the product $A_1$ on the column c(A);
- divide B by k, and write the integer part of the quotient $B_1$ on the column c(B);

... and so on with the new numbers $A_1$ and $B_1$ ,
until we get a $B_i < k$ on the column c(B);
Then:
- write another column c(r), on the right side of c(B), such that: for each number of column c(B), which may be a multiple of k plus the rest r (where r = 0, 1, 2, ..., k-1), the corresponding number on c(r) will be r;
- multiply each number of column A by its corresponding r of c(r), and put the new products on another column c(P) on the right side of c(r);
- finally add all numbers of column c(P).

$A \times B$ = the sum of all numbers of c(P).

Remark that any multiplication of integer numbers can be done only by multiplication with 2, 3, ..., k, divisions by k, and additions.

This is a generalization of Russian multiplication (when k = 2); we call it Romanian Multiplication.

This special multiplication is useful when k is very small, the best values being for k = 2 (Russian multiplication - known since Egyptian time), or k = 3. If k is greater than or equal to min {10, B}, this multiplication is trivial (the obvious multiplication).

Example 1 (if we choose k=3):

73×97= ?

| x3 | /3 | | |
|---|---|---|---|
| c(A) | c(B) | c(r) | c(P) |
| 73 | 97 | 1 | 73 |
| 219 | 32 | 2 | 438 |
| 657 | 10 | 1 | 657 |
| 1971 | 3 | 0 | 0 |
| 5913 | 1 | 1 | 5913 |

7081 total

therefore: 73x97=7081.

Remark that any multiplication of integer numbers can be done only by multiplication with 2, 3, divisions by 3, and additions.

Example 2 (if we choose k = 4):

73x97= ?

| X4 | /4 | | |
|---|---|---|---|
| c(A) | c(B) | c(r) | c(P) |
| 73 | 97 | 1 | 73 |
| 292 | 24 | 0 | 0 |
| 1168 | 6 | 2 | 2336 |
| 4672 | 1 | 1 | 4672 |

7081 total

therefore: 73x97=7081.

Remark that any multiplication of integer numbers can be done only by multiplication with 2, 3, 4, divisions by 4, and additions.

Example 3 (if we choose k = 5):

73x97= ?

| X5 | /5 | | |
|---|---|---|---|

| c(A) | c(B) | c(r) | c(P) |
|---|---|---|---|
| 73 | 97 | 2 | 146 |
| 365 | 19 | 4 | 1460 |
| 1825 | 3 | 3 | 5475 |

7081 total

therefore: 73x97=7081.
Remark that any multiplication of integer numbers can be done only by multiplication with 2, 3, 4, 5, divisions by 5, and additions.

This special multiplication becomes less useful when k increases.
Look at another example (4), what happens when k = 10:
73x97= ?

| x10 | /10 | | |
|---|---|---|---|
| c(A) | c(B) | c(r) | c(P) |
| 73 | 97 | 7 | 511(=73x7) |
| 730 | 9 | 9 | 6570 (=730x9) |

7081 total

therefore: 73x97=7081.
Remark that any multiplication of integer numbers can be done only by multiplication with 2, 3, ..., 9, 10, divisions by 10, and additions --
hence we obtain just the obvious multiplication!

128) Division by $k^n$:
Another algorithm to divide an integer number A by $k^n$, where k, n are integers ≥ 2:
- write A and $k^n$ on two different vertical columns: c(A), respectively c($k^n$);
- divide A by k, and write the integer quotient $A_1$ on the column c(A);
- divide $k^n$ by k, and write the quotient $q_1 = k^{n-1}$ on the column c($k^n$);

... and so on with the new numbers $A_1$ and $q_1$,
until we get $q_n = 1$ ($= k^0$) on the column $c(k^n)$;
Then:
- write another column c(r), on the left side of c(A), such that:
  for each number of column c(A), which may be a multiple of k plus the rest r (where r = 0, 1, 2, ..., k-1),
  the corresponding number on c(r) will be r;
- write another column c(P), on the left side of c(r), in the following way: the element on line i (except the last line which is 0) will be $k^{n-1}$;
- multiply each number of column c(P) by its corresponding r of c(r), and put the new products on another column c(R) on the left side of c(P);
- finally add all numbers of column c(R) to get the final rest $R_n$,

while the final quotient will be stated in front of $c(k^n)$'s 1.
Therefore:

$$\frac{A}{k^n} = A_n \text{ and rest } R_n.$$

Remark that any division of an integer number by $k^n$ can be done only by divisions to k, calculations of powers of k, multiplications with 1, 2, ..., k-1, additions.

This special division is useful when k is small, the best values being when k is an one-digit number, and n large.
If k is very big and n very small, this division becomes useless.

Example 1 :
$1357/(2^7) = ?$

| | | | | | |
|---|---|---|---|---|---|
| | | | /2 | /2 | |
| c(R) | c(P) | c(r) | c(A) | $c(2^7)$ | |
| 1 | $2^0$ | 1 | 1357 | $2^7$ | line_1 |
| 0 | $2^1$ | 0 | 678 | $2^6$ | line_2 |
| 4 | $2^2$ | 1 | 339 | $2^5$ | line_3 |
| 8 | $2^3$ | 1 | 169 | $2^4$ | line_4 |
| 0 | $2^4$ | 0 | 84 | $2^3$ | line_5 |
| 0 | $2^5$ | 0 | 42 | $2^2$ | line_6 |
| 64 | $2^6$ | 1 | 21 | $2^1$ | line_7 |
| | | | 10 | $2^0$ | Last_line |
| 77 | | | | | |

Therefore: $1357/(2^7) = 10$ and rest 77.

Remark that the division of an integer number by any power of 2 can be done only by divisions to 2, calculations of powers of 2, multiplications and additions.

Example 2 :
$19495/(3^8) = ?$

| | | | /3 | /3 | |
|---|---|---|---|---|---|
| c(R) | c(P) | c(r) | c(A) | c($3^8$) | |
| 1 | $3^0$ | 1 | 19495 | $3^8$ | line_1 |
| 0 | $3^1$ | 0 | 6498 | $3^7$ | line_2 |
| 0 | $3^2$ | 0 | 2166 | $3^6$ | line_3 |
| 54 | $3^3$ | 2 | 722 | $3^5$ | line_4 |
| 0 | $3^4$ | 0 | 240 | $3^4$ | line_5 |
| 486 | $3^5$ | 2 | 80 | $3^3$ | line_6 |
| 1458 | $3^6$ | 2 | 26 | $3^2$ | line_7 |
| 4374 | $3^7$ | 2 | 8 | $3^1$ | line_8 |
| | | | 2 | $3^0$ | Last_line |
| 6373 | | | | | |

Therefore: $19495/(3^8) = 2$ and rest 6373.

Remark that the division of an integer number by any power of 3 can be done only by divisions to 3, calculations of powers of 3, multiplications and additions.

129) Let M be a number in a base b. All distinct digits of M are named generalized period of M.
(For example, if M = 104001144, its generalized period is g(M) = {0, 1, 4}.)

Of course, g(M) is included in {0, 1, 2, ..., b-1}.

130) The number of generalized periods of M is equal to the number of the groups of M such that each group contains all distinct digits of M. (For example, $n_g(M) = 2$ because $M = \underbrace{104}_{1}\ \underbrace{001144}_{2}$.)

131) Length of generalized period is equal to the number of its distinct digits. (For example, $l_g(M) = 3$.)

Questions:

a) Find $n_g$, $l_g$ for $p_n$, $n!$, $n^n$, $\sqrt[n]{n}$.

b) For a given $k \geq 1$, is there an infinite number of primes $p_n$, or $n!$, or $n^n$, or $\sqrt[n]{n}$ which have a generalized period of length k ?

Same question such that the number of generalized periods be equal to k.

c) Let $a_1, a_2, \ldots, a_h$ be distinct digits. Is there an infinite number of primes $p_n$, or $n!$, or $n^n$, or $\sqrt[n]{n}$ which have as a generalized period the set $\{ a_1, a_2, \ldots, a_h \}$ ?

132) Let $\{ x_n \}_{n \geq 1}$ be a sequence of integers, and $0 \leq k \leq 9$ a digit. The sequence of position is defined as follows:

$$U_n^{(k)} = U^{(k)}(x_n) = \begin{cases} \max \{i\}, \text{ if k is the } 10^i\text{-th digit of } x_n; \\ -1, \text{ otherwise.} \end{cases}$$

(For example: if $x_1 = 5$, $x_2 = 17$, $x_3 = 775$, and $k = 7$, then $U_1^{(7)} = U^{(7)}(x_1) = -1$, $U_2^{(7)} = 0$, $U_3^{(7)} = \max \{1, 2\} = 2$.)

a) Study $\{U^{(k)}(p_n)\}_n$, where $\{p_n\}_n$ is the sequence of primes. Convergence, monotony.

b) Same question for the sequences: $x_n = n!$ and $x_n = n^n$.

More generally: when $\{x_n\}_n$ is a sequence of rational numbers, and k belongs to N.

133) Criterion for coprimes:

If a, b are strictly positive coprime integers, then:

$$a^{F(b)+1} + b^{F(a)+1} \cong a + b \ (\bmod a \cdot b),$$

where F is Euler's totient.

134) Congruence function:

$L : Z^2 \longrightarrow Z$, $L(x, m) = (x + c_1) \ldots (x + c_{F(m)})$,

where F is Euler's totient, and all $c_i$, $1 \le i \le F(m)$, are modulo m primitive rest classes.

Reference:

Florentin Smarandache, "A numerical function in the congruence theory", in <Libertas Mathematica>, Texas State University, Arlington, Vol. XII, 1992, pp. 181-5.

135) Generalization of Euler's Theorem:

If a, m are integers, m = 0, then

$$a^{(F(m_s)+s)} \equiv a^s \ (\bmod m),$$

where F is Euler's totient, and $m_s$ and s are obtained by the following algorithm:

$$(0)\begin{cases} a = a_0 d_0 ; \quad (a_0 , m_0) = 1 \\ m = m_0 d_0 ; \quad d = 1 \end{cases}$$

$$(1)\begin{cases} d_0 = d^1_0 \, d_1 ; \quad (d^1_0 , m_1) = 1 \\ m_0 = m_1 d_1 ; \quad d_1 = 1 \end{cases}$$

.................................................

$$\begin{cases} d_{s-2} = d^1_{s-2} \, d_{s-1} ; \quad (d^1_{s-2} , m_{s-1}) = 1 \end{cases}$$

(s-1)

$m_{s-2} = m_{s-1}d_{s-1}$ ;  $d_{s-1} = 1$

(s) $\begin{cases} d_{s-1} = d^{1}_{s-1} d_s ; & (d^{1}_{s-1}, m_s) = 1 \\ m_{s-1} = m_s d_s ; d_s = 1. \end{cases}$

[This is a generalization of Euler's theorem on congruences.]

136) Smarandache simple functions:

For any positive prime number p one defines

$S_p(k)$ is the smallest integer such that $S_p(k)!$ is divisible by $p^k$ .

137) Smarandache function:

S(n) is defined as the smallest integer such that S(n)! is divisible by n (for n = 0).

If the canonical factorization of n is

$$n = p_1^{k_1} \ldots p_s^{k_s}$$

then S(n) = max { $S_{p_i}(k_i)$ }, where $S_{p_i}$ are Smarandache simple functions.

138) Prime Equation Conjecture:

Let $k > 0$ be an integer. There is only a finite number of solutions in integers p, q, x, y, each greater than 1, to the equation

$$x^p - y^q = k.$$

For $k = 1$ this was conjectured by Cassels (1953) and proved by Tijdeman (1976).

139) Generalized Prime Equation Conjecture:

Let $k >= 2$ be a positive integer. The Diophantine equation

$$y = 2x_1 x_2 \dots x_k + 1$$

has infinitely many solutions in distinct primes $y, x_1, x_2, \dots, x_k$.

(For example: when $k = 3$, $647 = 2x17x19+1$;
when $k = 4$, $571 = 2x3x5x19+1$, etc.)

140) Progressions:

How many primes do the following progressions contain:

a) $\{ a \cdot p_n + b \}$, $n = 1, 2, 3, \dots$, where $(a, b) = 1$ and $p_n$ is the n-th prime?

b) $\{ a^n + b \}$, n = 1, 2, 3, ..., where (a, b) = 1, and a is different from -1, 0, +1 ?

c) $\{ n^n + 1 \}$ and $\{ n^n - 1 \}$, n = 1, 2, 3, ... ?

Reference:

Florentin Smarandache, "Only problems, not solutions!", Xiquan Publishing House, Phoenix, Chicago, 1990, Problem 17.

141) Inequality:

$$n! > k^{n-k+1} \prod_{i=0}^{k-1} \left( \frac{n-i}{k} \right)!$$

for any non-null positive integers n and k.

If k =2 (for example), one obtains:

$$n! > 2^{n-1} \left( \frac{n-i}{2} \right)! \left( \frac{n}{2} \right)!$$

and if k = 3 (another example), one obtains:

$$n! > 3^{n-2} \left( \frac{n-2}{3} \right)! \left( \frac{n-1}{3} \right)! \left( \frac{n}{3} \right)!$$

Reference:

Florentin Smarandache, "Problèmes avec et sans ... problèmes!", Somipress, Fès, Morocco, 1983, Problèmes 7.88 & 7.89, pp. 110-1.

142) Divisibility Theorem:

If a and m are integers, and m > 0, then

$$( a^m - a ) ( m - 1 )!$$

is divisible by m.

Reference:

Florentin Smarandache, "Problèmes avec et sans ... problèmes!", Somipress, Fes, Morocco, 1983, Problème 7.140, pp. 173-4.

143) Dilemmas:

Is it true that for any question there exists at least an answer?

Reciprocally:

Is any assertion the result of at least a question?

Reference:

Florentin Smarandache, "Only problems, not solutions!", Xiquan Publishing House, Phoenix, Chicago, 1990, Problem 5, p. 8.

144) Surface points:

a) Let n be an integer $\geq 5$. Find a minimum number M(n) such that anyhow are chosen M(n) points in space, four by four non-coplanar, there exist n points among these which belong to the surface of a sphere.

b) Same question for an arbitrary space body (for example: cone, cube, etc.).

c) Same question in plane (for $n \geq 4$, and the points are chosen three by three non-colinear).

Reference:

Florentin Smarandache, "Only problems, not solutions!", Xiquan Publishing House, Phoenix, Chicago, 1990, Problem 9, p. 10.

145) Inclusion Problems:

a) Find a method to get the maximum number of circles of radius 1 included in a given planar figure, at most tangential two by two or tangential to the border of the planar figure.
Study the general problem when "circle" are replaced by an arbitrary planar figure.

b) Same question for spheres of radius 1 included in a given space body.
Study the general problem when "sphere" are replaced by an arbitrary space body.

Reference:

Florentin Smarandache, "Only problems, not solutions!", Xiquan Publishing House, Phoenix, Chicago, 1990, Problem 8, p. 10.

146) Convex Polyhedrons:

a) Given n points in space, four by four non-coplanar, find the maximum number M(n) of points which constitute the vertexes of a convex polyhedron.
Of course, $M(n) \geq 4$.

b) Given n points in space, four by four non-coplanar, find the minimum number $N(n) \geq 5$ such that: any N(n) points among these do not constitute the vertexes of a convex polyhedron.

Of course, N(n) may not exist.

147) Integral Points:

How many non-coplanar points in space can be drawn at integral distances each from other?

Is it possible to find an infinite number of such points?

148) Counter:

C(a, b) = how many digits of "a" the number b contains.

Study, for example, $C(1, p_n)$, where $p_n$ is the n-th prime.

Same for: C(1, n!), $C(1, n^n)$.

149) Maximum points:

Let d > 0. Question:

What is the maximum number of points included in a given planar figure (generally: in a space body) such that the distance between any two points is greater or equal than d ?

150) Minimum Points:

Let d > 0. Question:

What is the minimum number of points $\{A_1, A_2, ... \}$ included in a given planar figure (generally: in a space body) such that if another point A is included in that figure then there exists at least an $A_i$ with the distance $|AA_i| < d$ ?

House, Phoenix, Chicago, 1990, Problem 2, p. 6.

151 Increasing Repeated Compositions:

Let g be a function, g : N ---> N, such that g(n) > n for all natural n.
An increasing repeated composition related to g and a given positive number m is defined as below:

$F_g : N \longrightarrow N$, $F_g(n) = k$, where k is the smallest integer such that $g(...g(n)...) \geq m$ (g is composed k times).

Study, for example, $F_s$, where s is the function that associates to each non-null positive integer n the sum of its positive divisors.

152) Decreasing Repeated Compositions:

Let g be a non-constant function, g : N ---> N, such that g(n) <= n for all natural n.
An decreasing repeated composition related to g is defined as below:

$f_s : N \longrightarrow N$, $f_s(n) = k$, where k is the smallest integer such that g(...g(n)...) = constant (g is composed k times).

Study, for example, $f_d$, where d is the function that associates to each non-null positive integer n the number of its positive divisors.
In this particular case, the constant is 2.
Same for $\prod(n)$= the number of primes not exceeding n,

and for p(n) = the largest prime factor of n,
and for o(n) = the number of distinct prime factors of n.

153) Coloration Conjecture:

Anyhow all points of an m-dimensional Euclidian space are colored with a finite number of colors, there exists a color which fulfills all distances.

154) Primes:

Let $a_1, ..., a_n$, be distinct digits, $1 \leq n \leq 9$.
How many primes can we construct from all these digits only (eventually repeated) ?

(More generally: when $a_1, ..., a_n$, and n are positive integers.)
Conjecture: Infinitely many!

Reference:
Florentin Smarandache, "Only problems, not solutions!", Xiquan Publishing House, Phoenix, Chicago, 1990, Problem 3, p. 7.

155) Prime Number Theorem:

There exist an infinite number of primes which contain given digits, $a_1, a_2, ..., a_m$, in the positions $i_1, i_2, ..., i_m$,
with $i_1, i_2, ..., i_m \geq 0$ (the "i-th position" is the $10^i$-th digit).
(Of course, if $i_m = 0$, then $a_m$ must be odd and different from 5.)

References:
Florentin Smarandache, Query 762 (and Answer 762), in <Mathematics Magazine>, April, 1990.
Florentin Smarandache, "Only problems, not solutions!", Xiquan Publishing House, Phoenix, Chicago, 1990, Problem 10, p. 11.

156) Cardinality Theorem:

For any positive integers $n \geq 1$ and $m \geq 3$, find the maximum number S(n ,m) such that:
the set {1, 2, 3, ..., n} has a subset A of cardinality S(n ,m) with the property that A contains no m-term arithmetic progression.
S(n ,m) is called the cardinality number.

Study it.

References:
Florentin Smarandache, "Arithmetic progressions: Problem 88-5", in <The Mathematical Intelligencer>, Vol. 11, No. 1, 1989;
E. Kurt Tekolste (Wayne, PA), "Solution to Problem 88-5", in <The Mathematical Intelligencer>, Vol. 11, No. 1, 1989.

157) Concatenated Natural Sequence:

1,22,333,4444,55555,666666,7777777,88888888,999999999,
10101010101010101010,1111111111111111111111,
121212121212121212121212,13131313131313131313131313,
1414141414141414141414141414,151515151515151515151515151515,...

158) Concatenated Prime Sequence (called Smarandache-Wellin numbers):
2, 23, 235, 2357, 235711, 23571113, 2357111317, 235711131719, 23571113171923, ...

159) Back Concatenated Prime Sequence:
2, 32, 532, 7532, 117532, 13117532, 1713117532, 191713117532, 23191713117532, ...

Conjecture: There are infinitely many primes among the first sequence numbers!

160) Concatenated Odd Sequence:
1, 13, 135, 1357, 13579, 1357911, 135791113, 13579111315, 1357911131517, ...

161) Back Concatenated Odd Sequence:
1, 31, 531, 7531, 97531, 1197531, 131197531, 15131197531, 1715131197531, ...

Conjecture: There are infinitely many primes among these numbers!

162) Concatenated Even Sequence:
2, 24, 246, 2468, 246810, 24681012, 2468101214, 246810121416, ...

163) Back Concatenated Even sequence:
2, 42, 642, 8642, 108642, 12108642, 1412108642, 161412108642, ...

Conjecture: None of them is a perfect power!

164) Concatenated S-Sequence {generalization}:
Let s1, s2, s3, s4, ..., sn, ... be an infinite integer sequence (noted by S).

Then:

$$s1, \overline{s1s2}, \overline{s1s2s3}, \overline{s1s2s3s4}, \overline{s1s2s3s4...sn}, ...$$

is called the Concatenated S-sequence,

$$s1, \overline{s2s1}, \overline{s3s2s1}, \overline{s4s3s2s1}, \overline{sn...s4s3s2s1}, ...$$

is called the Back Concatenated S-sequence.

Questions: a) How many terms of the Concatenated S-sequence belong to the initial S-sequence?

b) Or, how many terms of the Concatenated S-sequence verify the relation of other given sequences?

The first three cases are particular.

Look now at some other examples, when S is the sequence of squares, cubes, Fibonacci respectively (and one can go so on):

165) Concatenated Square Sequence:
1, 14, 149, 14916, 1491625, 149162536, 14916253649, 1491625364964, ...

166) Back Concatenated Square Sequence:
1, 41, 941, 16941, 2516941, 362516941, 49362516941, 6449362516941, ...

How many of them are perfect squares?

167) Concatenated Cubic Sequence:
1, 18, 1827, 182764, 182764125, 182764125216, 182764125216343, ...

168) Back Concatenated Cubic Sequence:
1, 81, 2781, 642781, 125642781, 216125642781, 343216125642781, ...

How many of them are perfect cubes?

169) Concatenated Fibonacci Sequence:
1, 11, 112, 1123, 11235, 112358, 11235813, 1123581321, 112358132134, ...

170) Back Concatenated Fibonacci Sequence:
1, 11, 211, 3211, 53211, 853211, 13853211, 2113853211, 342113853211, ...

Does any of these numbers is a Fibonacci number?

171) Power Function:

SP(n) is the smallest number m such that $m^m$ is divisible by n.

The following sequence SP(n) is generated:
1, 2, 3, 2, 5, 6, 7, 4, 3, 10, 11, 6, 13, 14, 15, 4, 17, 6, 19, 10, 21, 22, 23, 6, 5, 26, 3, 14, 29, 30, 31, 4, 33, 34, 35, 6, 37,38, 39, 20, 41, 42, ...

Remarks:

If p is prime, then SP(p) = p.
If r is square free, then SP(r) = r.

If n = $(p_1^{s_1}) \times \ldots \times (p_k^{s_k})$ and all $s_i \le p_i$ , then SP(n) =n.

If n = $p^s$, where p is prime, then:

$$SP(n) = \begin{cases} p, & \text{if } 1 \le s \le p; \\ p^2, & \text{if } p+1 \le s \le 2p^2; \\ p^3, & \text{if } 2p^{2+1} \le s \le 3p^3; \\ \ldots\ldots & \\ p^t, & \text{if } (t-1)p^{t-1} \le s \le tp^t . \end{cases}$$

172) Reverse Sequence:

1,21,321,4321,54321,654321,7654321,87654321,987654321,10987654321, 1110987654321,121110987654321,...

173) Multiplicative Sequence:

2,3,6,12,18,24,36,48,54,...

General definition: if $m_1$ , $m_2$ , are the first two terms of the sequence, then $m_k$ , for $k \ge 3$, is the smallest number equal to the product of two previous distinct terms.

All terms of rank ≥ 3 are divisible by $m_1$, and $m_2$ .

In our case the first two terms are 2, respectively 3.

174) Wrong Numbers:

(A number n = $\overline{a_1 a_2 \cdots a_k}$ , of at least two digits, with the property:
the sequence $a_1$, $a_2$, ..., $a_k$ , $b_{k+1}$, $b_{k+2}$, ... (where $b_{k+i}$ is the product of the previous k terms, for any i ≥1) contains n as its term.)
The author conjectured that there is no wrong number (!)
Therefore, this sequence is empty.

175) Impotent Numbers:
2,3,4,5,7,9,11,13,17,19,23,25,29,31,37,41,43,47,49,53,59,61,...

(A number n those proper divisors product is less than n.)
Remark: this sequence is { p, $p^2$ ;  where p is a positive prime }.

176) Random Sieve:
1,5,6,7,11,13,17,19,23,25,29,31,35,37,41,43,47,53,59,...

General definition:
- choose a positive number $u_1$ at random;
- delete all multiples of all its divisors, except this number;
- choose another number $u_2$ greater than $u_1$ among those remaining;
- delete all multiples of all its divisors, except this second number;
... and so on.

The remaining numbers are all coprime two by two.

The sequence obtained $u_k$ , k ≥1,  is less dense than the prime number sequence, but it tends to the prime number sequence as k tends to infinite. That's why this sequence may be important.

In our case, $u_1$ = 6, $u_2$ = 19, $u_3$ = 35, ... .

177) Non-Multiplicative Sequence:

General definition: let $m_1$ , $m_2$ , ..., $m_k$ be the first k given terms of the sequence, where $k \geq 2$;
then $m_i$ , for $i \geq k+1$, is the smallest number not equal to the product of k previous distinct terms.

178) Non-Arithmetic Progression:
1,2,4,5,10,11,13,14,28,29,31,32,37,38,40,41,64,...

General definition: if $m_1$ , $m_2$, are the first two terms of the sequence, then $m_k$ , for $k \geq 3$, is the smallest number such that no3-term arithmetic progression is in the sequence.
In our case the first two terms are 1, respectively 2.

Generalization: same initial conditions, but no i-term arithmetic progression in the sequence ( for a given $i \geq 3$ ).

179) Prime Product Sequence:

2,7,31,211,2311,30031,510511,9699691,223092871,6469693231,200560490131, 7420738134811,304250263527211,...

$P_n = 1 + p_1p_2 \ldots p_n$ , where $p_k$ is the k-th prime.

Question: How many of them are prime?

180) Square Product Sequence:

2,5,37,577,14401,518401,25401601,1625702401,131681894401,1316818940001, 1593350922240001,...

$S_n = 1 + s_1s_2 \ldots s_n$ , where $s_k$ is the k-th square number.

Question: How many of them are prime?

181) Cubic Product Sequence:
2,9,217,13825,1728001,373248001,128024064001,65548320768001,...

$C_n = 1 + c_1c_2 \ldots c_n$ , where $c_k$ is the k-th cubic number.

Question: How many of them are prime?

182) Factorial Product Sequence:
2,3,13,289,34561,24883201,125411328001,5056584744960001,...

$F_n = 1 + f_1 f_2 \ldots f_n$ , where $f_k$ is the k-th factorial number.

Question: How many of them are prime?

183) U-Product Sequence {generalization}:
Let $u_n$ , $n \geq 1$, be a positive integer sequence. Then we define a sequence as follows:

$U_n = 1 + u_1 u_2 \ldots u_n$ .

Reference:
F. Smarandache, "Properties of the Numbers", University of Craiova Archives, 1975;
[see also Arizona State University Special Collection, Tempe, Arizona, USA].

184-190) Sequences of Sub-Sequences:

For each of the following 7 sequences:

a) Crescendo Sub-Sequences:

1, 1, 2, 1, 2, 3, 1, 2, 3, 4, 1, 2, 3, 4, 5, 1, 2, 3,4, 5, 6,
1, 2, 3, 4, 5, 6, 7, 1, 2, 3, 4, 5, 6, 7, 8, . . .

b) Decrescendo Sub-sequences:

1, 2, 1, 3, 2, 1, 4, 3, 2, 1, 5, 4, 3, 2, 1, 6, 5, 4,3, 2, 1,
7, 6, 5, 4, 3, 2, 1, 8, 7, 6, 5, 4, 3, 2, 1, . . .

c) Crescendo Pyramidal Sub-sequences:

1, 1, 2, 1, 1, 2, 3, 2, 1, 1, 2, 3, 4, 3, 2, 1,
1, 2, 3, 4, 5, 4, 3, 2, 1, 1, 2, 3, 4, 5, 6, 5, 4, 3, 2, 1, . . .

d) Decrescendo Pyramidal Sub-sequences:

5, 4, 3, 2, 1, 2, 3, 4, 5,     6, 5, 4, 3, 2, 1, 2, 3, 4, 5, 6,   . . .

e) Crescendo Symmetric Sub-sequences:

1, 1,   1, 2, 2, 1,   1, 2, 3, 3, 2, 1,   1, 2, 3, 4, 4, 3, 2,1,
1, 2, 3, 4, 5, 5, 4, 3, 2, 1,   1, 2, 3, 4, 5, 6, 6, 5, 4, 3, 2,1, . . .

f) Decrescendo Symmetric Sub-sequences:

1, 1,   2, 1, 1, 2,   3, 2, 1, 1, 2, 3,   4, 3, 2, 1, 1, 2, 3,4,
5, 4, 3, 2, 1, 1, 2, 3, 4, 5,   6, 5, 4, 3, 2, 1, 1, 2, 3, 4, 5,6,   . . .

g) Permutation Sub-sequences:

1, 2,   1, 3, 4, 2,   1, 3, 5, 6, 4, 2,   1, 3, 5, 7, 8, 6, 4,2,
1, 3, 5, 7, 9, 10, 8, 6, 4, 2,   1, 3, 5, 7, 9, 10, 8, 6, 4, 2,. . .

Find a formula for the general term of each sequence.

Solutions:

For purposes of notation in all problems, let

a(n)

denote the n-th term in the complete sequence and

b(n)

the n-th subsequence. Therefore, a(n) will be a number and b(n) a sub-sequence.

a) Clearly, b(n) contains n terms. Using a well-known summation formula, at the end of b(n) there would be a total of

$\frac{n(n+1)}{2}$ terms.

Therefore, since the last number of b(n) is n,

$$a\left(\frac{n(n+1)}{2}\right)=n.$$

Finally, since this would be the terminal number in the sub-sequence

b(n) = 1, 2, 3, . . . , n

the general formula is

$$a\left(\left(\frac{n(n+1)}{2}\right)-i\right)=n-i$$

for $n \geq 1$ and $0 \leq i \leq n - i$.

b) With modifications for decreasing rather than increasing, the proof is essentially the same. The final formula is

$$a\left(\left(\frac{n(n+1)}{2}\right)-i\right)=1+i$$

for $n \geq 1$ and $0 \leq i_7 \leq n - 1$.

c) Clearly, b(n) has 2n - 1 terms. Using the well-known formula of summation

$1 + 3 + 5 + . . . + (2n - 1) = n^2$.

the last term of b(n) is in position $n^2$ and $a(n^2) = 1$. The largest number in b(n) is n, so counting back n - 1 positions, they increase in value by one each step until n is reached.

$a(n^2 - i) = 1 + i$, for $0 \leq i \leq n-1$.

After the maximum value at n-1 positions back from $n^2$, the values decrease by one. So at the nth position back, the value is n-1, at the (n-1)-th position back the value is n-2 and so forth.

$a(n^2 - n - i) = n - i - 1$

for $0 \leq i \leq n - 2$.

d) Using similar reasoning

$a(n^2) = n$ for $n \geq 1$

and

$a(n^2 - i) = n - i, \quad \text{for } 0 \leq i \leq n-1$

$a(n^2 - n - i) = 2 + i, \text{ for } 0 \leq i \leq n-2.$

e) Clearly, b(n) contains 2n terms. Applying another well-known summation formula

$2 + 4 + 6 + \ldots + 2n = n(n+1), \text{ for } n \geq 1.$

Therefore, a(n(n+1)) = 1. Counting backwards n-1 positions, each term decreases by 1 up to a maximum of n.

$a((n(n+1))-i) = 1 + i, \text{ for } 0 \leq i \leq n-1.$

The value n positions down is also n and then the terms decrease by one back down to one.

$a((n(n+1))-n-i) = n - i, \text{ for } 0 \leq i \leq n-1.$

f) The number of terms in b(n) is the same as that for (e). The only difference is that now the direction of increase/decrease is reversed.

$a((n(n+1))-i) = n - i, \text{ for } 0 \leq i \leq n-1.$

$a((n(n+1))-n-i) = 1 + i, \text{ for } 0 \leq i \leq n-1.$

g) Given the following circular permutation on the first
n integers.

$$\Psi_n = \begin{vmatrix} 1 & 2 & 3 & 4 & \ldots & n-2 & n-1 & n \\ 1 & 3 & 5 & 7 & \ldots & 6 & 4 & 2 \end{vmatrix}$$

Once again, b(n) has 2n terms. Therefore,

$a(n(n+1)) = 2.$

Counting backwards n-1 positions, each term is two larger than the successor

$a((n(n+1))-i) = 2 + 2i, \quad \text{for } 0 \leq i \leq n-1.$

The next position down is one less than the previous and after that, each term is again two less the successor.

$a((n(n+1))-n-i) = 2n - 1 - 2i$, for $0 \le i \le n-1$.

As a single formula using the permutation

$a((n(n+1)-i) = \Psi_n(2n-i)$, for $0 \le i \le 2n-1$.

191) $S^2$ function (numbers):
1, 2, 3, 2, 5, 6, 7, 4. 3, 10, 11, 6, 13, 14, 15, 4, 17, 6,19, 10,
21, 22, 23, 12, 5, 26, 9, 14, 29, 30, 31, 8, 33, ...

( S2(n) is the smallest integer m such that $m^2$ is divisible by n)

192) $S^3$ function (numbers):
1, 2, 3, 2, 5, 6, 7, 8, 3, 10, 11, 6, 13, 14, 15, 4, 17, 6,19, 10,
21, 22, 23, 6, 5, 26, 3, 14, 29, 30, 31, 4, 33, ...

( S3(n) is the smallest integer m such that $m^3$ is divisible by n)

193) Anti-Symmetric Sequence:
11,1212,123123,12341234,1234512345,123456123456,12345671234567,
1234567812345678,123456789123456789,12345678910123456789 10,
1234567891011123456789 1011,123456789101112123456789101112,...

194-202) Recurrence Type Sequences:

A. 1,2,5,26,29,677,680,701,842,845,866,1517,458330,458333,
458354,...
( ss2(n) is the smallest number, strictly greater than the previous one, which is the squares sum of two previous distinct terms of the sequence;
in our particular case the first two terms are 1 and 2. )

Recurrence definition:

1) The number a belongs to SS2;

2) If b, c belong to SS2, then $b^2 + c^2$ belongs to SS2too;

3) Only numbers, obtained by rules 1) and/or 2) applied a finite number of times, belong to SS2.

The sequence (set) SS2 is increasingly ordered.

[ Rule 1) may be changed by: the given numbers a1, a2,..., ak , where $k \geq 2$, belong to SS2. ]

B. 1,1,2,4,5,6,16,17,18,20,21,22,25,26,27,29,30,31,36,37,38,40,41, 42,43,45,46,...

( ss1(n) is the smallest number, strictly greater than the previous one (for $n \geq 3$), which is the squares sum of one or more previous distinct terms of the sequence; in our particular case the first term is 1. )

Recurrence definition:

1) The number a belongs to SS1;

2) If $b_1, b_2, ..., b_k$ belong to SS1, where $k \geq 1$, then $b_1^2 + b_2^2 + ... + b_k^2$ belongs to SS1 too;

3) Only numbers, obtained by rules 1) and/or 2) applied a finite number of times, belong to SS1.

The sequence (set) SS1 is increasingly ordered.

[ Rule 1) may be changed by: the given numbers $a_1, a_2,...,$ $a_k$ , where $k \geq 1$, belong to SS1. ]

C. 1,2,3,4,6,7,8,9,11,12,14,15,16,18,19,21,...

( nss2(n) is the smallest number, strictly greater than the previous one, which is NOT the squares sum of two previous distinct terms of the sequence; in our particular case the first two terms are 1 and 2. )

Recurrence definition:

1) The numbers $a \leq b$ belong to NSS2;

2) If b, c belong to NSS2, then $b^2 + c^2$ DOES NOT belong to NSS2; any other numbers belong to NSS2;

3) Only numbers, obtained by rules 1) and/or 2) applied a finite number of times, belong to NSS2.

The sequence (set) NSS2 is increasingly ordered.

[ Rule 1) may be changed by: the given numbers $a_1, a_2,...,$ $a_k$ , where $k \geq 2$, belong to NSS2. ]

D. 1,2,3,6,7,8,11,12,15,16,17,18,19,20,21,22,23,24,25,26,27,28, 29,30,31,32,33,34,35,38,39,42,43,44,47,...

( nss1(n) is the smallest number, strictly greater than the previous distinct terms of the sequence; in our particular case the first term is 1. )

Recurrence definition:

1) The number a belongs to NSS1;

2) If $b_1$, $b_2$, ..., $b_k$ belong to NSS1, where $k \geq 1$, then $b_1^2 + b_2^2 + ... + b_k^2$ DO NOT belong to NSS1; any other numbers belong to NSS1;

3) Only numbers, obtained by rules 1) and/or 2) applied a finite number of times, belong to NSS1.

The sequence (set) NSS1 is increasingly ordered.

[ Rule 1) may be changed by: the given numbers $a_1$, $a_2$,..., $a_k$ , where $k \geq 1$, belong to NSS1. ]

E. 1,2, 9, 730,737, 389017001,389017008,389017729,...

( cs2(n) is the smallest number, strictly greater than the previous one, which is the cubes sum of two previous distinct terms of the sequence; in our particular case the first two terms are 1 and 2. )

Recurrence definition:

1) The numbers $a \leq b$ belong to CS2;

2) If c, d belong to CS2, then $c^3 + d^3$ belongs to CS2too;

3) Only numbers, obtained by rules 1) and/or 2) applied a finite number of times, belong to CS2.

The sequence (set) CS2 is increasingly ordered.

[ Rule 1) may be changed by: the given numbers $a_1$, $a_2$,..., $a_k$ , where $k \geq 2$, belong to CS2. ]

F. 1,1, 2, 8,9,10, 512,513,514,520,521,522,729,730,731,737,738, 739,1241,...

( cs1(n) is the smallest number, strictly greater than the previous one (for $n \geq 3$), which is the cubes sum of one or more previous distinct terms of the sequence; in our particular case the first term is 1. )

Recurrence definition:

1) The number a belongs to CS1.
2) If $b_1, b_2, ..., b_k$ belong to CS1, where $k \geq 1$, then $b_1^3 + b_2^3 + ... + b_k^3$ belongs to CS2 too.
3) Only numbers, obtained by rules 1) and/or 2) applied a finite number of times, belong to CS1.
The sequence (set) CS1 is increasingly ordered.
[ Rule 1) may be changed by: the given numbers $a_1, a_2,...,$ .

G. 1,2,3,4,5,6,7,8,10,11,12,13,14,15,16,17,18,19,20,21,22,23, 24,25,26,27,29,30,31,32,33,34,36,37,38,...

( ncs2(n) is the smallest number, strictly greater than the previous one, which is NOT the cubes sum of two previous distinct terms of the sequence;
in our particular case the first two terms are 1 and 2. )

Recurrence definition:

1) The numbers $a \leq b$ belong to NCS2;
2) If c, d belong to NCS2, then $c^3 + d^3$ DOES NOT belong to NCS2; any other numbers do belong to NCS2;
3) Only numbers, obtained by rules 1) and/or 2) applied a finite number of times, belong to NCS2.
The sequence (set) NCS2 is increasingly ordered.
[ Rule 1) may be changed by: the given numbers $a_1, a_2,...,$ $a_k$ , where $k \geq 2$, belong to NCS2. ]

H. 1,2,3,4,5,6,7,10,11,12,13,14,15,16,17,18,19,20,21,22,23,24, 25,26,29,30,31,32,33,34,37,38,39,...

( ncs1(n) is the smallest number, strictly greater than the previous one, which is NOT the cubes sum of one or more previous distinct terms of the sequence;
in our particular case the first term is 1. )

Recurrence definition:

1) The number a belongs to NCS1.
2) If $b_1, b_2, ..., b_k$ belong to NCS1, where $k \geq 1$, then $b_1^2 + b_2^2 + ... + b_k^2$ DO NOT belong to NCS1.
3) Only numbers, obtained by rules 1) and/or 2) applied a finite number of times, belong to NCS1.
The sequence (set) NCS1 is increasingly ordered.
[ Rule 1) may be changed by: the given numbers $a_1, a_2,...,$ $a_k$ , where $k \geq 1$, belong to NCS1. ]

I. General-Recurrence (Positive) Type Sequence:

General (positive) recurrence definition:

Let $k \geq j$ be natural numbers, and $a_1, a_2, ..., a_k$ given elements, and R a j-relationship (relation among j elements).

Then:

1) The elements $a_1, a_2, ..., a_k$ belong to SGPR.

2) If $m_1, m_2, ..., m_j$ belong to SGPR, then $R(m_1, m_2, ..., m_j)$ belongs to SGPR too.

3) Only elements, obtained by rules 1) and/or 2) applied a finite number of times, belong to SGPR.

The sequence (set) SGPR is increasingly ordered.

Method of construction of the (positive) general recurrence sequence:

- level 1: the given elements $a_1, a_2, ..., a_k$ belong to SGPR;
- level 2: apply the relationship R for all combinations of j elements among $a_1, a_2, ..., a_k$;
  the results belong to SGPR too;
  order all elements of levels 1 and 2 together;
  ............................................................
- level i+1:
  if $b_1, b_2, ..., b_m$ are all elements of levels 1, 2,..., i-1,
  and $c_1, c_2, ..., c_n$ are all elements of level i,
  then apply the relationship R for all combinations of j elements among $b_1, b_2, ..., b_m, c_1, c_2, ..., c_n$ such that at least a element is from the level i;
  the results belong to SRPG too;
  order all elements of levels i and i+1 together;

and so on ...

J. General-Recurrence (Negative) Type Sequence:

General (negative) recurrence definition:

Let $k \geq j$ be natural numbers, and $a_1, a_2, ..., a_k$ given elements, and R a j-relationship (relation among j elements).

Then:

1) The elements $a_1, a_2, ..., a_k$ belong to SGNR.

2) If $m_1, m_2, ..., m_j$ belong to SGNR, then $R(m_1, m_2, ..., m_j)$ does NOT belong to SGNR; any other elements do belong to SGNR.

3) Only elements, obtained by rules 1) and/or 2) applied a finite number of times, belong to SGNR.

The sequence (set) SGNR is increasingly ordered.

Method of construction of the (negative) general recurrence

sequence:

- level 1: the given elements $a_1, a_2, ..., a_k$ belong to SGNR;
- level 2: apply the relationship R for all combinations of j elements among $a_1, a_2, ..., a_k$;
  none of the results belong to SGNR;
  the smallest element, strictly greater than all $a_1, a_2, ..., a_k$, and different from the previous results, belongs to SGNR;
  order all elements of levels 1 and 2 together;
  ...........................................................
- level i+1:
  if $b_1, b_2, ..., b_m$ are all elements of levels 1, 2,..., i-1,
  and $c_1, c_2, ..., c_n$ are all elements of level i,
  then apply the relationship R for all combinations of j elements among $b_1, b_2, ..., b_m, c_1, c_2, ..., c_n$ such that at least a element is from the level i;
  none of the results belong to SRNG;
  the smallest element, strictly greater than all previous elements, and different from the previous results, belongs to SGNR;
  order all elements of levels i and i+1 together;

and so on ...

203- 215) Partition Type Sequences:

A. 1,1,1,2,2,2,2,3,4,4,...

( How many times is n written as a sum of non-null squares, disregarding the terms order;
for example:

$9 = 1^2 + 1^2 + 1^2 + 1^2 + 1^2 + 1^2 + 1^2 + 1^2 + 1^2$
$= 1^2 + 1^2 + 1^2 + 1^2 + 1^2 + 2^2$
$= 1^2 + 2^2 + 2^2$
$= 3^2$,

therefore ns(9) = 4. )

B. 1,1,1,1,1,1,1,2,2,2,2,2,2,2,2,3,3,3,3,3,3,3,3,4,4,4,5,5,
5,5,5,6,6,...

( How many times is n written as a sum of non-null cubes, disregarding the terms order;
for example:

$9 = 1^3 + 1^3 + 1^3 + 1^3 + 1^3 + 1^3 + 1^3 + 1^3 + 1^3$
$= 1^3 + 2^3$,
therefore $n_c(9) = 2$. )

C. General Partition Type Sequence:
Let f be an arithmetic function, and R a k-relation among numbers.
How many times can n be expressed under the form of;
$n = R ( f_1 (n ), f_2 (n ), ..., f_k (n ))$,
for some k and $n_1 , n_2 , ..., n_k$ such that $n_1 + n_2 + ... + n_k = n$ ?

[ Particular cases when $f(x) = x^2$, or $x^3$, or x!, or $x^x$, etc.
and the relation R is the trivial addition of numbers, or
multiplication, etc. ]

216) Reverse Sequence:
1,21,321,4321,54321,654321,7654321,87654321,987654321,10987654321,
1110987654321,121110987654321,13121110987654321,1413121110987654321,
151413121110987654321,16151413121110987654321,...

217) Non-Geometric Progression:
1,2,3,5,6,7,8,10,11,13,14,15,16,17,19,21,22,23,24,26,27,29,30,31,33,34,35,
37,38,39,40,41,42,43,45,46,47,48,50,51,53,...

General definition: if $m_1$ , $m_2$ , are the first two terms of the sequence, then $m_k$ , for $k \geq 3$, is the smallest number such that no 3-term geometric progression is in the sequence.

In our case the first two terms are 1, respectively 2.

Generalization: if $m_1 , m_2, ..., m_{i-1}$ are the first i-1 terms of the sequence, $i \geq 3$, then $m_k$ , for $k \geq i$, is the smallest number such that no i-term geometric progression is in the sequence.

218) Unary Sequence:

11,111,11111,1111111,11111111111,1111111111111,11111111111111111,
1111111111111111111,11111111111111111111111,11111111111111111111111111111,
1111111111111111111111111111111,...

u(n) = $\overline{1,1,...1,}$, $p_n$ digits of "1", where $p_n$ is the n-th prime.

The old question: are there an infinite number of primes belonging to the sequence?

219) No-prime-digit Sequence:

1,4,6,8,9,10,11,1,1,14,1,16,1,18,19,0,1,4,6,8,9,0,1,4,6,8,9,40,41,42,4,44,
4,46,4,48,49,0,...

(Take out all prime digits of n.)

Is it any number that occurs infinitely many times in this sequence ?
(for example 1, or 4, or 6, or 11, etc.).

Solution by Dr. Igor Shparlinski,
It seems that: if, say n has already occurred, then for example n3, n33, n333, etc. gives infinitely many repetitions of this number.

220) No-square-digits Sequence:

2,3,5,6,7,8,2,3,5,6,7,8,2,2,22,23,2,25,26,27,28,2,3,3,32,33,3,35,36,37,38,
3,2,3,5,6,7,8,5,5,52,52,5,55,56,57,58,5,6,6,62,...

(Take out all square digits of n.)

221-222) Polygons and Polyhedrons.

A) Polygons:

Definition:

Let $V_1$, $V_2$, ..., $V_n$ be an n-sided polygon, $n \geq 3$, and S a given set of m elements, $m \geq n$. In each vertex $V_i$, $1 \leq i \leq n$, of the polygon one put at most an element of S, and on each side $s_i$, $1 \leq i \leq n$, of the polygon one put $e_i \geq 0$ elements from A. Let R be a given relationship of two or more elements on S.
All elements of S should be put on the polygon's sides and vertexes such that the relationship R on each side give the same result.

What connections should be among the number of sides of the polygon, the set S (what elements and how many), and the relation R in order for this problem to have solution(s) ?

B) Polyhedrons:

I) Polyhedrons (With Edge Points):
Similar definition generalized in a 3-dimensional space:
the points are set on the vertexes, and edges only of the polyhedron (not on inside faces).

II) Polyhedrons (With Face Points):
Similar definition generalized in a 3-dimensional space:
the points are set on the vertexes and inside faces of the polyhedron (not on the edges).

III) Polyhedrons (With Edge/Face Points):
Similar definition generalized in a 3-dimensional space:
the points are set on the vertexes, on edges, and on inside faces of the polyhedron.

IV) Polyhedrons (With Edge Points):
Similar definition generalized in a 3-dimensional space:
the points are set on the edges, and inside faces only of the polyhedron (not in vertexes).

What connections should be among the number of vertexes/faces/edges of the polyhedron, the set S (what elements and how many), and the relation R in order for this problem to have solution(s) ?
[v + f = e + 2, Euler's Result, where v = number of vertexes of the polyhedron, f = number of faces, e = number of edges.]

Examples: I)

a) For an equilateral triangle, the set N6 = {1,2,3,4,5,6}, the addition, and on each side to have three elements exactly, and in each vertex an element, one gets:

```
      1
   6     5      the sum on each of the three sides is 9;
 2    4     3
```

[the minimum elements 1,2,3 are put in the vertexes]

```
      6
   1     2      the sum on each of the three sides is 12;
 5    3     4
```

[the maximum elements 6,5,4 are put in the vertexes]

There are no other possible combinations to keep (the sum constant).

Therefore:
The Triangular Index SGI(3) = (9, 12; 2),
which means: the minimum value, the maximum value, the total number of combinations respectively.

b) For a square, the set N12 = {1,2,3,4,5,6, ..., 12}, the addition, and on each side to have four elements exactly, and in each vertex an element.

What is the Square Index SGI(4) ?

c) For a regular pentagon, the set N20 = {1,2,3,4,5,6, ..., 20}, the addition, and on each side to have five elements exactly, and in each vertex an element.

What is the Pentagonal Index SGI(5) ?

d) For a regular hexagon, the set N30 = {1,2,3,4,5,6, ..., 30}, the addition, and on each side to have six elements exactly.

What is the Hexagonal Index SGI(6) ?

e) And generally speaking:
For an n-sided regular polygon,
the set $N(n^2-n) = \{1,2,3,4,5,6, ..., n^2-n\}$,
the addition, and on each side to have n elements exactly.

What is the n-Sided Regular Polygon Index SGI(n) ?

Examples: II)

a) For a regular tetrahedron (4 vertexes, 4 triangular faces, 6 edges), the set N10 = {1,2,3,...,10}, the addition, and on each edge to have 3 elements exactly (in each vertex one element) -- to get the same sum.

b) For a regular hexahedron (8 vertexes, 6 square faces, 12 edges), the set N32 = {1,2,3,...,32}, the addition, and on each edge to have 4 elements exactly (in each vertex one element) -- to get the same sum.

c) For a regular octahedron (6 vertexes, 8 triangular faces, 12 edges),

the set N18 = {1,2,3,...,18}, the addition, and on each edge to have 3 elements exactly (in each vertex one element) -- to get the same sum.

d) For a regular dodecahedron (14 vertexes, 12 square faces, 24 edges), the set N62 = {1,2,3,...,62}, the addition, and on each edge to have 4 elements exactly (in each vertex one element) -- to get the same sum.

Compute the Polyhedron Index SHI(v,e,f) = (m, M; nc)
in each case,
where v = number of vertexes of the polyhedron, e = number of edges,
f = number of faces,
and m = minimum sum value, M = maximum sum value, nc = the total number of combinations.

e) Or other versions for the previous particular polyhedrons:
- putting elements on inside faces and vertexes only (not on edges);
- or putting elements on inside faces and edges and vertexes too;
- or putting elements on inside faces and edges only (not on vertexes).

223) Lucky Numbers:

A number is called a Lucky Number if by a wrong operation one gets to the true result!

The wrong operation should somehow be similar to a correct operation, see our students' common math mistakes:

$$\sqrt{a} - \sqrt{b} = \sqrt{a - b}$$

of course avoiding the trivial cases (if b = 0 or a = b, in the previous case, for example).

Question: Are There Finitely or Infinitely Lucky Numbers?

224) The Lucky Method/Algorithm/Operation/etc.

Is said to be any incorrect method or algorithm or operation etc. which leads to a correct result. The wrong calculation should be fun, somehow similarly to the students' common mistakes, or to produce confusions or paradoxes.

Can someone give an example of a Lucky Derivation, or Integration, or Solution to a Differential Equation?

225) The Sequence of 3's and 1's containing many primes:
31, 331, 3331, 33331, 333331, 3333331, 33333331, 333333331, ...

226) The Sequence of 3's and 1's containing many primes:
31, 31331, 313313331, 31331333133331, 3133133313333133331, 313313331333313333313333331, 31331333133331333331333333133333331, ...

227) Generalization:
Let f, g, h, ... be n integer functions. By concatenation of all of them one gets a composed sequence of the form:

$$\overline{f(1)g(1)h(1)\ldots}, \ldots \overline{f(n)g(n)h(n)\ldots},$$

where some function value(s) f(i) or g(i) or h(i), for $i \geq 1$, may be empty.

228) Recurrence of Smarandache type Functions:
A function

$$f: N \to N$$

such that if (a, b) = 1, then f(a·b) = max {f(a), f(b)}.
(Definition introduced by Sabin Tabirca, England.)

For example:
a) f(n) = 1
b) g(n) = max {p, where p is prime and p divides n},

229-231) G Add-On Sequence (I)

"Personal Computer World" Numbers Count of February 1997 presented some of the Smarandache Sequences and related open problems.

Let G = $\{g_1, g_2, ..., g_k, ...\}$ be an ordered set of positive integers with a given property G.

Then the corresponding G Add-On Sequence is defined through

$$SG = \{a_i : a_1 = g_1,\ a_k = a_{k-1} 10^{1+\log_{10}(g_k)} + g_k,\ k \geq 1\}.$$

H. Ibstedt studied some particular cases of this sequence, that he has presented at the First International Conference on Smarandache type Notions in Number Theory, University of Craiova, Romania, August 21-24, 1997.

a) Examples of G Add-On Sequences (II)

The following particular cases are studied:

a1) Odd Sequence is generated by choosing

$G = \{1, 3, 5, 7, 9, 11, ...\}$, and it is:

1, 13, 135, 1357, 13579, 1357911, 13571113, ... .

Using the elliptic curve prime factorization program we find the first five prime numbers among the first 200 terms of this sequence, i.e. the ranks 2, 15, 27, 63, 93.

But are they infinitely or finitely many?

a2) Even Sequence is generated by choosing

$G = \{2, 4, 6, 8, 10, 12, ...\}$, and it is:

2, 24, 246, 2468, 246810, 24681012, ... .

Searching the first 200 terms of the sequence we didn't find any n-th perfect power among them, no perfect square, nor even of the form 2p, where p is a prime or pseudo-prime.

Conjecture: There is no n-th perfect power term!

a3) Prime Sequence is generated by choosing

$G = \{2, 3, 5, 7, 11, 13, 17, ...\}$, and it is:

2, 23, 235, 2357, 235711, 23571113, 2357111317, ... .

Terms #2 and #4 are primes; terms #128 (of 355 digits) and #174 (of 499 digits) might be, but H. Ibstedt couldn't check -- among the first 200 terms of the sequence.

Question: Are there infinitely or finitely many such primes?

(H. Ibstedt)

232) Non-Arithmetic Progression (I)

"Personal Computer World" Numbers Count of February 1997 presented some of the Smarandache Sequences and related open problems.

One of them defines the t-Term Non-Arithmetic Progression as the set:

$\{ a_i : a_i$ is the smallest integer such that $a_i > a_{i-1}$ ,
and there are at most t-1 terms in an arithmetic progression$\}$.

A QBASIC program is designed to implement a strategy for building a such progression, and a table for the 65 first terms of the non-arithmetic progressions for t=3 to 15 is given.

(H. Ibstedt)

233) Concatenation Type Sequences

Let $s_1 , s_2 , s_3 , ..., s_n , ...$ be an infinite integer sequence (noted by S). Then the Concatenation is defined as:

$s_1$, $\overline{s_1 s_2}$, $\overline{s_1 s_2 s_3}$, … .

I search, in some particular cases, how many terms of this concatenated S-sequence belong to the initial S-sequence.
(H. Ibstedt)

234) Construction of Elements of the Square-Partial-Digital Subsequence

The Square-Partial-Digital Subsequence (SSPDS) is the sequence of square integers which admit a partition for which each segment is a square integer. An example is $506^2$ = 256036, which has partition 256/0/36. C. Ashbacher showed that SSPDS is infinite by exhibiting two infinite families of elements. We will extend his results by showing how to construct infinite families of elements of SSPDS containing desired patterns of digits.

Unsolved Question 1:
441 belongs to SSPDS, and his square $441^2$ = 194481 also belongs to SSPDS. Can an example be found of integers m, $m^2$, $m^4$ all belonging to SSPDS?

Unsolved Question 2:
It is relatively easy to find two consecutive squares in SSDPS, i.e. $12^2$ = 144 and $13^2$ = 169.
Does SSDPS also contain three or more consecutive squares?
What is the maximum length?

235) Prime-Digital Sub-Sequence

"Personal Computer World" Numbers Count of February 1997 presented some of the Smarandache Sequences and related open problems.
One of them defines the Prime-Digital Sub-Sequence as the ordered set of primes whose digits are all primes:
2, 3, 5, 7, 23, 37, 53, 73, 223, 227, 233, 257, 277, ... .
We used a computer program in Ubasic to calculate the first 100 terms of the sequence. The 100-th term is 33223.
Sylvester Smith [1] conjectured that the sequence is infinite. In this paper we will prove that this sequence is in fact infinite.
(H. Ibstedt)

India, Vol. 15E, No. 1, 1996, pp. 101-107.

236-237) Special Expressions.

a) Perfect Powers in Special Expressions (I)

How many primes are there in the Expression:

$$x^y + y^x,$$

where gcd(x, y) = 1 ? [J. Castillo & P. Castini]
K. Kashihara announced that there are only finitely many numbers of the above form which are products of factorials.
In this note we propose the following conjecture:
Let a, b, and c three integers with ab nonzero. Then the equation:
$ax^y + by^x = cz^n$, with x, y, $n \geq 2$, and gcd(x, y) = 1,
has finitely many solutions (x, y, z, n).
And we prove some particular cases of it.
(F. Luca)

b) Products of Factorials in Special Expressions (II)

J. Castillo ["Mathematical Spectrum", Vol. 29, 1997/8, 21] asked how many primes are there in the n-Expression:

$$x1^{x2} + x2^{x3} + \ldots + xn^{x1},$$

where n > 1, x1, x2, ..., xn > 1, and gcd (x1, x2, ..., xn) = 1 ?
[This is a generalization of the 2-Expression: $x^y + y^x$.]
In this note we announce a lower bound for the size of the largest prime divisor of an expression of type $ax^y + b\,y^x$, where ab is nonzero, x, $y \geq 2$, and gcd (x, y) = 1.
(F. Luca)

238) The General Periodic Sequence
Let S be a finite set, and $f : S \longrightarrow S$ be a function defined for all elements of S.
There will always be a periodic sequence whenever we repeat the composition of the function f with itself more times than card(S), accordingly to the box principle of Dirichlet.
[The invariant sequence is considered a periodic sequence whose period length has one term.]

Thus the General Periodic Sequence is defined as:

$a_1 = f(s)$, where s is an element of S;
$a_2 = f(a_1) = f(f(s))$;
$a_3 = f(a_2) = f(f(a_1)) = f(f(f(s)))$;
and so on.

We particularize S and f to study interesting cases of this type of sequences.
(M. R. Popov)

239-244) Periodic Sequences.

a) The Two-Digit Periodic Sequence (I)

Let N1 be an integer of at most two digits and let N1' be its digital reverse. One defines the absolute value N2 = abs (N1 - N1'). And so on: N3 = abs (N2 - N2'), etc. If a number N has one digit only, one considers its reverse as Nx10 (for example: 5, which is 05, reversed will be 50). This sequence is periodic.
Except the case when the two digits are equal, and the sequence becomes:

N1, 0, 0, 0, ...

the iteration always produces a loop of length 5, which starts on the second or the third term of the sequence, and the period is 9, 81, 63, 27, 45 or a cyclic permutation thereof.

(The next periodic sequences are extracted from this paper too).

b) The n-Digit Periodic Sequence (II)

Let N1 be an integer of at most n digits and let N1' be its digital reverse. One defines the absolute value N2 = abs (N1 - N1').

And so on: N3 = abs (N2 - N2'), etc. If a number N has less than n digits, one considers its reverse as N'x($10^k$), where N' is the reverse of N and k is the number of missing digits, (for example: the number 24 doesn't have five digits, but can be written as 00024, and reversed will be 42000). This sequence is periodic according to Dirichlet's box principle.

The 3-Digit Periodic Sequence (domain 100 ≤ N1≤ 999):
- there are 90 symmetric integers, 101, 111, 121, ..., for which N2 = 0;
- all other initial integers iterate into various entry points of the same periodic subsequence (or a cyclic permutation thereof) of five terms:

99, 891, 693, 297, 495.

The 4-Digit Periodic Sequence (domain 1000≤ N1≤ 9999):
- the largest number of iterations carried out in order to reach the first member of the loop is 18, and it happens for N1 = 1019;
- iterations of 8818 integers result in one of the following loops (or a cyclic permutation thereof): 2178, 6534; or 90, 810, 630, 270, 450; or 909, 8181, 6363, 2727, 4545; or 999, 8991, 6993, 2997, 4995;
- the other iterations ended up in the invariant 0.

(H. Ibstedt)

c) The 5-Digit and 6-Digit Periodic Sequences (III)

Let N1 be an integer of at most n digits and let N1' be its digital reverse. One defines the absolute value N2 = abs (N1 - N1'). And so on: N3 = abs (N2 - N2'), etc. If a number N has less than n digits, one considers its reverse as N'x(10^k), where N' is the reverse of N and k is the number of missing digits, (for example: the number 24 doesn't have five digits, but can be written as 00024, and reversed will be 42000). This sequence is periodic according to Dirichlet's box principle, leading to invariant or a loop.

The 5-Digit Periodic Sequence (domain 10000 ≤ N1 ≤ 99999):
- there are 920 integers iterating into the invariant 0 due to symmetries;
- the other ones iterate into one of the following loops (or a cyclic permutation of these): 21978, 65934; or 990, 8910, 6930, 2970, 4950; or 9009, 81081, 63063, 27027, 45045; or 9999, 89991, 69993, 29997, 49995.

The 6-Digit Periodic Sequence (domain 100000 ≤ N1 ≤ 999999):
- there are 13667 integers iterating into the invariant 0 due to symmetries;
- the longest sequence of iterations before arriving at the first loop member is 53 for N1 = 100720;
- the loops have 2, 5, 9, or 18 terms.

d) The Subtraction Periodic Sequences (IV)

Let c be a positive integer. Start with a positive integer N, and let N' be its digital reverse. Put N1 = abs(N1' - c), and let N1' be its digital reverse. Put N2 = abs (N1' - c), and let N2' be its digital reverse. And so on. We shall eventually obtain a repetition.
For example, with c = 1 and N = 52 we obtain the sequence: 52, 24, 41, 13, 30, 02, 19, 90, 08, 79, 96, 68, 85, 57, 74, 46, 63, 35, 52, ... . Here a repetition occurs after 18 steps, and the length of the repeating cycle is 18.

First example: c = 1, 10≤ N ≤999.
Every other member of this interval is an entry point into one of five cyclic periodic sequences (four of these are of length 18, and one of length 9).
When N is of the form 11k or 11k-1, then the iteration process results in 0.

Second example: 1 ≤ c ≤ 9, 100 ≤ N ≤ 999.
For c = 1, 2, or 5 all iterations result in the invariant 0 after, sometimes, a large number of iterations.
For the other values of c there are only eight different possible values for the length of the loops, namely 11, 22, 33, 50, 100, 167, 189, 200.
For c = 7 and N = 109 we have an example of the longest loop obtained: it has 200 elements, and the loop is closed after 286 iterations.
(H. Ibstedt)

e) The Multiplication Periodic Sequences (V)

Let c > 1 be a positive integer. Start with a positive integer N, multiply each digit x of N by c and replace that digit by the last digit of cx to give N1. And so on. We shall eventually obtain a repetition.
For example, with c = 7 and N = 68 we obtain the sequence:

68, 26, 42, 84, 68, ... .

Integers whose digits are all equal to 5 are invariant under the given operation after one iteration.

One studies the One-Digit Multiplication Periodic Sequences only.
(For c of two or more digits the problem becomes more complicated.)

If c = 2, there are four term loops, starting on the first or second term.

If c = 3, there are four term loops, starting with the first term.

If c = 4, there are two term loops, starting on the first or second term (could be called Switch or Pendulum).

If c = 5 or 6, the sequence is invariant after one iteration.

If c = 7, there are four term loops, starting with the first term.

If c = 8, there are four term loops, starting with the second term.

If c = 9, there are two term loops, starting with the first term (pendulum).
(H. Ibstedt)

f) The Mixed Composition Periodic Sequences (VI)

Let N be a two-digit number. Add the digits, and add them again if the sum is greater than 10. Also take the absolute value of their difference. These are the first and second digits of N1. Now repeat this.
For example, with N = 75 we obtain the sequence: 75, 32, 51, 64, 12, 31, 42, 62, 84, 34, 71, 86, 52, 73, 14, 53, 82, 16, 75, ... .
There are no invariants in this case. Four numbers: 36, 90, 93, and 99 produce two-element loops. The longest loops have 18 elements. There also are loops of 4, 6, and 12 elements.
(H. Ibstedt)

There will always be a periodic (invariant) sequence whenever we have a function $f : S \longrightarrow S$, where S is a finite set,
and we repeat the function f more times than card(S).
Thus the General Periodic Sequence is defined as:

a1 = f(s), where s is an element of S;
a2 = f(a1) = f(f(s));
a3 = f(a2) = f(f(a1)) = f(f(f(s)));
and so on.

245) New Sequences: The Family of Metallic Means

The family of Metallic Means (whom most prominent members are the Golden Mean, Silver Mean, Bronze Mean, Nickel Mean, Copper Mean, etc.) comprises every quadratic irrational number that is the positive solution of one of the algebraic equations

$$x^2 - nx - 1 = 0 \quad \text{or} \quad x^2 - x - n = 0,$$

where n is a natural number.
All of them are closely related to quasi-periodic dynamics, being therefore important basis of musical and architectural proportions. Through the analysis of their common mathematical properties, it becomes evident that they interconnect different human fields of knowledge, in the sense defined in "Paradoxist Mathematics".
Being irrational numbers, in applications to different scientific disciplines, they have to be approximated by ratios of integers -- which is the goal of this paper.
(Vera W. de Spinadel)

246-251) Some New Functions in the Number Theory.
Let S(n) be the Smarandache Function.
One defines:

S1 : N-{0,1} ⟶ N, S1(n) = 1/S(n);
and S2 : N* ⟶ N, S2(n) = S(n)/n
which verify the Lipschitz condition
Other functions::
S3 : N-{0,1} ⟶ N, S3(n) = n/S(n),
and Fs : N-{0,1} ⟶ N,

$$Fs(x) = \sum_{i=1}^{\pi(x)} ( S( p_i^x ),$$

where $p_i$ are the prime numbers not greater than x and

π(x) is the number of primes less than or equal to x;

Θ : N* ⟶ N,

$\Theta(x) = \sum S( p_i^x )$, where $p_i$ are prime numbers which divide x;

$\bar{\Theta}$ : N* ⟶ N,

$\bar{\Theta}(x) = \sum S( p_i^x )$, where $p_i$ are prime numbers not greater than x which do not divide x;

(V. Seleacu, S. Zanfir)

252) Erdös-Smarandache Numbers:
2, 3, 5, 6, 7, 10, 11, 13, 14, 15, 17, 19, 20, 21, 22, 23, 26, 28, 29, 30, 31, 33, 34, 35, ... .
Solutions to the Diophantine equation P(n)=S(n), where P(n) is the largest prime factor which divides n, and S(n) is the classical Smarandache function. [S. Tabarca]

253) Analogues of the Smarandache function:
a(n) is the smallest number m such that $n \leq m!$
1,1,2,3,3,3,3,4,4,4,4,4,4,4,4,4,4,4,4,4,4,4,4,4,5, … .

254) A Generalized Smarandache Palindrome (GSP) is a number of the concatenated form: $a_1a_2...a_na_n...a_2a_1$ with $n \geq 1$, or $a_1a_2...a_{n-1}a_na_{n-1}...a_2a_1$ with $n \geq 2$, where all $a_1, a_2, ..., a_n$ are positive integers of various number of digits.

Examples:

a) 1235656312 is a GSP because we can group it as (12)(3)(56)(56)(3)(12), i.e. ABCCBA.

b) 23523 is also a GSP since we can group it as (23)(5)(23), i.e. ABA.

It has been proven that 1234567891010987654321, which is a GSP, is a prime (see http://www.kottke.org/notes/0103.html, and the Prime Curios site).
Conjecture: There are infinitely many primes which are GSP.

255-260) Special Relationships (edited by M. Bencze).
Let $\{ a_n \}$, $n \geq 1$ be a sequence of numbers, and p, q integers $\geq 1$. Then we say that the

terms

$a_{k+1}, a_{k+2}, \ldots, a_{k+p}, a_{k+p+1}, a_{k+p+2}, \ldots, a_{k+p+q}$

verify a *special p-q relationship* if

$a_{k+1} <> a_{k+2} <> \ldots <> a_{k+p} = a_{k+p+1} <> a_{k+p+2} <> \ldots <> a_{k+p+q}$

where "<>" may be any arithmetic, algebraic, or analytic operation (generally a binary law on $\{ a_1 , a_2 , a_3 , \ldots \}$).

If this relationship is verified for any $k \geq 1$ (i.e. by all terms of the sequence), then

$\{ a_n \}$, n >= 1 is called a *Special p-q-<> sequence*

where "<>" is replaced by "additive" if <> is +, "multiplicative" if <> is *, etc. [according to the operation (<>) used].

As a particular case, we can easily see that Fibonacci/Lucas sequence ($a_n + a_{n+1} = a_{n+2}$ ), for $n \geq 1$, is a *Special 2-1 additive sequence.*

A Tribonacci sequence ( $a_n + a_{n+1} + a_{n+2} = a_{n+3}$ ), $n \geq 1$ is a Special 3-1 additive sequence. Etc.

M. Bencze considered the sequence of Smarandache numbers,

1, 2, 3, 4, 5, 3, 7, 4, 6, 5, 11, 4, 13, 7, 5, 6, 17, . . . ,

i.e. for each n the smallest number S(n) such that S(n)! is divisible by n (the values of the Smarandache Function), and raised the questions:

(a) How many quadruplets verify a Special 2-2 additive relationship i.e.

$S(n+1) + S(n+2) = S(n+3) + S(n+4)$?

He found that: $S(6) + S(7) = S(8) + S(9)$, $3 + 7 = 4 + 6$;

$S(7) + S(8) = S(9) + S(10)$, $7 + 4 = 6 + 5$;

$S(28) + S(29) = S(30) + S(31)$, $7 + 29 = 5 + 31$.

M. Bencze asked if there exist a finite or infinite number ?

(b) He also asked how many quadruplets verify a Special 2-2-subtractive relationship,

i.e.

$S(n+1) - S(n+2) = S(n+3) - S(n+4)$?

He found: $S(1) - S(2) = S(3) - S(4)$, $1 - 2 = 3 - 4$;

$S(2) - S(3) = S(4) - S(5)$, $2 - 3 = 4 - 5$;

$S(49) - S(50) = S(51) - S(52)$, $14 - 10 = 17 - 13$.

(c) M. Bencze also asked how many sextuplets verify a Special 3-3 additive relationship, i.e.

$S(n+1) + S(n+2) + S(n+3) = S(n+4) + S(n+5) + S(n+6)$ ?

He found: $S(5) + S(6) + S(7) = S(8) + S(9) + S(10)$, $5 + 3 + 7 = 4 + 6 + 5$.

Charles Ashbacher designed a computer program that calculates the Smarandache Function's values, therefore he may be able to add more solutions to mine.

261-264) More General Special Relationship (edited by M. Bencze):

If $f_p$ is a p-ary relation and $g_q$ is a q-ary relation, both of them defined on $\{ a_1, a_2, a_3, \ldots \}$, then

$a_{i1}, a_{i2}, \ldots, a_{ip}, a_{j1}, a_{j2}, \ldots, a_{jq}$

verify a *Special $f_p$ - $g_q$ - relationship* if

$f_p(a_{i1}, a_{i2}, \ldots, a_{ip}) = g_q(a_{j1}, a_{j2}, \ldots, a_{jq})$.

If this relationship is verified by all terms of the sequence, then $\{a_n\}$, $n \geq 1$ is called a *Special $f_p$ -$g_q$ -sequence*.

M. Bencze asked about the properties of Special $f_p$ -$g_q$ - relationships for well-known sequences (perfect numbers, Ulam numbers, abundant numbers, Catalan numbers, Cullen numbers, etc.).

For example: a Special 2-2-additive, subtractive, or multiplicative relationship, etc.

If $f_p$ is a p-ary relation on $\{ a_1, a_2, a_3, \ldots \}$ and

$f_p(a_{i1}, a_{i2}, \ldots, a_{ip}) = f_p(a_{j1}, a_{j2}, \ldots, a_{jp})$

for all $a_{ik}, a_{jk}$, where $k = 1, 2, \ldots, p$, and for all $p \geq 1$, then $\{a_n\}$, $n \geq 1$, is called a *Special perfect f- sequence*.

If not all p-plets $(a_{i1}, a_{i2}, \ldots, a_{ip})$ and $(a_{j1}, a_{j2}, \ldots, a_{jp})$ verify the $f_p$ relation, or not for all $p \geq 1$, the relation $f_p$ is verified, then $\{a_n\}$, $n \geq 1$ is called a *Special partial perfect f-sequence*.

An example: a *Special partial perfect additive sequence*:

1, 1, 0, 2, -1, 1, 1, 3, -2, 0, 0, 2, 1, 1, 3, 5, -4, -2, -1, 1, -1, 1, 1, 3, 0, 2, . . .

This sequence has the property that

$$a_1 + a_2 + \ldots + a_p = a_{p+1} + a_{p+2} + \ldots + a_{2p} \text{ for all } p \geq 1.$$

It is constructed in the following way:

$$a_1 = a_2 = 1$$

$$a_{2p+1} = a_{p+1} - 1$$

$$a_{2p+2} = a_{p+1} + 1$$

for all $p \geq 1$.

(a) M. Bencze asked for a general expression of $a_n$ (as function of n)?

It is periodical, or convergent, or bounded?

(b) He demanded for the design of other Special perfect (or partial perfect) sequences.

See a multiplicative sequence of this type.
265-268) P-digital Subsequences.
Let $\{a_n\}$, $n \geq 1$, be a sequence defined by a property (or a relationship involving its terms) P.
One screens this sequence, selecting only its terms whose digits hold the property (or relationship involving the digits) P.

For example:

(a) *Special square-digital subsequence:*

0, 1, 4, 9, 49, 100, 144, 400, 441, . . .

i.e. from 0, 1, 4, 9, 16, 25, 36, ..., $n^2$, ... one chooses only the terms whose digits are all perfect squares (therefore only 0, 1, 4, and 9).

Disregarding the square numbers of the concatenated form:

$$\overline{N0\ \ .\ \ .\ \ .\ \ 0},$$

{ 2k zeros }

where N is also a perfect square, M. Bencze asked how many other numbers belong to this sequence?

(b) *Special cube-digital subsequence:*

0, 1, 8, 1000, 8000, . . .

i.e. from 0, 1, 8, 27, 64, 125, 216, . . . , $n^3$, . . . one chooses only the terms whose digits are all perfect cubes (therefore only 0, 1 and 8).

Similar question, disregarding the cube numbers of the form

$$\overline{M0\ \ .\ \ .\ \ .\ \ 0},$$

{ 3k zeros }

where M is a perfect cube.

(c) *Special prime digital subsequence:*

2, 3, 5, 7, 23, 37, 53, 73, . . .

i.e. the prime numbers whose digits are all primes (they are called Smarandache-Wellin primes).

Conjecture: this sequence is infinite.

In the same general conditions of a given sequence, one screens it selecting only its terms whose groups of digits hold the property (or relationship involving the groups of digits) P.
[ A group of digits may contain one or more digits, but not the whole term.]

The new sequence obtained is called:
269-275) Special P-partial digital subsequence.

Similar examples:

(a) *Special square-partial-digital subsequence:*

49, 100, 144, 169, 361, 400, 441, . . .

i.e. the square members that is to be partitioned into groups of digits which are also perfect squares.

(169 can be partitioned as $16 = 4^2$ and $9 = 3^2$, etc.) Disregarding the square numbers of the form

$$\overline{N0\ .\ .\ .\ 0},$$
$$\{ \text{ 2k zeros } \}$$

where N is also a perfect square, how many other numbers belong to this sequence?

(b) *Special cube-partial digital subsequence:*

1000, 8000, 10648, 27000, . . .

i.e. the cube numbers that can be partitioned into groups of digits which are also perfect cubes.
(10648 can be partitioned as $1 = 1^3$, $0 = 0^3$, $64 = 4^3$, and $8 = 2^3$).

Same question: disregarding the cube numbers of the form:

$$\overline{M0\ .\ .\ .\ 0},$$
$$\{ \text{ 3k zeros } \}$$

where M is also a perfect cube, how many other numbers belong to this sequence?

(c) *Special prime-partial digital subsequence:*

23, 37, 53, 73, 113, 137, 173, 193, 197, . . .

i.e. prime numbers, that can be partitioned into groups of digits which are also prime,

(113 can be partitioned as 11 and 3, both primes).

Conjecture: this sequence is infinite.

(d) *Special Lucas-partial digital subsequence*

123, . . .

i.e. the sum of the two first groups of digits is equal to the last group of digits, and the whole number belongs to Lucas numbers:

2, 1, 3, 4, 7, 11, 18, 29, 47, 76, 123, 199, . . .

(beginning at 2 and $L(n+2) = L(n+1) + L(n)$, $n > 1$) ( 123 is partitioned as 1, 2 and 3, then $3 = 2 + 1$).

Is 123 the only Lucas number that verifies a Special type partition?

Study some Special P-(partial)-digital subsequences associated to:
- Fibonacci numbers (M. Bencze didn't find any Fibonacci number verifying a Special type partition; is there any?
- Smith numbers, Eulerian numbers, Bernoulli numbers, Mock theta numbers, Smarandache type sequences, etc.

M. Bencze remarked that some sequences may not be smarandachely partitioned (i.e. their associated Smarandache type subsequences might be empty).

If a sequence $\{a_n\}$, $n \geq 1$, is defined by $a_n = f(n)$ (a function of n), then *Special f-digital subsequence* is obtained by screening the sequence and selecting only its terms that can be partitioned in two groups of digits $g_1$ and $g_2$ such that

$g_2 = f(g_1)$.

For example:

(a) If $a_n = 2n$, $n \geq 1$, then

*Special even-digital subsequence is:*

12, 24, 36, 48, 510, 612, 714, 816, 918, 1020, 1122, 1224, . . .

(i.e. 714 can be partitioned as $g_1 = 7$, $g_2 = 14$, such that $14 = 2 \cdot 7$, etc. )

(b) *Special lucky-digital subsequence*

37, 49, . . .

(i.e. 37 can be partitioned as 3 and 7, and $L_3 = 7$; the lucky numbers are

1, 3, 7, 9, 13, 15, 21, 25, 31, 33, 37, 43, 49, 51, 63, . . .

276-284) Other Special Definitions and Conjectures (edited by M. Bencze).

For $n \geq 2$, let A be a set of $n^2$ elements, and *l* an n-ary law defined on A.
As a generalization of the XVIth - XVIIth centuries magic squares, one presents the *Special magic square of order n* which is a square array of rows of elements of A arranged so that the law *l* applied to each horizontal and vertical row and diagonal gives the same result.

If A is an arithmetic progression and *l* the addition of n numbers, then many magic squares have been found. Look at Durer's 15/4 engraving "Melancholia's" one:

| | | | |
|---|---|---|---|
| 16 | 3 | 2 | 13 |
| 5 | 10 | 11 | 8 |
| 9 | 6 | 7 | 12 |
| 4 | 15 | 14 | 1 |

(1) Can you find such a magic square of order at least 3 or 5, when A is a set of prime numbers and *l* the addition?

(2) Same question when A is a set of square numbers, or cube numbers, or special numbers [for example: Fibonacci or Lucas numbers, triangular numbers, Smarandache quotients (i.e. q(m) is the smallest k such that mk is a factorial), etc.].

A similar definition for the *Special magic cube of order n*, where the elements of A are arranged in the form of a cube of length n:

(a) either each element inside of a unitary cube (that the initial cube is divided in).

(b) either each element on a surface of a unitary cube.

(c) either each element on a vertex of a unitary cube.

(d) Study similar questions for this case, which is more complex.

An interesting law may be

$l(a_1, a_2, \ldots, a_n) = a_1 + a_2 - a_3 + a_4 - a_5 + \ldots$

Special prime conjecture:

Any odd number can be expressed as the sum of two primes minus a third prime (not including the trivial solution
p = p + q - q when the odd number is the prime itself).

For example:

1 = 3 + 5 - 7 = 5 + 7 - 11 = 7 + 11 - 17 = 11 + 13 - 24 = . . .

3 = 5 + 11 - 13 = 7 + 19 - 23 = 17 + 23 - 37 = . . .

5 = 3 + 13 - 11 = . . .

7 = 11 + 13 - 17 = . . .

9 = 5 + 7 - 3 = . . .

11 = 7 + 17 - 13 = . . .

(a) Is this conjecture equivalent to Goldbach's conjecture (any odd number > = 9 can be expressed as a sum of three primes - finally solved 15 by Vinogradov in 1937 for any odd number greater than $_{3}3^{15}$ )?

(b) The number of times each odd number can be expressed as a sum of two primes minus a third prime are called *Smarandache prime conjecture numbers*. None of them are known!

(c) Write a computer program to check this conjecture for as many positive odd numbers as possible.

(e) There are infinitely many numbers that cannot be expressed as the difference between a cube and a square (in absolute value).

They are called *bad numbers* (!)

For example: 5, 6, 7, 10, 13, 14, . . . are probably such bad numbers (F. Smarandache has conjectured that), while

1, 2, 3, 4, 8, 9, 11, 12, 15, . . . are not, because

$1 = | 2^3 - 3^2 |$

$2 = | 3^3 - 5^2 |$

$3 = | 1^3 - 2^2 |$

$4 = | 5^3 - 11^2 |$

$8 = | 1^3 - 3^2 |$

$9 = | 6^3 - 15^2 |$

$11 = | 3^3 - 4^2 |$

$12 = | 13^3 - 47^2 |$

$15 = | 4^3 - 7^2 |$, etc.

(a) Write a computer program to get as many non Smarandache bad numbers (it's easier this way!) as possible,

i.e. find an ordered array of a's such that

$a = | x^3 - y^2 |$, for x and y integers $\geq 1$.

285) Anti-Prime Function is defined as follows:

$P : \mathrm{N} \mapsto \{0, 1\}$, with

$$P(n) = \begin{cases} 0, & \text{if } n \text{ prime}; \\ 1, & \text{otherwise}. \end{cases}$$

For example $P(2) = P(3) = P(5) = P(7) = P(11) = \ldots = 0$, whereas $P(0) = P(1) = P(4) = P(6) = \ldots = 1$.

More general:
$P_k : N^k \mapsto \{0, 1\}$, where k is an integer $\geq 2$, and

$$P_k(n_1, n_2, \ldots, n_k) = \begin{cases} 0, & \text{if all } n_1, n_2, \ldots, n_k \text{ are prime numbers;} \\ 1, & \text{otherwise.} \end{cases}$$

286) Anti-Coprime Function is similarly defined:
$C_k : N^k \mapsto \{0, 1\}$, where k is an integer $\geq 2$, and

$$C_k(n_1, n_2, \ldots, n_k) = \begin{cases} 0, & \text{if all } n_1, n_2, \ldots, n_k \text{ are coprime numbers;} \\ 1, & \text{otherwise.} \end{cases}$$

287) Special Functional Iteration of First Kind:

Let $f: A \mapsto A$ be a function, such that $f(x) \leq x$ for all x, and

$$\min \{f(x), x \in A\} \geq m_0 \neq -\infty.$$

Let f have $p \geq 1$ fix points: $m_0 \leq x_1 < x_2 < \ldots < x_p$.
[The point x is called fix if $f(x) = x$.]
Then

$SI1_f(x)$ = the smallest number of iterations k such that
f(f(...f(x)...)) = constant
{f iterated k times}.

Example:
Let $n > 1$ be an integer, and d(n) be the number of positive divisors of n, $d: N \mapsto N$.
Then $SI1_d(n)$ is the smallest number of iterations k such that $d(d(...d(n)...)) = 2$,
{d iterated k times }
because $d(n) < n$ for $n > 2$, and the fix points of the function d are 1 and 2.
Thus $SI1_d(6) = 3$, because $d(d(d(6))) = d(d(4)) = d(3) = 2$ = constant.
$SI1_d(5) = 1$, because $d(5) = 2$.

288) Special Functional Iteration of Second Kind:

Let $g: A \mapsto A$ be a function, such that $g(x) > x$ for all x, and let $b > x$. Then:
$SI2_g(x, b)$ = the smallest number of iterations k such that
$g(g(...g(x)...)) \geq b$
{ g iterated k times }.

Example:
Let $n > 1$ be an integer, and $\sigma(n)$ be the sum of positive divisors of n (1 and n included), $\sigma : N \mapsto N$.
Then $SI2_\sigma(n, b)$ is the smallest number of iterations k such that
$\sigma(\sigma(...\sigma(n)...)) \geq b$
{ σ iterated k times }
because $\sigma(n) > n$ for $n > 1$.
Thus $SI2_\sigma(4, 11) = 3$, because $\sigma(\sigma(\sigma(4))) = \sigma(\sigma(7)) = \sigma(8) = 15 \geq 11$.

289) Special Functional Iteration of Third Kind:
Let $h: A \mapsto A$ be a function, such that $h(x) < x$ for all x, and let $b < x$. Then:
$SI3_h(x, b)$ = the smallest number of iterations k such that
$h(h(...h(x)...)) \leq b$
{ h iterated k times }.

Example:
Let n be an integer and gd(n) be the greatest divisor of n, less than n,

gd: N* ↣ N*. Then gd(n) < n for n > 1.
$SI3_{gd}(60, 3) = 4$, because gd(gd(gd(gd(60)))) = gd(gd(gd(30))) =
gd(gd(15)) = gd(5) = 1 ≤ 3.

290) Smarandacheials

Let n > k ≥ 1 be two integers. Then the Smarandacheial is defined as:

$$!n!_k = \prod_{\substack{0<|n-k\cdot i|\le n \\ i\in N}} (n-k\cdot i)$$

By example:

In the case k=3:

$$!n!_3 = \prod_{\substack{0<|n-3i|\le n \\ i=0,\,1,\,2,\,\dots}} (n-3i) = n(n-3)(n-6)\dots .$$

Thus for n=7 one has $!7!_3$ = 7(7-3)(7-6)(7-9)(7-12) = 7(4)(1)(-2)(-5) = 280.

More general:

Let n > k ≥ 1 be two integers and m ≥ 1 another integer. Then the generalized Smarandacheial is defined as:

$$!n!m_k = \prod_{\substack{0<|n-k\cdot i|\le m \\ i\in N}} (n-k\cdot i)$$

For example:

$!7!9_2$ = 7(7-2)(7-4)(7-6)(7-8)(7-10)(7-12)(7-14)(7-16) = 7(5)(3)(1)(-1)(-3)(-5)(-7)(-9)
= -99225.

www.gallup.unm.edu/~smarandache/Smarandacheials.htm.

F. Smarandache, "Back and Forth k-Factorials", Arizona State Univ., Special Collections, 1972.

291) Smarandache-Kurepa Function:

For p prime, SK(p) is the smallest integer such that !SK(p) is divisible by p, where !SK(p) = 0! + 1! + 2! + ... + (p-1)!

For example:

| p | 2 | 3 | 7 | 11 | 17 | 19 | 23 | 31 | 37 | 41 | 61 | 71 | 73 | 89 |
|---|---|---|---|---|---|---|---|---|---|---|---|---|---|---|
| SK(p) | 2 | 4 | 6 | 6 | 5 | 7 | 7 | 12 | 22 | 16 | 55 | 54 | 42 | 24 |

292) Smarandache-Wagstaff Function:

For p prime, SW(p) is the smallest integer such that W(SW(p)) is divisible by p, where W(p) = 1! + 2! + ... + (p)!

For example:

| p | 3 | 11 | 17 | 23 | 29 | 37 | 41 | 43 | 53 | 67 | 73 | 79 | 97 |
|---|---|---|---|---|---|---|---|---|---|---|---|---|---|
| SW(p) | 2 | 4 | 5 | 12 | 19 | 24 | 32 | 19 | 20 | 20 | 7 | 57 | 6 |

293) Smarandache Ceil Function of Second Order:

($S_2$ (n) = m, where m is the smallest positive integer for which n divides $m^2$.)

2, 4, 3, 6, 4, 6, 10, 12, 5, 9, 14, 8, 6, 20, 22, 15, 12, 7, 10, 26, 18, 28, 30, 21, 8, 34, 12, 15, 38, 20, 9, 42, 44, 30, 46, 24, 14, 33, 10, 52, 18, 28, 58, 39, 60, 11, 62, 25, 42, 16, 66, 45, 68, 70, 12, 21, 74, 30, 76, 51, 78, 40, 18, 82, 84, 13, 57, 86, ...

294) Smarandache Ceil Function of Third Order:
($S_3$ (n) = m, where m is the smallest positive integer for which n divides $m^3$.)
2, 2, 3, 6, 4, 6, 10, 6, 5, 3, 14, 4, 6, 10, 22, 15, 12, 7, 10, 26, 6, 14, 30, 21, 4, 34, 6, 15, 38, 20, 9, 42, 22, 30, 46, 12, 14, 33, 10, 26, 6, 28, 58, 39, 30, 11, 62, 5, 42, 8, 66, 15, 34, 70, 12, 21, 74, 30, 38, 51, 78, 20, 18, 82, 42, 13, 57, 86, ...

295) Smarandache Ceil Function of Fourth Order:
($S_4$ (n) = m, where m is the smallest positive integer for which n divides $m^4$.)
2, 2, 3, 6, 2, 6, 10, 6, 5, 3, 14, 4, 6, 10, 22, 15, 6, 7, 10, 26, 6, 14, 30, 21, 4, 34, 6, 15, 38, 10, 3, 42, 22, 30, 46, 12, 14, 33, 10, 26, 6, 14, 58, 39, 30, 11, 62, 5, 42, 4, 66, 15, 34, 70, 6, 21, 74, 30, 38, 51, 78, 20, 6, 82, 42, 13, 57, 86, ...

296) Smarandache Ceil Function of Fifth Order:
($S_5$ (n) =m, where m is the smallest positive integer for which n divides $m^5$.)
2, 2, 3, 6, 2, 6, 10, 6, 5, 3, 14, 2, 6, 10, 22, 15, 6, 7, 10, 26, 6, 14, 30, 21, 4, 34, 6, 15, 38, 10, 3, 42, 22, 30, 46, 6, 14, 33, 10, 26, 6, 14, 58, 39, 30, 11, 62, 5, 42, 4, 66, 15, 34, 70, 6, 21, 74, 30, 38, 51, 78, 10, 6, 82, 42, 13, 57, 86, ...

297) Smarandache Ceil Function of Sixth Order:
($S_6$ (n) = m, where m is the smallest positive integer for which n divides $m^6$.)
2, 2, 3, 6, 2, 6, 10, 6, 5, 3, 14, 2, 6, 10, 22, 15, 6, 7, 10, 26, 6, 14, 30, 21, 2, 34, 6, 15, 38, 10, 3, 42, 22, 30, 46, 6, 14, 33, 10, 26, 6, 14, 58, 39, 30, 11, 62, 5, 42, 4, 66, 15, 34, 70, 6, 21, 74, 30, 38, 51, 78, 10, 6, 82, 42, 13, 57, 86, ...

298) Smarandache Ceil Functions of k-th Order:
Sk(n) is the smallest integer for which n divides Sk(n)^k.

299) Pseudo-Smarandache Function:
Z(n) is the smallest integer such that 1 + 2 + ... + Z(n) is divisible by n.

For example:

| n | 1 | 2 | 3 | 4 | 5 | 6 | 7 |
|---|---|---|---|---|---|---|---|
| Z(n) | 1 | 3 | 2 | 3 | 4 | 3 | 6 |

300) Smarandache Near-To-Primordial Function:
SNTP(n) is the smallest prime such that either p* - 1, p* , or p* + 1 is divisible by n,
where p* , of a prime number p, is the product of all primes less than or equal to p.

For example:

| n | 1 | 2 | 3 | 4 | 5 | 6 | 7 | 8 | 9 | 10 | 11 | ... | 59 ... |
|---|---|---|---|---|---|---|---|---|---|---|---|---|---|
| SNTP(n) | 2 | 2 | 2 | 5 | 3 | 3 | 3 | 5 | ? | 5 | 11 | ... | 13 ... |

301) Smarandache - Fibonacci triplets:
(An integer n such that S(n) = S(n-1) + S(n-2) where S(k) is the Smarandache function: the smallest number k such that S(k)! is divisible by k.)
11, 121, 4902, 26245, 32112, 64010, 368140, 415664, 2091206, 2519648, 4573053, 7783364, 79269727, 136193976, 321022289, 445810543, 559199345, 670994143, 836250239, 893950202, 937203749, 1041478032, 1148788154, ...

Remarks:
It is not known if this sequence has infinitely or finitely many terms.
H. Ibstedt and C. Ashbacher independently conjectured that there are infinitely many.
H. I. found the biggest known number: 19 448 047 080 036.

302) Smarandache-Radu duplets
(An integer n such that between S(n) and S(n+1) there is no prime [S(n) and S(n + 1) included], where S(k) is the Smarandache function: the smallest number k such that S(k)! is divisible by k.)
224, 2057, 265225, 843637, 6530355, 24652435, 35558770, 40201975, 45388758, 46297822, 67697937, 138852445, 157906534, 171531580, 299441785, 551787925, 1223918824, 1276553470, 1655870629, 1853717287, 1994004499, 2256222280, ...

Remarks:
It is not known if this sequence has infinitely or finitely many terms.
H. Ibstedt conjectured that there are infinitely many.
H. I. found the biggest known number:
270 329 975 921 205 253 634 707 051 822 848 570 391 313 !

303) Concatenated Fibonacci Sequence:
1, 11, 112, 1123, 11235, 112358, . . . .

How many primes does it contain?

304) Uniform Sequences:

Let n be an integer not equal to zero and $d_1$ , $d_2$ , . . . , $d_r$ digits in a base B (of course $r < B$). Then: multiples of n, written with digits $d_1$ , $d_2$ , . . . , $d_r$ only (but all r of them), in base B, increasingly ordered, are called uniform sequences.

Let's see some examples in base 10 for one digit.

(a) Multiples of 7 written with digit 1 only:

111111, 111111,111111, 111111,111111,111111, 111111,111111,111111,111111, ...

(b) Multiples of 7 written with digit 2 only:

222222, 222222222222, 222222222222222222, 222222222222222222222222, ...

(c) Multiples of 79365 written with digit 5 only:

555555, 555555555555, 555555555555555555, 555555555555555555555555, ...

There are cases where the uniform sequence may be empty (impossible):

(d) Multiples of 79365 written with digit 6 only (because any multiple of 79365 will end in 0 or 5.

As a consequence, if there exists at least a multiple m of n, written with digits $d_1$ ,$d_2$ , . . ., $d_r$ only, in base B, then there exists an infinite number of multiples of n (they have the form:

$m, \overline{mm}, \overline{mmm}, \overline{mmmm}, \dots$

305-311) Special Operation Sequences (edited by S. Smith):

Let E be an ordered set of elements, E = { $e_1$ , $e_2$ , . . . } and O a set of binary operations well-defined for these elements.

Then: $a_1$ is an element of { $e_1$ ,$e_2$ , . . . }.

$a_{n+1} = \min \{ e_1 O_1 e_2 O_2 \dots O_n e_{n+1} \} > a_n$ , for $n > 1$.

where all $O_i$ are operations belonging to O, is called the Special operation sequence.

Some examples:

(a) When E is the natural number set, and O is formed by the four arithmetic operations: +, -, *, /.

Then: $a_1 = 1$
$a_{n+1} = \min \{ 1\ O_1\ 2\ O_2 \ldots O_n\ (n+1) \} > a_n$ , for $n > 1$,
(therefore, all $O_i$ may be chosen among addition, subtraction, multiplication or division in a convenient way).

Questions: Find this *arithmetic operation infinite sequence*. Is it possible to get a general expression formula for this sequence (which starts with 1, 2, 3, 5, 4,?

(b) A finite sequence
$a_1 = 1$
$a_{n+1} = \min \{ 1\ O_1\ 2\ O_2 \ldots O_{98}\ 99 \} > a_n$
for $n > 1$, where all $O_i$ are elements of { +, -, *, / }.
Same questions for this *arithmetic operation finite sequence*.

(c) Similarly for *algebraic operation infinite sequence*
$a_1 = 1$
$a_{n+1} = \min \{ 1\ O_1\ 2\ O_2 \ldots O_n\ (n+1) \} > a_n$ for $n > 1$,
where all $O_i$ are elements of { +, -, *, /, **, ythrtx }
( X**Y means $X^Y$, and "ythrtx" means the y-th root of x).

The same questions become harder but more exciting.

(d) Similarly for algebraic operation finite sequence:
$a_1 = 1$
$a_{n+1} = \min \{ 1\ O_1\ 2\ O_2 \ldots O_{98}\ 99\} > a_n$ , for $n > 1$,
where all $O_i$ are elements of { +, -, *, /, **, ythrtx }
( X**Y means $X^Y$ and "ythrtx" means the y-th root of x).

Same questions.
More generally: one replaces "binary operations" by "$K_i$ -ary operations"
where all $K_i$ are integers $\geq 2$).
Therefore, $a_i$ is an element of $\{ e_1 , e_2 , \ldots \}$,
$a_{n+1} = \min\{ 1\ O_1^{(K_1)} 2 O_1^{(K_1)} \ldots O_1^{(K_1)}\ K_1\ O_2^{(K_2)}(K_2+1 O_2^{(K_2)} \ldots O_2^{(K_2)} (K_1+K_2-1) \ldots (n+2-K_r) O_r^{(K_r)} \ldots O_r^{(K_r)}(n+1)\} > a_n$, for $n > 1$,
where $O_1^{(K_1)}$ is $K_1$ – ary, $O_2^{(K_2)}$ is $K_2$ – ary,
and of course $K_1 + (K_2 - 1) + \ldots + (K_r - 1) = n+1$.

Remark: The questions are much easier when O = { +, -}; study the operation type sequences in this case.

(9) *Operation sequences at random*:
Same definitions and questions as for the previous sequences, except that
$a_{n+1} = \{ e_1\ O_1\ e_2\ O_2 \ldots O_n\ e_{n+1} \} > a_n$ , for $n > 1$,
(i.e. it's no "min" any more, therefore $a_{n+1}$ will be chosen at random, but greater than $a_n$ , for any n > 1). Study these sequences with a computer program for random variables (under weak conditions).

312) A function $f: N^* \mapsto N^*$ is called S-multiplicative

if (a, b) = 1 implies that $f(a \cdot b) = \max\{ f(a), f(b) \}$.

(S. Tabarca)

313) General Theorem of Characterization of N Primes Simultaneously

Let $p_{ij}$, $1 \le i \le n$, $1 \le j \le m_i$, be coprime integers two by two,

and let $a_i$, $r_i$ be integers such that $a_i$ and $r_i$ are coprime for all i.

The following conditions are considered for all i:

(i) $p_{i1}, \ldots, p_{im_i}$ are simultaneously prime iff $c_i \cong 0 \pmod{r_i}$.

Then the following statements are equivalent:

a) The numbers $p_{ij}$ , $1 \le i \le n$, $1 \le j \le m_i$, are simultaneously prime.

b) $(R/D) \sum_{i=1}^{n} (a_i c_i / r_i) \cong 0 \pmod{R/D}$, where $R = \prod_{i=1}^{n} r_i$ and D is a divisor of R.

c) $(P/D) \sum_{i=1}^{n} (a_i c_i / \prod_{j=1}^{mi} p_{ij}) \cong 0 \pmod{P/D}$, where $P = \prod_{,j=1}^{n,mi} r_i$ and D is a divisor of P.

d) $\sum_{i=1}^{n} (a_i c_i (P / \prod_{j=1}^{mi} p_{ij}) \cong 0 \pmod{P}$, where $P = \prod_{,j=1}^{n,mi} r_i$.

e) $\sum_{i=1}^{n} (a_i c_i / \prod_{j=1}^{mi} p_{ij})$ is an integer.

314) Infinite Product.
It is defined as:

$$\prod_{n>=1} 1/a(n)$$

where a(n) is any of the above sequences, subsequences, or functions, or any other infinite product involving such sequences, subsequences, or functions. Some of them will lead to nice constants.

315) Jung's Theorem Generalized to Space:
Let's have n points in the 3D-space such that the maximum distance between any two of them is d. Then: there exists a sphere of radius

$$r \le d\frac{\sqrt{6}}{4}$$

which contains in interior or on surface all these points.

316-325) Back and Forth Summants

Let $n>k\geq 1$ be two integers. Then a Back and Forth Summant is defined as:

$$S(n, {}_{k}) = \sum_{\substack{0<|n-k\cdot i|\leq n \\ i=0,1,2,\ldots}} (n-k\cdot i) \quad \text{[for signed numbers]}$$

$$S|n, {}_{k}| = \sum_{\substack{0<|n-k\cdot i|\leq n \\ i=0,1,2,\ldots}} |n-k\cdot i| \quad \text{[for absolute value numbers]}$$

which are duals and semi-duals respectively of Smarandacheials.

$S(n, {}_1)$ and $S(n, {}_2)$ with corresponding $S|n, {}_1|$ and $S|n, {}_2|$ are trivial.

a) In the case k=3:

$$S(n, {}_{3}) = \sum_{\substack{0<|n-3i|\leq n \\ i=0,1,2,\ldots}} (n-3i) = n+(n-3)+(n-6)+\ldots \text{ ; [for signed numbers].}$$

$$S|n, {}_{3}| = \sum_{\substack{0<|n-3i|\leq n \\ i=0,1,2,\ldots}} |n-3i| = n+|n-3|+|n-6|+\ldots \text{ ; [for absolute value numbers].}$$

Thus $S(7, {}_3) = 7+(7-3)+(7-6)+(7-9)+(7-12) = 7+(4)+(1)+(-2)+(-5) = 5$; [for signed numbers].
Thus $S|7, {}_3| = 7+|7-3|+|7-6|+|7-9|+|7-12| = 7+4+1+2+5 = 19$; [for absolute value numbers].

The sequence is $S(n, {}_3)$: 3, 2, 0, 5, 3, 0, 7, 4, 0, 9, 5, 0, … ; [for signed numbers].
The sequence is $S|n, {}_3|$: 7, 12, 18, 19, 27, 36, 37, 48, **… ;** [for absolute value numbers].

b) In the case k=4:

$$S(n, {}_{4}) = \sum_{\substack{0<|n-4i|\leq n \\ i=0,1,2,\ldots}} (n-4i) = n+(n-4)+(n-8)\ldots \text{ ; [for signed numbers].}$$

$$S|n, {}_{4}| = \sum_{\substack{0<|n-4i|\leq n \\ i=0,1,2,\ldots}} |n-4i| = n+|n-4|+|n-8|\ldots \text{ ; [for absolute value numbers].}$$

Thus $S(9, {}_4) = 9+(9-4)+(9-8)+(9-12)+(9-16) = 9+(5)+(1)+(-3)+(-7) = 5$; for signed

numbers.
Thus $S|9, {}_4| = 9+|9-4|+|9-8|+|9-12|+|9-16| = 9+5+1+3+7 = 25$; [for absolute value numbers].

The sequence is $S(n, {}_4) =$ 3, 0, 4, 0, 5, 0, 6, 0, 7, 0, 8, 0, 9, 0, 10, 0, 11, … .
The sequence is $S|n, {}_4|$ = 9, 16, 16, 24, 25, 36, 36, 48, 49, 64, 64, 80, 81, 100, 100, **… .**

c) In the case k=5:

$$S(n, {}_5) = \sum_{\substack{0<|n-5i|\leq n \\ i=0,1,2,\ldots}} (n-5i) = n+(n-5)+(n-10)\ldots .$$

$$S|n, {}_5| = \sum_{\substack{0<|n-5i|\leq n \\ i=0,1,2,\ldots}} |n-5i| = n+|n-5|+|n-10|\ldots .$$

Thus $S(11, {}_5) = 11+(11-5)+(11-10)+(11-15)+(11-20) = 11+6+1+(-4)+(-9) = 5$.
Thus $S|11, {}_5| = 11+|11-5|+|11-10|+|11-15|+|11-20| = 11+6+1+4+9 = 31$.

The sequence is $S(n, {}_5)$: 3, 6, 2, 6, 0, 5, 10, 3, 9, 0, 7, 14, 4, 12, 0, … .
The sequence is $S|n, {}_5|$: 11, 12, 20, 20, 30, 31, 32, 33, 45, 60, 61, 62, 80, 80, 100, **… .**

More general:
Let $n>k\geq 1$ be two integers and $m\geq 0$ another integer.
Then the Generalized Back and Forth Summant is defined as:

$$\mathbf{S(}n, m, {}_k) = \sum_{\mathbf{i}=0, 1, 2, \ldots, \text{floor}[(\mathbf{n+m})/k]} (n-\mathbf{k \cdot i}) \quad \text{[for signed numbers].}$$

$$\mathbf{S|n}, m, {}_k| = \sum_{\mathbf{i}=0, 1, 2, \ldots, \text{floor}[(\mathbf{n+m})/k]} |n-\mathbf{k \cdot i}| \quad \text{[for absolute value numbers].}$$

For examples:
$S(7, 9, {}_2) = 7+(7-2)+(7-4)+(7-6)+(7-8)+(7-10)+(7-12)+(7-14)+(7-16)$
$= 7+(5)+(3)+(1)+(-1)+(-3)+(-5)+(-7)+(-9) = -2$.
$S|7, 3, {}_2| = 7+|7-2|+|7-4|+|7-6|+|7-8|+|7-10| = 7+5+3+1+1+3 = 20$.

CONTENTS

42) (Inferior) Cube Part
43) (Superior) Cube Part
44) (Inferior) Factorial Part
45) (Superior) Factorial Part
46) Inferior Fractional f-Part of x
47) Superior Fractional f-Part of x
48) Fractional Prime Part
49) Fractional Square Part
50) Fractional Cubic Part
51) Fractional Factorial Part
52) Smarandacheian Complements
53) Square Complements
54) Cubic Complements
55) m-Power Complements
56) Double Factorial Complements
57) Prime Additive Complements
58) Prime Base
59) Square Base
60) m-Power Base (generalization)
61) Factorial Base
62) Double Factorial Base
63) Triangular Base
64) Generalized Base
65) Smarandache Numbers
66) Smarandache Quotients
67) Double Factorial Numbers
68) Primitive Numbers (of power 2)
69) Primitive Numbers (of power 3)
70) Primitive Numbers (of power p, p prime) {generalization}
71) Smarandache Functions of the First Kind
72) Smarandache Functions of the Second Kind
73) Smarandache Function of the Third Kind
74) Pseudo-Smarandache Numbers
75) Square Residues
76) Cubical Residues
77) m-Power Residues (generalization)
78) Exponents (of power 2)
79) Exponents (of power 3)
80) Exponents (of power p) {generalization}
81) Pseudo-Primes of First Kind
82) Pseudo-Primes of Second Kind
83) Pseudo-Primes of Third Kind
84) Almost Primes of First Kind
85) Almost Primes of Second Kind

86) Pseudo-Squares of First Kind
87) Pseudo-Squares of Second Kind
88) Pseudo-Squares of Third Kind
89) Pseudo-Cubes of First Kind
90) Pseudo-Cubes of Second Kind
91) Pseudo-Cubes of Third Kind
92) Pseudo-m-Powers of First Kind
93) Pseudo-m-powers of second kind
94) Pseudo-m-Powers of Third Kind
95) Pseudo-Factorials of First Kind
96) Pseudo-Factorials of Second Kind
97) Pseudo-Factorials of Third Kind
98) Pseudo-Divisors of First Kind
99) Pseudo-Divisors of Second Kind
100) Pseudo-Divisors of Third Kind
101) Pseudo-Odd Numbers of First Kind
102) Pseudo-odd Numbers of Second Kind
103) Pseudo-Odd Numbers of Third Kind
104) Pseudo-Triangular Numbers
105) Pseudo-Even Numbers of First Kind
106) Pseudo-Even Numbers of Second Kind
107) Pseudo-Even Numbers of Third Kind
108) Pseudo-Multiples of First Kind (of 5)
109) Pseudo-Multiples of Second Kind (of 5)
110) Pseudo-Multiples of Third kind (of 5)
111) Pseudo-Multiples of first kind of p (p is an integer $\geq$2) {Generalizations}
112) Pseudo-Multiples of Second Kind of p (p is an integer $\geq$ 2)
113) Pseudo-multiples of third kind of p (p is an integer $\geq$ 2)
114) Constructive Set (of digits 1, 2)
115) Constructive Set (of digits 1,2,3)
116) Generalized Constructive Set
117) Square Roots
118) Cubical Roots
119) m-Power Roots
120) Numerical Carpet
121) Goldbach-Smarandache Table
122) Smarandache-Vinogradov Table
123) Smarandache-Vinogradov Sequence
124) Smarandache Paradoxist Numbers
125) Non-Smarandache Numbers
126) The Paradox of Smarandache numbers
127) Romanian Multiplication
128) Division by $k^n$
129) Generalized Period

130) Number of Generalized Periods
131) Length of Generalized Period
132) Sequence of Position
133) Criterion for Coprimes
134) Congruence Function
135) Generalization of Euler's Theorem
136) Smarandache Simple functions
137) Smarandache Function
138) Prime Equation Conjecture
139) Generalized Prime Equation Conjecture
140) Progressions
141) Inequality
142) Divisibility Theorem
143) Dilemmas
144) Surface Points
145) Inclusion Problems
146) Convex Polyhedrons
147) Integral Points
148) Counter
149) Maximum Points
150) Minimum Points
151 Increasing Repeated Compositions
152) Decreasing Repeated Compositions
153) Coloration Conjecture
154) Primes
155) Prime Number Theorem
156) Cardinality Theorem
157) Concatenate Natural Sequence
159) Back Concatenated Prime Sequence
160) Concatenated Odd Sequence
161) Back Concatenated Odd Sequence
162) Concatenated Even Sequence
163) Back Concatenated Even Sequence
164) Concatenated S-Sequence {generalization}
165) Concatenated Square Sequence
166) Back Concatenated Square Sequence
167) Concatenated Cubic Sequence
168) Back Concatenated Cubic Sequence
169) Concatenated Fibonacci Sequence
170) Back Concatenated Fibonacci Sequence
171) Power Function
172) Reverse Sequence
173) Multiplicative Sequence
174) Wrong Numbers

175) Impotent Numbers
176) Random Sieve
177) Non-Multiplicative Sequence
178) Non-Arithmetic Progression
179) Prime Product Sequence
180) Square Product Sequence
181) Cubic Product Sequence
182) Factorial Product Sequence
183) U-Product Sequence {generalization}
184-190) Sequences of Sub-Sequences
191) $S^2$ function (numbers)
192) $S^3$ function (numbers)
193) Anti-Symmetric Sequence
194-202) Recurrence Type Sequences
203- 215) Partition Type Sequences
216) Reverse Sequence
217) Non-Geometric Progression
218) Unary Sequence
219) No-Prime-Digit Sequence
220) No-square-digits Sequence
221-222) Polygons and Polyhedrons
223) Lucky Numbers
224) Lucky Method/Algorithm/Operation/etc.
225) Sequence of 3's and 1's
226) Sequence of 3's and 1's
227) General Sequence of Concatenation of Given Digits
228) Recurrence of Smarandache type Functions
229-231) G Add-On Sequence (I)
232) Non-Arithmetic Progression (I)
233) Concatenation Type Sequences
234) Construction of Elements of the Square-Partial-Digital Subsequence
235) Prime-Digital Sub-Sequence
236-237) Special Expressions
238) General Periodic Sequence
239-244) Periodic Digit Sequences
245) New Sequences The Family of Metallic Means
246-251) Other New Functions in the Number Theory
252) Erdös-Smarandache Numbers
253) Analogues of the Smarandache Function
254) Generalized Smarandache Palindrome (GSP)
255-260) Special Relationships
261-264) More General Special Relationship
265-268) P-digital Subsequences
269-275) Special P-Partial Digital Subsequence

Over 300 sequences and many unsolved problems and conjectures related to them are presented herein.

The book contains definitions, unsolved problems, questions, theorems corollaries, formulae, conjectures, examples, mathematical criteria, etc. ( on integer sequences, numbers, quotients, residues, exponents, sieves, pseudo-primes/squares/cubes/factorials, almost primes, mobile periodicals, functions, tables, prime/square/factorial bases, generalized factorials, generalized palindromes, etc. ).